\documentclass{article}

\usepackage{amssymb} 
\usepackage{amsmath} 
\usepackage{latexsym} 
\usepackage{theorem} 

\usepackage{epsfig}

\def\vect#1{\mbox{\boldmath $#1$}} 
\def\tanufrac#1#2{\mbox{\footnotesize$\displaystyle\frac{\,#1\,}{\,#2\,}$}}
\def\e{\varepsilon}
\def\p{\prime}

\theorembodyfont{\itshape} 

\newtheorem{theorem}{Theorem}[section]
\newtheorem{lemma}[theorem]{Lemma}

\newtheorem{corollary}[theorem]{Corollary}
\newtheorem{proposition}[theorem]{Proposition}

\newtheorem{definition}[theorem]{Definition}

\newtheorem{problem}[theorem]{Problem}
\newenvironment{remark}{{\noindent \bf Remark:}}{}
\newenvironment{proof}{{\par\addvspace{0.1cm}\noindent \bf P{\footnotesize ROOF. }}}{\hfill$\Box$} 
\newenvironment{proofoftheorem}{{\par\addvspace{0.1cm}\noindent \bf P{\footnotesize ROOF OF} T{\footnotesize HEOREM}\hspace{0.1cm}}}{\hfill$\Box$} 
\newenvironment{proofofproposition}{{\par\addvspace{0.1cm}\noindent \bf P{\footnotesize ROOF OF} P{\footnotesize ROPOSITION}\hspace{0.1cm}}}{\hfill$\Box$}

\title{\bf The configuration space of planar spidery linkages}
\author{Jun O'Hara}

\begin{document}

\maketitle

Dedicated to Professor Yukio Matsumoto on his 60th birthday

\begin{abstract} The configuration space of the mechanism of a planar robot is studied. 
We consider a robot which has $n$ arms such that each arm is of length $1+1$ and has a rotational joint in the middle, and that the endpoint of the $k$-th arm is fixed to $Re^{\frac{2(k-1)\pi}ni}$. 
Generically, the configuration space is diffeomorphic to an orientable closed surface. 
Its genus is given by a topological way and a Morse theoretical way. 
The homeomorphism types of it when it is singular is also given. 
\end{abstract}



\section{Introduction}

\begin{figure}[htbp]
\begin{center}
\includegraphics[width=.65\linewidth]{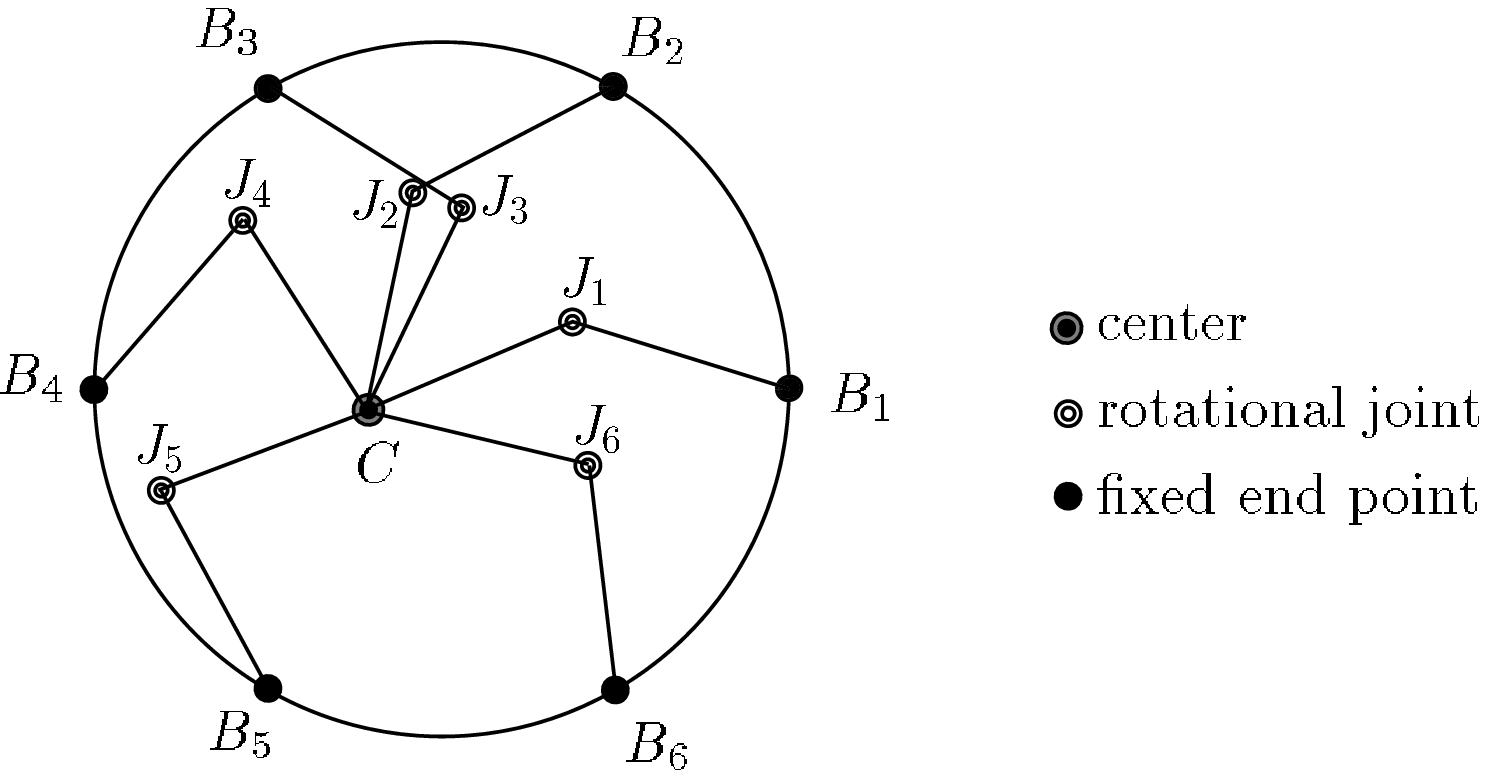}
\caption{A robot with $6$ arms}
\label{fig_6_arms}
\end{center}
\end{figure}

We study the configuration space of the linkage of a robot which can move only in a plane. 
We consider a robot which has $n$ arms such that each arm is of length $1+1$ and has a rotational joint in the middle, and that the endpoints of the arms are fixed to $n$ equally located points in a circle of radius $R$ (Figure \ref{fig_6_arms}). 
We assume that its arms and joints can intersect each other. 
Let us call this robot a ``{\em spider}'' and denote the configuration space of the spiders with $n$ arms of radius $R$ by ${\mathcal{M}}_n(R)$. 

Let $\vect x$ be a point in ${\mathcal{M}}_n(R)$ that corresponds to a spider such that none of the arms is stretched-out nor folded. 
all the angles at the joints of the arms belong to $(0,\pi)$. 
The configuration of a spider, if it is close to the above mentioned, is determined by the position of the body. 
Therefore, the neighbourhood of a generic point in ${\mathcal{M}}_n(R)$ is of dimension $2$. 

In this paper we show that ${\mathcal{M}}_n(R)$ is generically diffeomorphic to an orientable closed surface, and give the genus in terms of $n$ and $R$ (Theorem \ref{thm_M_n}) by a topological and a Morse theoretical methods. 
We also give the topological type of ${\mathcal{M}}_n(R)$ when it is not a surface (Theorem \ref{thm_M_n_sing}). 

\smallskip
The puzzle which was mentioned in Dror Bar-Natan's talk at Siegen, Germany in 2001 provoked the author to study the generalization. 
Dror Bar-Natan asked what ${\mathcal{M}}_6(R)$ is for a big $R$. 
The answer was given by calculating the Euler numer: 
the configuration space can be obtained by gluing $2^6$ hexagons at the edges in such a way that four hexagons meet at each vertex. 
D. Eldar's home page \cite{E} explicitly shows the illustration of this idea. 

\smallskip
When the number of arms of the spider is equal to $2$, by joining the two endpoints, we can identify the space ${\mathcal{M}}_2(R)$ with the configuration space of the planar quasi-equilateral pentagons that have edges of lengths $1,1,1,1,2R$. 
Especially, when $R=\frac12$ the space ${\mathcal{M}}_2(\frac12)$ is equal to the configuration space of the planar equilateral pentagons. 
It was proved to be homeomorphic to an orientable closed connected surface $\varSigma_4$ of genus $4$ by Havel (\cite{Hav}) and by Kamiyama (\cite{Ka}). 
Generally, the space ${\mathcal{M}}_2(R)$ was proved to be homeomorphic to $\varSigma_4$ if $0<2R<2$ and to $S^2$ if $2<2R<4$ by Toma (\cite{To}). 
This result is a special case of the theorem of Kapovich and Millson \cite{Ka-Mi1} et.al., where the genus of the configuration space of pentagons with edges of length $a_1, \cdots, a_5$ was given in terms of $a_1, \cdots, a_5$ when the configuration space becomes a smooth manifold. 
The singular cases, ${\mathcal{M}}_2(0)$ and ${\mathcal{M}}_2(1)$, were studied in \cite{To}. 

There have been a great number of works on the topology of the configuration spaces of polygons, linkages, and mechanisms from various viewpoints. 
For example, topological argument for general cases can be found in \cite{SSB}. 

\bigskip
This is a revised version of ArXiv math.GT/0505462 (2005). 
Because of the limit of the capacity, we had to reduce the number and the sizes of the figures. 
The full figure version can be found at 

http://www.comp.metro-u.ac.jp/\textasciitilde ohara/download/mod\_sp\_spider.pdf

\section{Main results} 

Let us give an explicit definition of ${\mathcal{M}}_n(R)$. 
We assume $n\ge 2$ in what follows. 
Let $C(x,y)$ denote the ``body" of the spider. 
Let 
\begin{equation}\label{coord_B_k}
B_k=(u_k,v_k)={\left(R\cos\frac{2(k-1)\pi}{n}, R\sin\frac{2(k-1)\pi}{n}\right)}\end{equation}
be the $k$-th fixed endpoint and $J_k(p_k,q_k)$ the joint of the $k$-th arm $(k=1,\cdots, n)$. 
We denote the vector $\overrightarrow{J_kC}$ by $\vect a_k$ and $\overrightarrow{B_kJ_k}$ by $\vect b_k$. 
All the vectors are considered row vectors. 

\begin{definition} \rm Let $R$ be a constant with $0\le R\le 2$. 
Define $f_i:\mathbb{R}^{2n+2}\to\mathbb{R}$ for $i=1,\cdots, 2n$ by 
$$\begin{array}{rccl}
f_{2k-1}(x,y,p_1,q_1,\cdots,p_n,q_{n})\!\!&\!\!=\!\!&\!\!|J_kC|^2-1\!\!&\!\!=(x-p_k)^2+(y-q_k)^2-1,\\
f_{2k}(x,y,p_1,q_1,\cdots,p_n,q_{n})\!\!&\!\!=\!\!&\!\!|B_kJ_k|^2-1\!\!&\!\!=(p_k-u_k)^2+(q_k-v_k)^2-1,
\end{array}$$
and $F:\mathbb{R}^{2n+2}\to
\mathbb{R}^{2n}$ by 
$$F
=(f_1,\cdots,f_{2n}).$$

The {\em configuration space of the spiders with $n$ arms of radius $R$}, ${\mathcal{M}}_n(R)$, is given by 
\begin{eqnarray}\label{def}
{\mathcal{M}}_n(R)\!\!&\!\!=\!\!&\!\!\displaystyle \left\{
	(C,J_1,\cdots,J_n)\in\left(\mathbb{R}^2\right)^{n+1}:\vert J_kC\vert =
	\vert B_kJ_k\vert =1 \> (k=1,\cdots,n)
	\right\} \nonumber\\[2mm]
	\!\!&\!\!=\!\!&\!\!\left\{\vect x=(x,y,p_1,q_1,\cdots,p_n,q_{n})\in\mathbb{R}^{2n+2}:
	f_i(\vect x)=0 \hspace{0.2cm} (1\le i\le 2n)	\right\}\nonumber\\[1mm]
	\!\!&\!\!=\!\!&\!\!F^{-1}(\vect 0).
\end{eqnarray}
\end{definition}
%
%
%
\begin{definition} \rm 
Let $d_n$ be the maximum distance between $n$-th roots of unity: 
$$
d_n=\left\{
\begin{array}{cl}
2 & \>\textrm{ if $n$ is even},\\[1mm]
\sqrt{2-2\cos\frac{2m\pi}{2m+1}} & \>\textrm{ if $n$ is odd, $n=2m+1$}. 
\end{array}
\right.
$$
Put 
$$R_n=\frac{2}{d_n}.$$
\end{definition}
\begin{lemma}\label{d_n}
A spider with $n$ arms can have both stretced-out arms and a folded arm {\rm (}hence the body is located at some $B_k${\rm )} if and only if $R=R_n$. 
It can have folded arms if and only if $R\le R_n$. 
\end{lemma}
\begin{theorem}\label{thm_M_n}
The configuration space of the spiders with $n$ arms of radius $R$, ${\mathcal{M}}_n(R)$, is diffeomorphic to a connected orientable closed surface $\varSigma_g$ if $R$ satisfies 
\begin{equation}\label{cond_R_non-sing}
0<R<2 \hspace{0.2cm}\textrm{and}\hspace{0.2cm}  R\ne R_n.
\end{equation}
The genus $g$ is given by 
$$
g=\left\{
\begin{array}{lll}
1-2^{n-1}+n2^{n-3}+n2^{n-1}\!\!&\!\!=1+(5n-4)2^{n-3} & \>\textrm{ if }\>\> 0<R<R_n,\\[1mm]
1-2^{n-1}+n2^{n-3}\!\!&\!\!=1+(n-4)2^{n-3} & \>\textrm{ if }\>\> R_n<R<2.
\end{array}
\right.
$$
\end{theorem}
\begin{theorem}\label{thm_M_n_sing} 
The topological type of the configuration space ${\mathcal{M}}_n(R)$ when it is not diffeomorphic to a surface is given as follows: 
\begin{enumerate}
\item If $R=0$, ${\mathcal{M}}_n(0)$ can be decomposed as ${\mathcal{M}}_n(0)=S^1\times{\mathcal{M}}^{\p}_n(0)$. 
Let $\displaystyle \stackrel{2^{n-1}}{\vee} S^1$ be a bouquet of $2^{n-1}$ circles with base point $P$. 
Then the sapce ${\mathcal{M}}^{\p}_n(0)$ can be obtained from the union of $\displaystyle \stackrel{2^{n-1}}{\vee} S^1$ and $T^{n-1}$ by gluing $2^{n-1}$ points of $\displaystyle \stackrel{2^{n-1}}{\vee} S^1\setminus\{P\}$ each of which belongs to mutually distinct circle to $2^{n-1}$ distinct points in $T^{n-1}$ {\rm (Figures \ref {conf_quadrilateral} and \ref{M'_3-0})}. 
\begin{figure}[htbp]
\begin{center}
\begin{minipage}{.3\linewidth}
\begin{center}
\includegraphics[width=\linewidth]{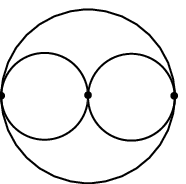}
\caption{${\mathcal{M}}^{\p}_2(0)=T^1\cup (S^1\vee S^1)/\sim$}
\label{conf_quadrilateral}
\end{center}
\end{minipage}
\hskip 1.0cm
\begin{minipage}{.55\linewidth}
\begin{center}
\includegraphics[width=\linewidth]{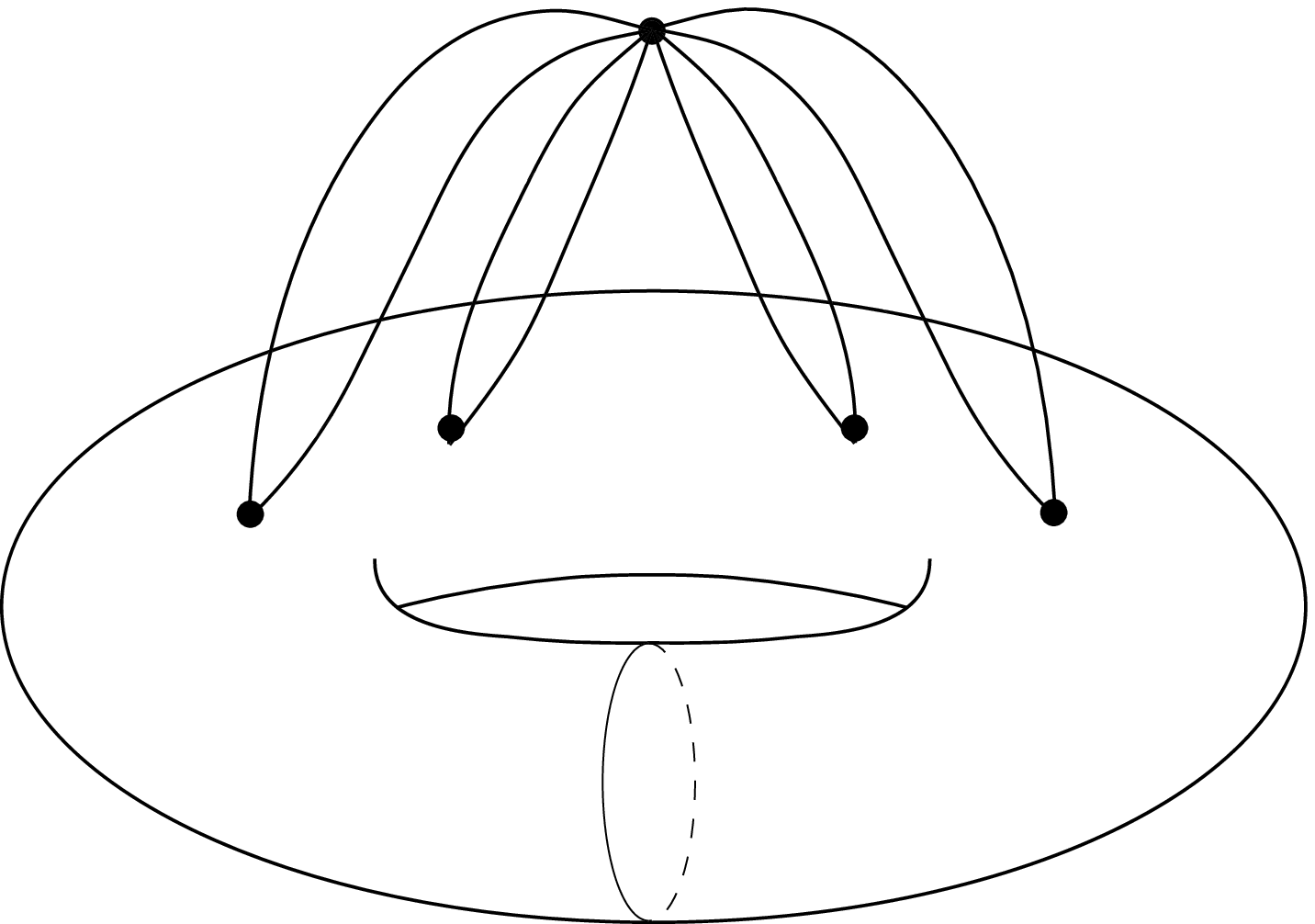}
\caption{${\mathcal{M}}^{\p}_3(0)=T^2\cup\stackrel{4}{\vee} S^1/\sim$}
\label{M'_3-0}
\end{center}
\end{minipage}
\end{center}
\end{figure}
\item If $n$ is even and $R=R_n=1$, ${\mathcal{M}}_n(1)$ can be obtained from 
$${\mathcal{M}}_n(R^{\p})\approx \varSigma_{1-2^{n-1}+n2^{n-3}+n2^{n-1}}\hspace{1cm}(0<R^{\p}<1)$$
by pinching $n2^{n-1}$ $1$-handles in the middle {\rm (Figure \ref{pinch_Sigma_4})}. 
It can also be obtained from 
$$
{\mathcal{M}}_n(R^{\p\p})\approx \varSigma_{1-2^{n-1}+n2^{n-3}}\hspace{1cm}(1<R^{\p\p}<2)
$$
by identifying $n2^{n-1}$ pairs of points respectively. 
%
\begin{figure}[htbp]
\begin{center}
\includegraphics[width=.55\linewidth]{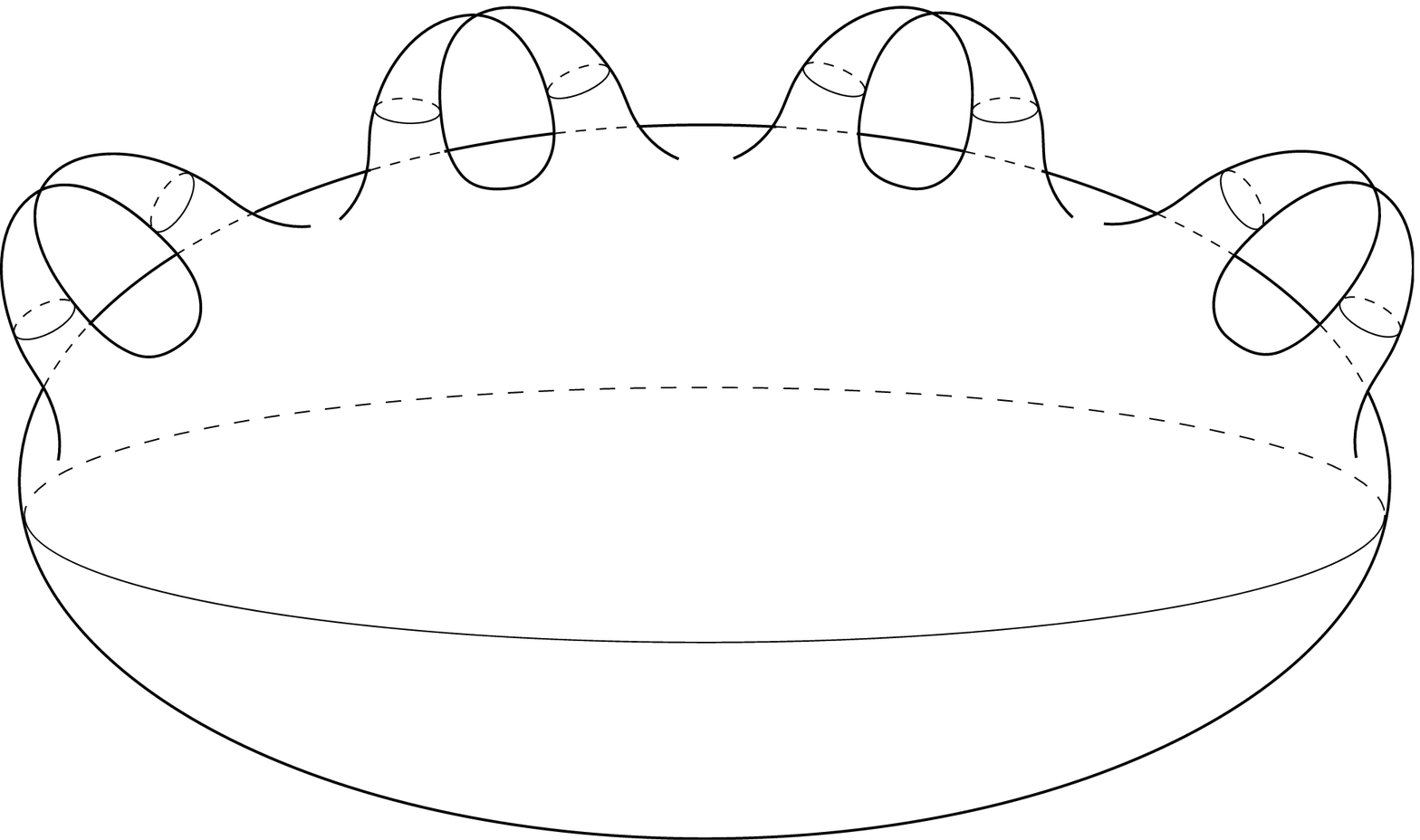}
\caption{${\mathcal{M}}_2(1)$ appears in a continuous family $\{{\mathcal{M}}_2(R)\}_{0<R<2}$ between $\varSigma_4\approx{\mathcal{M}}_2(R^{\p})$ $(0<R^{\p}<1)$ and $S^2\approx{\mathcal{M}}_2(R^{\p\p})$ $(1<R^{\p\p}<2)$. Four $1$-handles of $\varSigma_4$ are pinched at the middle.}
\label{pinch_Sigma_4}
\end{center}
\end{figure}
\item If $n$ is odd and $R=R_n$, ${\mathcal{M}}_n(R_n)$ can be obtained from 
$$
M={\mathcal{M}}_n(R^{\p\p})\approx \varSigma_{1-2^{n-1}+n2^{n-3}}\hspace{1cm}(1<R^{\p\p}<2)
$$
as follows. 
Replace $n2^{n-2}$ mutually disjoint discs $D_i$ by the same number of copies $\Delta_i$ of the space illustrated in {\rm Figure \ref{stitched_cap}}. 
Let $S_i$ and $T_i$ be the endpoints of the arc of $\Delta_i$ along which the surface is stitched up. 
Join $S_{2j}$ and $S_{2j+1}$, and $T_{2j}$ and $T_{2j+1}$ by mutually disjoint curves which do not intersect with 
$$
\left(M\setminus\bigcup_iD_i\right)\cup\bigcup_i\left(\Delta_i\setminus\{S_i, T_i\}\right),
$$
which produces ${\mathcal{M}}_n(R_n)$. 

\begin{figure}[htbp]
\begin{center}
\includegraphics[width=.4\linewidth]{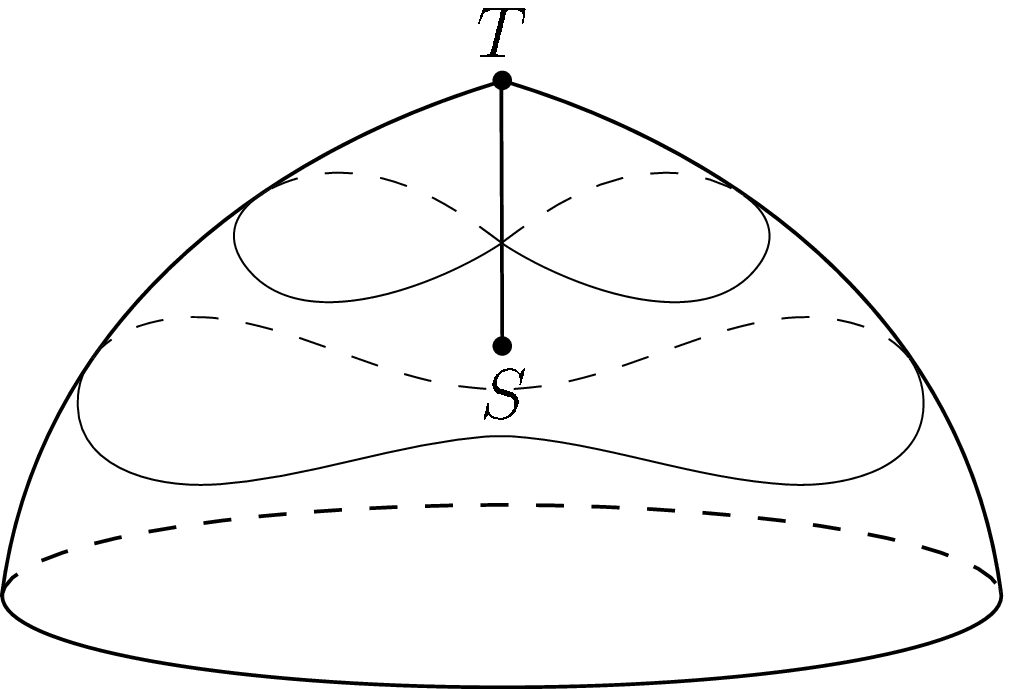}
\caption{A disc $\Delta$ which is stitched up along an arc between $S$ and $T$}
\label{stitched_cap}
\end{center}
\end{figure}
Join pairs of $[\xi_i]$ and $[\xi_{i+n2^{n-3}}]$, and pairs of $[\eta_i]$ and $[\eta_{i+n2^{n-3}}]$ by mutually disjoint arcs which intersect $M$ only at their endpoints. 
It produces ${\mathcal{M}}_n(R_n)$. 
\item If $R=2$, ${\mathcal{M}}_n(2)$ consists of one point. 
\item If $R>2$, ${\mathcal{M}}_n(R)$ is an empty set. 
\end{enumerate}
\end{theorem}

\smallskip
\begin{remark}\label{symmetry} 
The configuration space ${\mathcal{M}}_n(R)$ admits the symmetry group which is the semidirect product of the dihedral group of order $n$ (rigidly moving the $B_i$'s) and $(\mathbb{Z}/2)^n$ (interchanging $\vect a_k=\overrightarrow{J_kC}$ and $\vect b_k=\overrightarrow{B_kJ_k}$). 
We will use the symmetry of the dihedral group. 

When $R=0$ the symmetry group is the semidirect product of $O(2)$ (rotation and reflection) and $(\mathbb{Z}/2)^n$. 
We will use the symmetry of $S^1$. 
\end{remark}

\section{Proof for the non-singular case}
Theorem \ref{thm_M_n} is the consequence of the following Propositions. 
We always assume that $R$ satisfies $0<R<2$ and $R\ne R_n$ in this section. 
\begin{proposition}\label{rank_of_Jacobian}
Suppose $0<R<2$. 
The rank of Jacobian matrix $\partial F(\vect x)$ of $F$ at $\vect x$ is smaller than $2n$ if and only if $R=R_n$ and one of the following conditions is satisfied: 
\begin{enumerate}
\item[{\rm (i)}] $n$ is even and the spider has a folded arm and a stretched-out arm both of which lie on a same line. 
\item[{\rm (ii)}] $n$ is odd and the spider has a folded arm and two stretced-out arms. 
\end{enumerate}
In any case, the body is located at some $B_k$. 
\end{proposition}
%
%
\begin{corollary}\label{conf_sp=surface}
If $0<R<R_n$ or $R_n<R<2$ any connected component of ${\mathcal{M}}_n(R)=F^{-1}(\vect 0)$ is an orientable closed $2$-dimensional submanifold of $\mathbb{R}^{2n+2}$. 
\end{corollary}
\begin{proposition}\label{connectedness}
If $0<R<R_n$ or $R_n<R<2$ then the configuration space of the spiders with $n$ arms of radius $R$, ${\mathcal{M}}_n(R)$, is arcwise connected. 
\end{proposition}

The genus of ${\mathcal{M}}_n(R)=F^{-1}(\vect 0)$ is determined by calculating the Euler number. 
We have a topological and a Morse theoretical ways to do it. 
\begin{proposition}\label{topological_proof_of_the_genus}
{\rm (1)} If $R_n<R<2$ then ${\mathcal{M}}_n(R)$ admits a cell decomposition; it can be obtained by gluing $2^n$ $n$-gons at their edges so that four $n$-gons meet at each vertex. 

{\rm (2)} If $0<R<R_n$ then ${\mathcal{M}}_n(R)$ can be obtained from $\varSigma=\varSigma_{1-2^{n-1}+n2^{n-3}}$ which is homeomorphic to ${\mathcal{M}}_n(R^{\p})$ with $R_n<R^{\p}<2$ as follows. 
Blow up $n2^n$ points of $\varSigma$, i.e. replace $n2^n$ points by the same number of $S^1$'s, each point of which corresponds to the direction of the approach of a point to the blown-up point. 
Pair the $S^1$'s up and glue each pair, which is equivalent to attaching $n2^{n-1}$ $1$-handles to $\varSigma$. 
\end{proposition}
\begin{proposition}\label{Morse_theoretical_proof_of_the_genus}
Suppose $0<R<R_n$ or $R_n<R<2$. 
Let $\psi:{\mathcal{M}}_n(R)\to\mathbb{R}$ be the height function of the body of the spider: 
$$\psi(x,y,p_1,q_1,\cdots, p_n,q_n)=y.$$
Then $\psi$ is a Morse function on ${\mathcal{M}}_n(R)$. 
The number of critical points and their indices of $\psi$ are given as follows. 
\begin{enumerate}
\item If $R_n<R<2$ then there are $2^{n-2}$ critical points of index $0$, $(n-2)2^{n-2}$ critical points of index $1$, and $2^{n-2}$ critical points of index $2$. 
\item If $0<R<R_n$ and $n$ is even then there are $2^{n-2}$ critical points of index $0$, $(n-2)2^{n-2}+n2^n$ critical points of index $1$, and $2^{n-2}$ critical points of index $2$. 
\item If $0<R<R_n$ and $n$ is odd then there are $2^{n-1}$ critical points of index $0$, $n2^{n-2}+n2^n$ critical points of index $1$, and $2^{n-1}$ critical points of index $2$. 
\end{enumerate}
\end{proposition}

\subsection{Proof of ${\mathcal{M}}_n(R)$ being an orientable surface}
\begin{lemma}\label{stretched-out_arms}
Suppose $0<R<2$. 
Then the following holds. 
\begin{enumerate}
\item If two arms of the spider are stretched out then the two arms are adjacent. 
\item The spider cannot have three or more arms stretched out. 
\end{enumerate}
\end{lemma}
\begin{proof} Let $\Gamma_R(\vect 0)$ denote the circle with center the origin and radius $R$, and $\Gamma_2(C)$ the circle with center $C$ (the body of the spider) and radius $2$. 

\smallskip
(1) Suppose $j$-th and $k$-th arms ($j-k\not\equiv \pm1$ (mod $n$)) are stretched out. 
Then one of the two open subarcs of $\Gamma_R(\vect 0)$ between $B_j$ and $B_k$ is outside the circle $\Gamma_2(C)$. 
It contains at least one fixed endpoint, say, $B_i$. 
Then $|B_iC|>2>R$, which is a contradiction. 

\smallskip
(2) Suppose the $i$-th, $j$-th, and $k$-th arms are stretched out. 
Then both $\Gamma_R(\vect 0)$ and $\Gamma_2(C)$ pass through $B_i$, $B_j$, and $B_k$. 
As there is a unique circle through three points, $\Gamma_R(\vect 0)=\Gamma_2(C)$, which contradicts the condition $R\ne 2$.  
\end{proof}

Put 
$$\vect a_k=\overrightarrow{J_kC}=(x-p_k,y-q_k), \>\>\vect b_k=\overrightarrow{B_kJ_k}=(p_k-u_k, q_k-v_k).$$
Then the Jacobian matrix $\partial F(\vect x)$ of $F$ at $\vect x\in\mathbb{R}^{2n+2}$ is given by 
\begin{equation}\label{Jacobian_matrix_of_F}
\partial F(\vect x)=
\left(\begin{array}{c}
\partial f_1(\vect x)\\
\partial f_2(\vect x)\\
\partial f_3(\vect x)\\
\partial f_4(\vect x)\\
\vdots\\
\partial f_{2n-1}(\vect x)\\
\partial f_{2n}(\vect x)
\end{array}\right)
=2\left(
\begin{array}{ccccc}
\> \vect a_1 \>&\> -\vect a_1 \> &&& \\
& \vect b_1 &&& \\
\> \vect a_2 \>&&\> -\vect a_2 \> && \\
&& \vect b_2 && \\
&&&\ddots &\\
\> \vect a_n \>&&&&\> -\vect a_n \> \\
&&&& \vect b_n  
\end{array}
\right). 
\end{equation}
%
We may denote $\partial f_k(\vect x)$ by $\partial f_k$. 

\begin{proofofproposition}\ref{rank_of_Jacobian}. 
(1) 
Suppose 
$$\sum_{k=1}^{2n}c_k\partial f_k=\vect 0 \hspace{0.5cm} \textrm{with} \hspace{0.5cm} (c_1, \cdots, c_{2n})\ne (0, \cdots, 0).$$ 

%
If $c_1=c_3=\cdots=c_{2n-1}=0$ then $c_2=c_4=\cdots=c_{2n}=0$, which is a contradiction. 
Therefore, at least one of $c_{2i-1}$'s does not vanish. 
Since $\vect a_i\ne\vect 0$, it implies that at least two of $c_{2i-1}$'s do not vanish. 
If $c_{2i-1}\ne 0$ then $c_{2i}\ne 0$; hence $\vect a_i=\pm \vect b_i$, i.e. the $i$-th arm is either stretched out or folded. 

Case I. Suppose the $i$-th arm is folded. 
Then the body $C$ of the spider is located at $B_i$. 
Therefore, there are no more folded arms. 
If there are more than two non-zero $c_{2j-1}$'s besides $c_{2i-1}$ then there are more than two stretced out arms, which contradicts Lemma \ref{stretched-out_arms}. 
Therefore, there are one or two stretched-out arms. 

If there are two stretched-out arms, they are from the farest $B_j$'s from $B_i$, which can occur if and only if $n$ is odd and $R=R_n$. 
In this case, $\partial f_k$'s are in fact linearly dependent. 
This corresponds to the case (ii). 

If there is only one stretched-out arm, it is from the unique farest $B_j$ from $B_i$, which can occur if and only if $n$ is even and $R=R_n=1$. 
Since there are no more non-zero $c_k$'s besides $c_{2i-1}, c_{2i}, c_{2j-1}$, and $c_{2j-1}$, $\partial f_k$'s are linearly dependent if and only if $\vect a_i=-\vect b_i=\pm \vect a_j=\pm \vect b_j$, in other words, the $i$-th and the $j$-th arms lie on the same line $\overline{B_iB_j}$. 
This corresponds to the case (i). 

Case II. Suppose there are no folded arms. 
It follows that there are exactly two non-zero $c_{2i-1}$'s and two stretched-out arms. 
Then $\partial f_k$'s are linearly dependent if and only if these two stretched-out arms lie on the same line, which contradicts the condition that $R\ne2$. 
\end{proofofproposition}
\subsection{Connectedness of the configuration space ${\mathcal{M}}_n(R)$}
The configuration of the spider is determined by two kinds of data; 
the position of the body $C(x,y)$, and the state of the $n$ arms. 

\begin{figure}[htbp]
\begin{center}
\includegraphics[width=.3\linewidth]{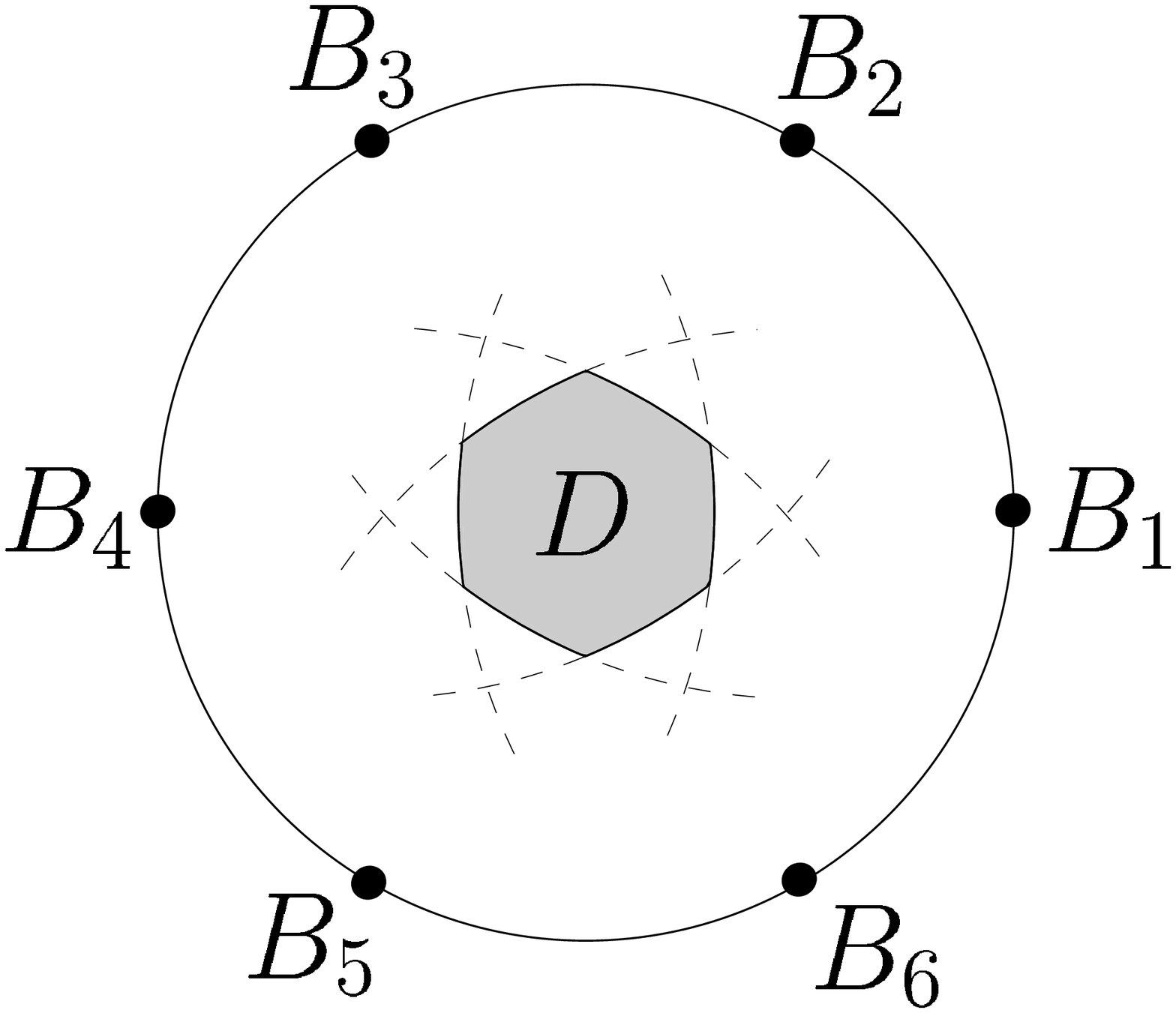}
\caption{The {\em curved hexagon} $D$}
\label{6arms-domain_a}
\end{center}
\end{figure}

The former is given by a point in the domain where the body can be located, which we shall call the {\em body domain}. 
In our case it is a ``{\sl curved $n$-gon}'' $D$ (Figure \ref{6arms-domain_a}) given by 
\begin{equation}\label{curved_n-gon}
D=\{C=(x,y):|CB_k|\le 2 \hspace{0.5cm}(1\le k\le n)\}. 
\end{equation}

Since $R\ne R_n$ the boundary $\partial D$ of $D$ does not contain $B_k$. 

\begin{figure}[htbp]
\begin{center}
\begin{minipage}{.25\linewidth}
\begin{center}
\includegraphics[width=\linewidth]{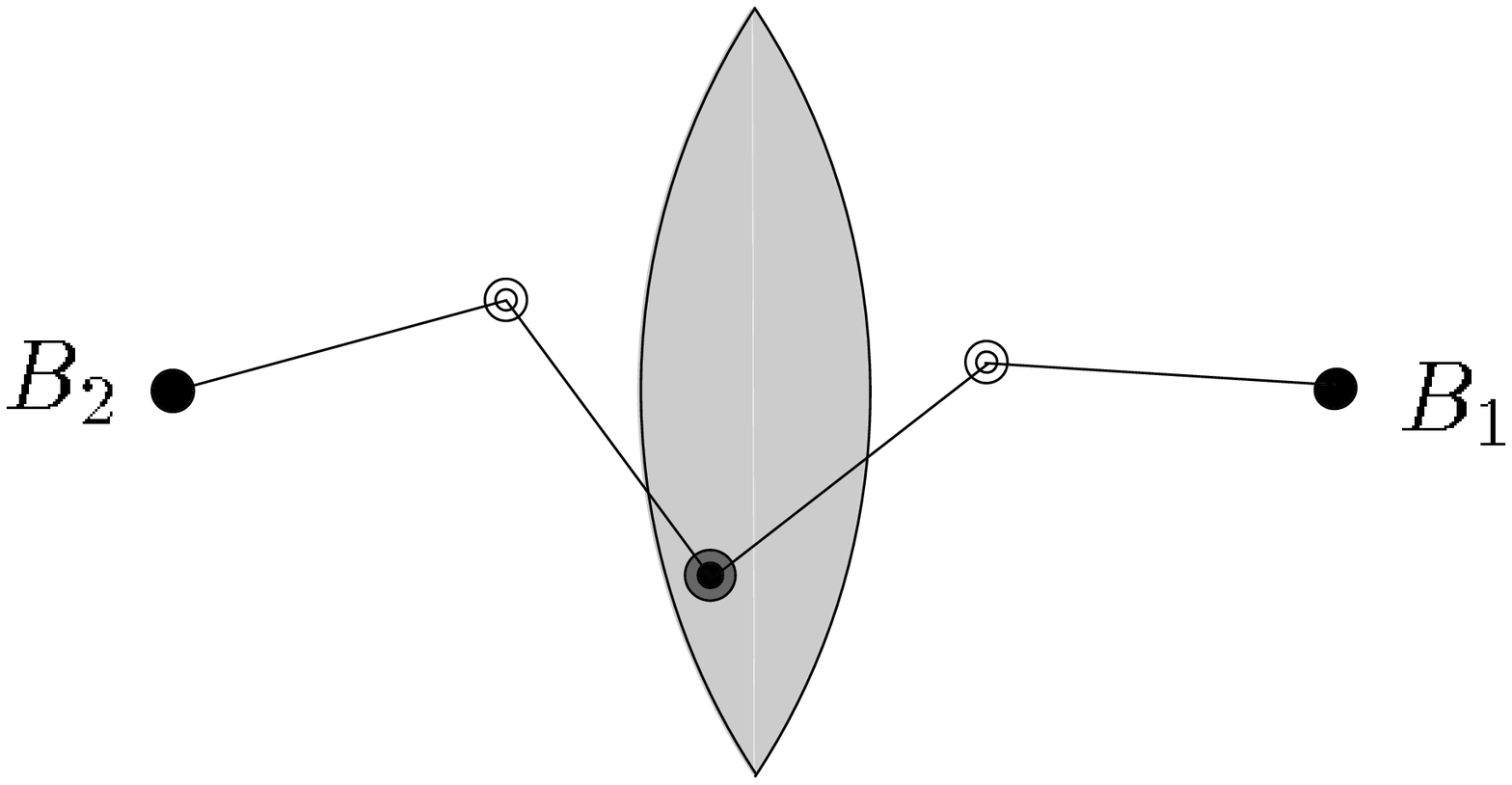}
\caption{}
\label{D+-_a}
\end{center}
\end{minipage}
\begin{minipage}{.25\linewidth}
\begin{center}
\includegraphics[width=\linewidth]{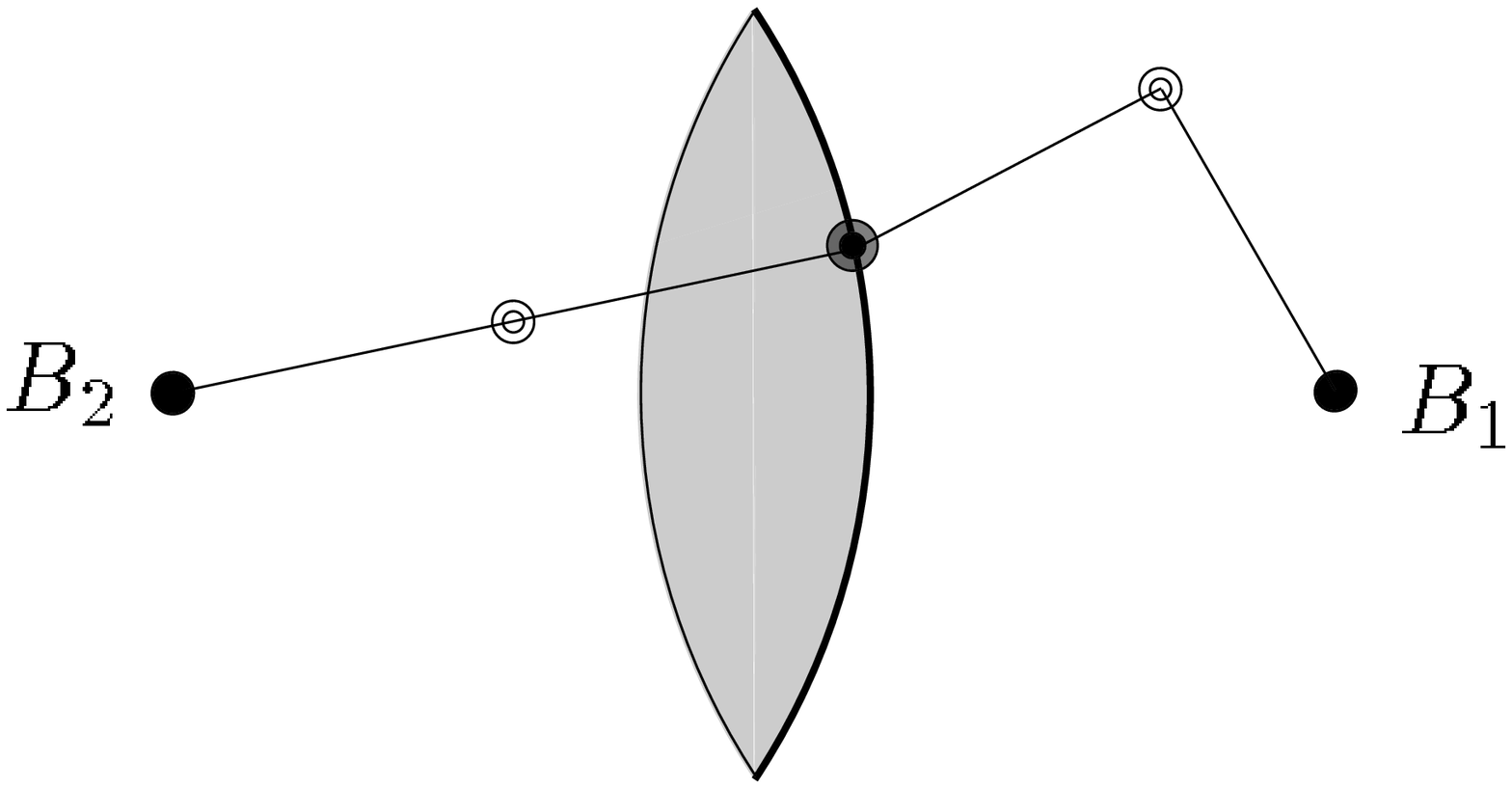}
\caption{}
\label{E+0_a}
\end{center}
\end{minipage}
\begin{minipage}{.25\linewidth}
\begin{center}
\includegraphics[width=\linewidth]{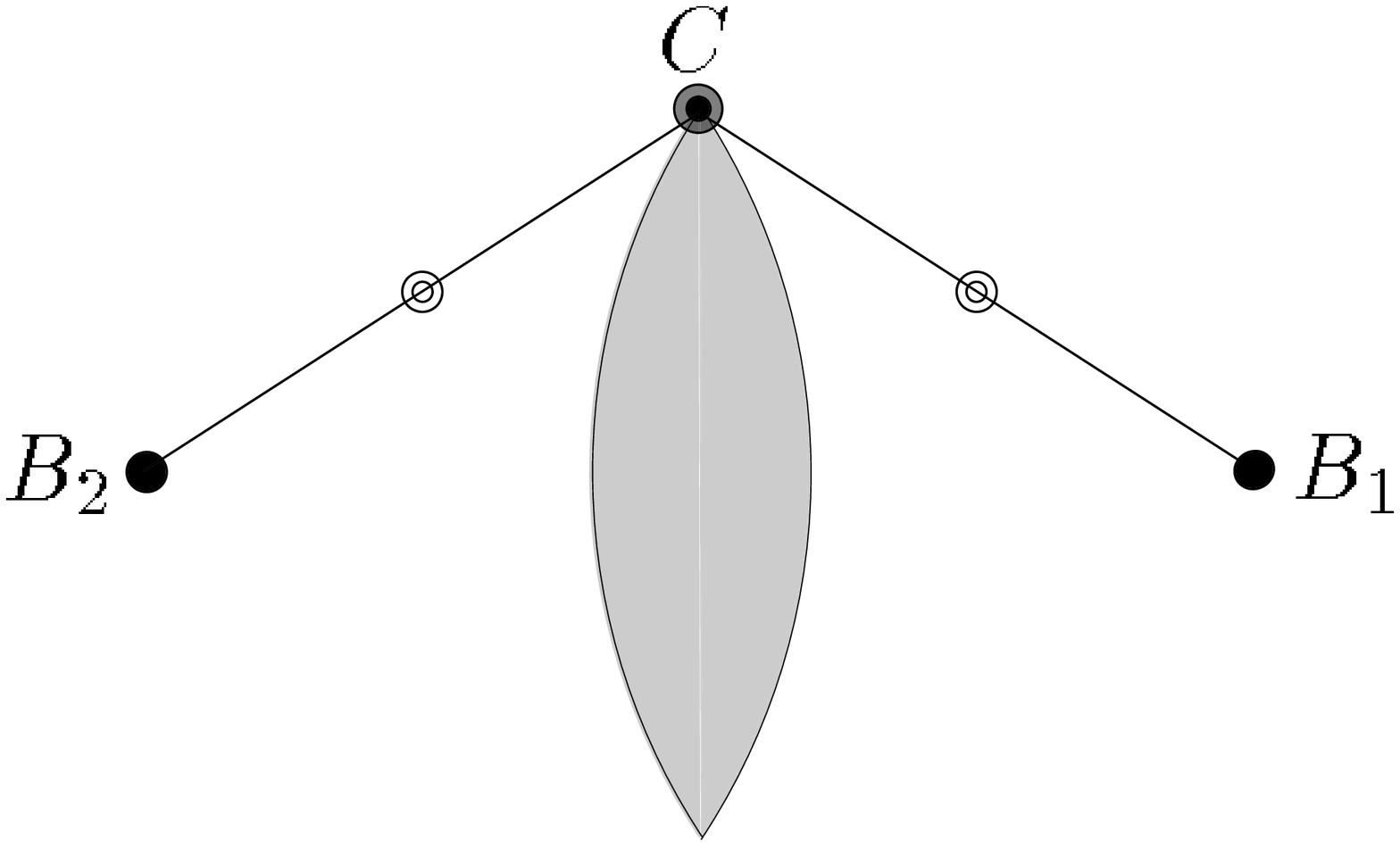}
\caption{}
\label{V00u_a}
\end{center}
\end{minipage}
\begin{minipage}{.215\linewidth}
\begin{center}
\includegraphics[width=\linewidth]{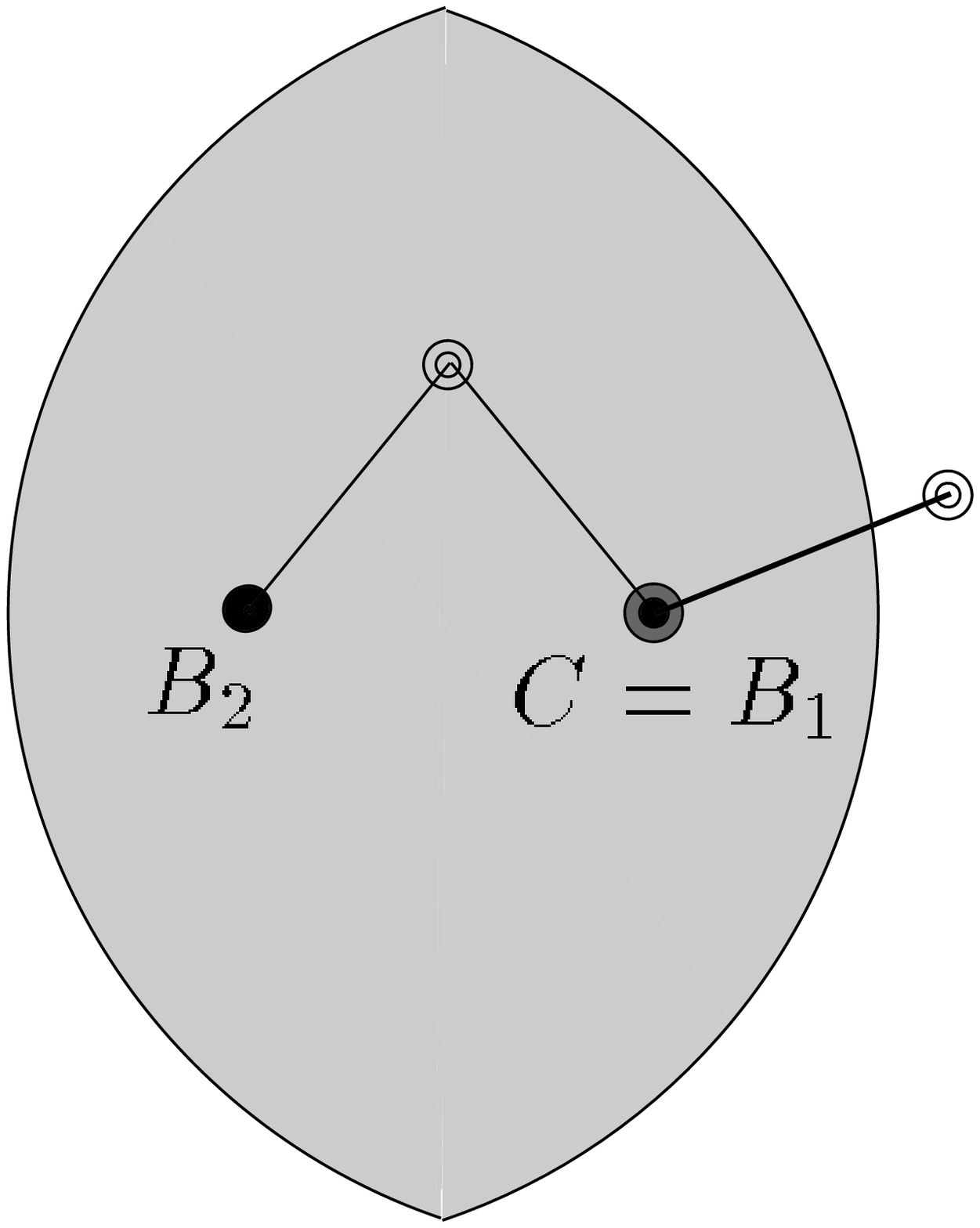}
\caption{}
\label{dilateral_R1_a}
\end{center}
\end{minipage}
\end{center}
\end{figure}

The latter depends on the position of the body of the spider. 
There are three (when $R_n<R<2$) or four (when $0<R<R_n$) mutually disjoint cases: 
\begin{enumerate}
\item The body of the spider is located in the interior of the curved $n$-gon $D$, but not at $B_k$. 
There are $2^n$ states how the arms are bended (Figure \ref{D+-_a}). 
\item Exactly one arm is stretched out. 
It occurs if and only if the body of the spider is located on an interior of an edge of $D$ (Figure \ref{E+0_a}). 
\item Exactly two arms are stretched out. 
It occurs if and only if the body of the spider is located at a vertex of $D$ (Figure \ref{V00u_a}). 
\item The body of the spider is located at $B_k$. 
The $k$-th arm, which is folded, can rotate around $B_k$ (Figure \ref{dilateral_R1_a}). 
It can occur only when $0<R<R_n$. 
\end{enumerate}

\begin{figure}[htbp]
\begin{center}
\includegraphics[width=.45\linewidth]{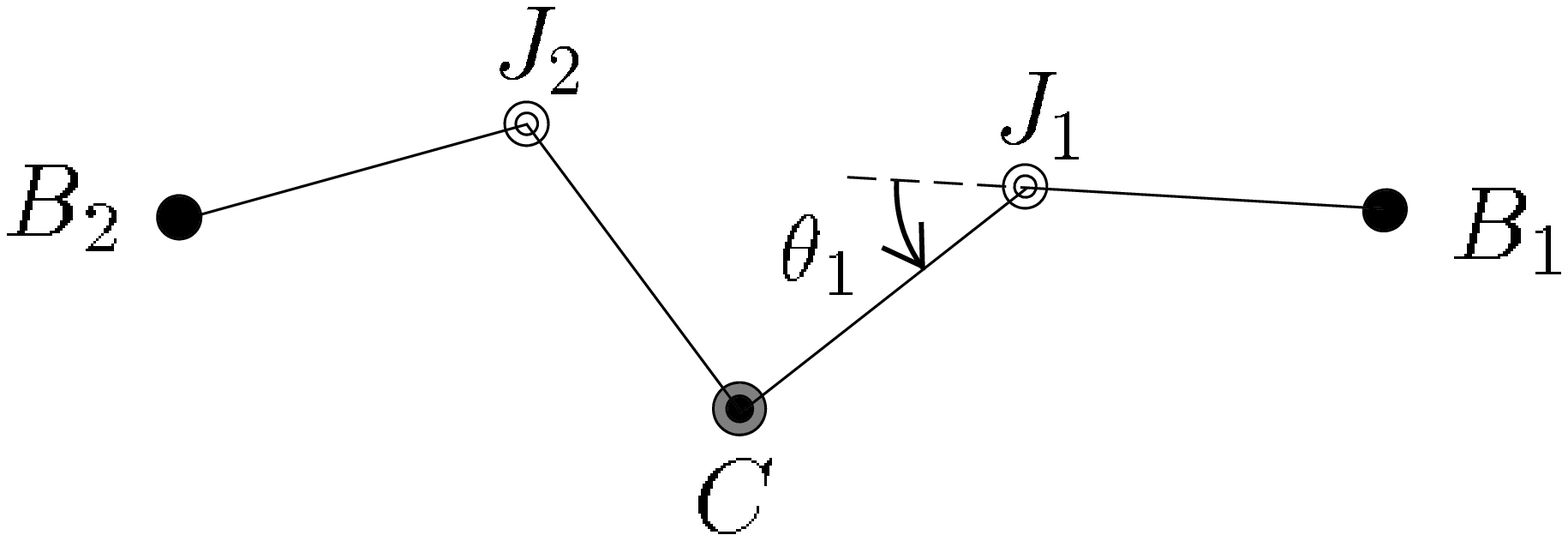}
\caption{$\varepsilon_1=+, \varepsilon_2=-$. The case when $n=2$. }
\label{fig_2arms+-}
\end{center}
\end{figure}

\begin{definition}\label{def_index_arm} \rm 
Let $\theta_k$ $(-\pi<\theta_k\le\pi)$ be the angle from $\overrightarrow{B_kJ_k}$ to $\overrightarrow{J_kC}$. 
The {\em index} of the $k$-th arm, $\varepsilon_k\in\{+,-,0,\infty\}$, is given by the signature of $\tan\tanufrac{\theta_k}{2}$, where $-\infty$ is identified with $\infty$ (Figure \ref{fig_2arms+-}). 
We say that the $k$-th arm is {\em positively bended} (or {\em negatively bended}) if its index $\varepsilon_k$ is $+$ (or respectively, $-$), {\em bended} if it is either positively or negatively bended. 
We note that it is {\em stretched-out} if $\varepsilon_k=0$, and {\em folded} if $\varepsilon_k=\infty$. 
\end{definition}
\begin{definition}\label{int_D} \rm \begin{enumerate}
\item We call $D\cap\{C=(x,y):|CB_k|=2\}$ the {\em $k$-th edge} of $D$. 
\item Define $\stackrel{\circ}D$ by 
$$
\stackrel{\circ}D=\{C=(x,y):0<|CB_k|<2 \hspace{0.5cm}(1\le k\le n)\}, 
$$
namely, 
$$
\stackrel{\circ}D=\left\{
	\begin{array}{ll}
	\textrm{Int} D\setminus\{B_1, \cdots, B_n\} & \textrm{when}\hspace{0.5cm} 0<R<R_n\\
	\textrm{Int} D & \textrm{when}\hspace{0.5cm} R_n<R<2, 
	\end{array}
\right.
$$
and call it the {\em open body domain}. 
It is the domain where the body of a spider can be located whose arms are all bended. 
\end{enumerate}
\end{definition}
\begin{remark}
The indices of the arms are kept invariant while the body of the spider moves around inside $\stackrel{\circ}D$. 
\end{remark}
\begin{proofofproposition}\ref{connectedness}
We show that any given spider $\vect x\in{\mathcal{M}}_n(R)$ can be deformed continuously to a fixed configuration $\vect x^{\mbox{\footnotesize$\vect 0$}}_+$ where the body is located at the origin and every arm has index $+$. 

First, deform the spider continuously so that each arm of it has index either $+$ or $-$. 
This can be done by moving the body a little bit to a point in $\stackrel{\circ}D$. 

Second, change the $-$ indices to $+$ one by one, by iteration of stretching out a negatively bended arm and then bending it again positively without changing the indices of the other arms. 
Suppose the $k$-th arm is negatively bended. 
Move the body through $\stackrel{\circ}D$ to a point in the inerior of the $k$-th edge of $D$, and then move it inward to make the $k$-th arm positively bended. 

Finally, move the body to the origin through $\stackrel{\circ}D$ to complete the proof. 
\end{proofofproposition}

\begin{remark}
The configuration space of an ``{\sl asymmetric spider}'' may be disconnected. 
It happens when there is an arm whose index cannot be changed by any continuous motion of the body (Figure \ref{asymmetric}). 
\begin{figure}[htbp]
\begin{center}
\includegraphics[width=.4\linewidth]{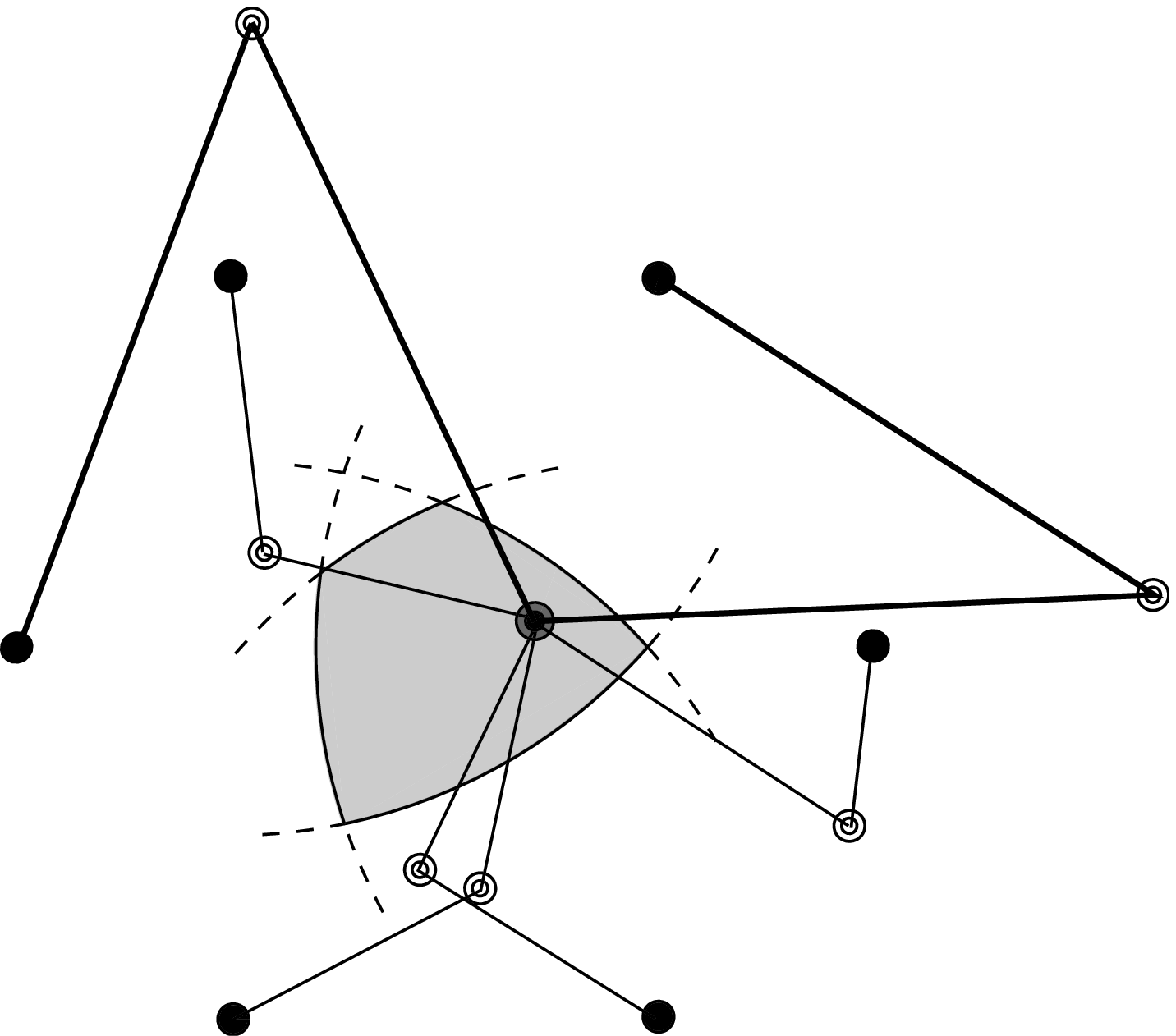}
\caption{An {\sl asymmetric spider}. Its second and fourth arms cannot be stretched out nor folded. }
\label{asymmetric}
\end{center}
\end{figure}
When $|B_kJ_k|\ne|J_kC|$ the open domain $\stackrel{\circ}D$ is replaced by the intersection of open annuli. 
\end{remark}

\subsection{Topological method to determine the genus}
\begin{proofofproposition}\ref{topological_proof_of_the_genus}. 

(i) {\bf The {\boldmath $R_n<R<2$} case. }

We remark that any arm of a spider cannot be folded in this case. 

We give a cell decomposition $\Delta$ of ${\mathcal{M}}_n(R)$ as follows. 

Let $\vect{\varepsilon}=\vect{\varepsilon}(\vect x)$ denote the multi-index of the arms of a spider $\vect x\in {\mathcal{M}}_n(R)$ (Definition \ref{def_index_arm}): 
$$\vect{\varepsilon}=(\e_1, \cdots, \e_n), \hspace{1cm} \e_k\in\{+,-,0\}. $$
Lemma \ref{stretched-out_arms} implies that $\vect{\varepsilon}(\vect x)$ contains at most two $0$'s, and if so, they are adjacent modulo $n$. 

Let ${\mathcal{I}}$ denote the set of the multi-indices of the points in ${\mathcal{M}}_n(R)$: 
\begin{equation}\label{f_multi-indices}
{\mathcal{I}}=\left\{\vect{\varepsilon}=(\e_1, \cdots, \e_n)\left|
\begin{array}{l}
\e_k\in\{+,-,0\}, \\
\sharp\{j:\e_j=0\}\le 2,\\
\textrm{If }\>\>\e_i=\e_j=0\>\>\textrm{then}\>\>j-i\equiv \pm 1 \>\>(\textrm{mod}\>\> n)
\end{array}\right.
\right\}.
\end{equation}
%
Let ${\mathcal{I}}_m\subset {\mathcal{I}}$ be the set of multi-indices of a spider with $2-m$ stretched-out arms $(0\le m\le 2)$: 
\begin{equation}\label{f_I_m}
{\mathcal{I}}_m=\left\{\vect{\varepsilon}=(\e_1, \cdots, \e_n)\in{\mathcal{I}}:\,\sharp\{j:\e_j=0\}=2-m\right\}. 
\end{equation}
Define $D_{\vect{\varepsilon}}$ $(\vect{\varepsilon}\in{\mathcal{I}}_2)$, $E_{\vect{\varepsilon}}$ $(\vect{\varepsilon}^{\p}\in{\mathcal{I}}_1)$, and $V_{\vect{\varepsilon}}$ $(\vect{\varepsilon}^{\p\p}\in{\mathcal{I}}_0)$ by
%
\begin{eqnarray}\label{}
D_{\vect{\varepsilon}}\!\!&\!\!=\left\{\vect x\in {\mathcal{M}}_n(R):\,\vect{\varepsilon}(\vect x)=\vect{\varepsilon}\right\} &  (\vect{\varepsilon}\in{\mathcal{I}}_2),\label{f_D_e}\\[1mm]
E_{\vect{\varepsilon}^{\p}}\!\!&\!\!=\left\{\vect x\in {\mathcal{M}}_n(R):\,\vect{\varepsilon}(\vect x)=\vect{\varepsilon}^{\p}\right\} &  (\vect{\varepsilon}^{\p}\in{\mathcal{I}}_1),\label{f_E_e}\\[1mm]
V_{\vect{\varepsilon}^{\p\p}}\!\!&\!\!=\left\{\vect x\in {\mathcal{M}}_n(R):\,\vect{\varepsilon}(\vect x)=\vect{\varepsilon}^{\p\p}\right\} &  (\vect{\varepsilon}^{\p\p}\in{\mathcal{I}}_0).\label{f_V_e}
\end{eqnarray}
%
\begin{remark}\rm 
 We may write $D_{\e_1 \cdots \e_n}$, $E_{\e_1 \cdots \e_n}$, etc. instead of $D_{(\e_1, \cdots, \e_n)}$, $E_{(\e_1, \cdots, \e_n)}$ etc. in Figures. 
\end{remark}

Each $D_{\vect{\varepsilon}}$ is homeomorphic to $\stackrel{\circ}D=\textrm{Int}D$ since any point in $D_{\vect{\varepsilon}}$ can be identified by the position of its body as the indices of the arms are constant on $D_{\vect{\varepsilon}}$ (Figure \ref{fig_D+-}). 
Hence it is a $2$-cell of ${\mathcal{M}}_n(R)$. 
Similarly, each $E_{\vect{\varepsilon}^{\p}}$ is homeomorphic to an open interval since the body is located in the interior of the $k$-th edge if $\e^{\p}_k=0$ (Figure \ref{fig_E+0}). 
Hence it is a $1$-cell. 
Each $V_{\vect{\varepsilon}^{\p\p}}$ consists of $0$-cell(s) (a point when $n\ge 3$ or a pair of points when $n=2$) (Figure \ref{fig_V00}). 

\begin{figure}[htbp]
\begin{center}
\begin{minipage}{.315\linewidth}
\begin{center}
\includegraphics[width=\linewidth]{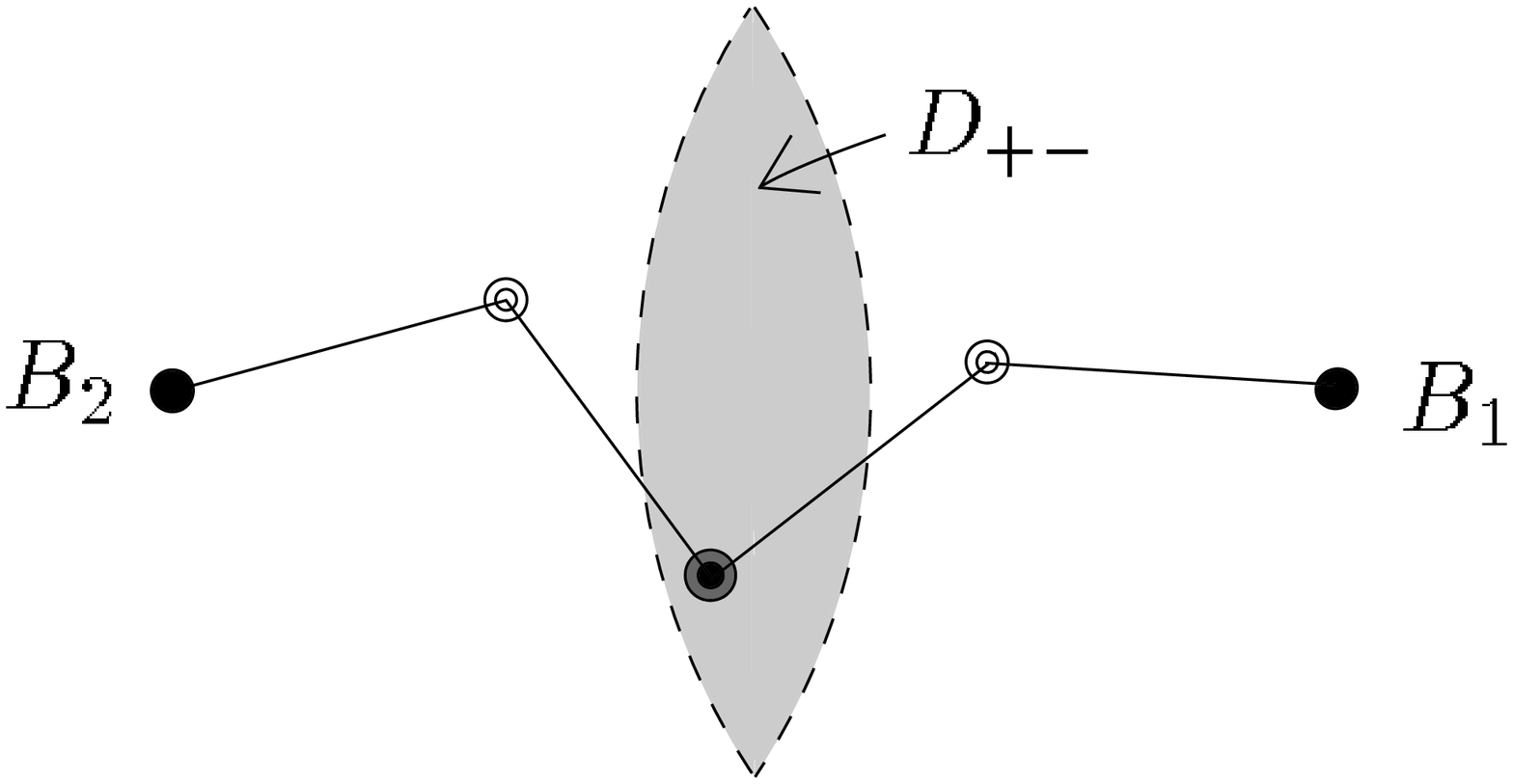}
\caption{The $2$-cell $D_{+-}$}
\label{fig_D+-}
\end{center}
\end{minipage}
\hskip 0.1cm
\begin{minipage}{.315\linewidth}
\begin{center}
\includegraphics[width=\linewidth]{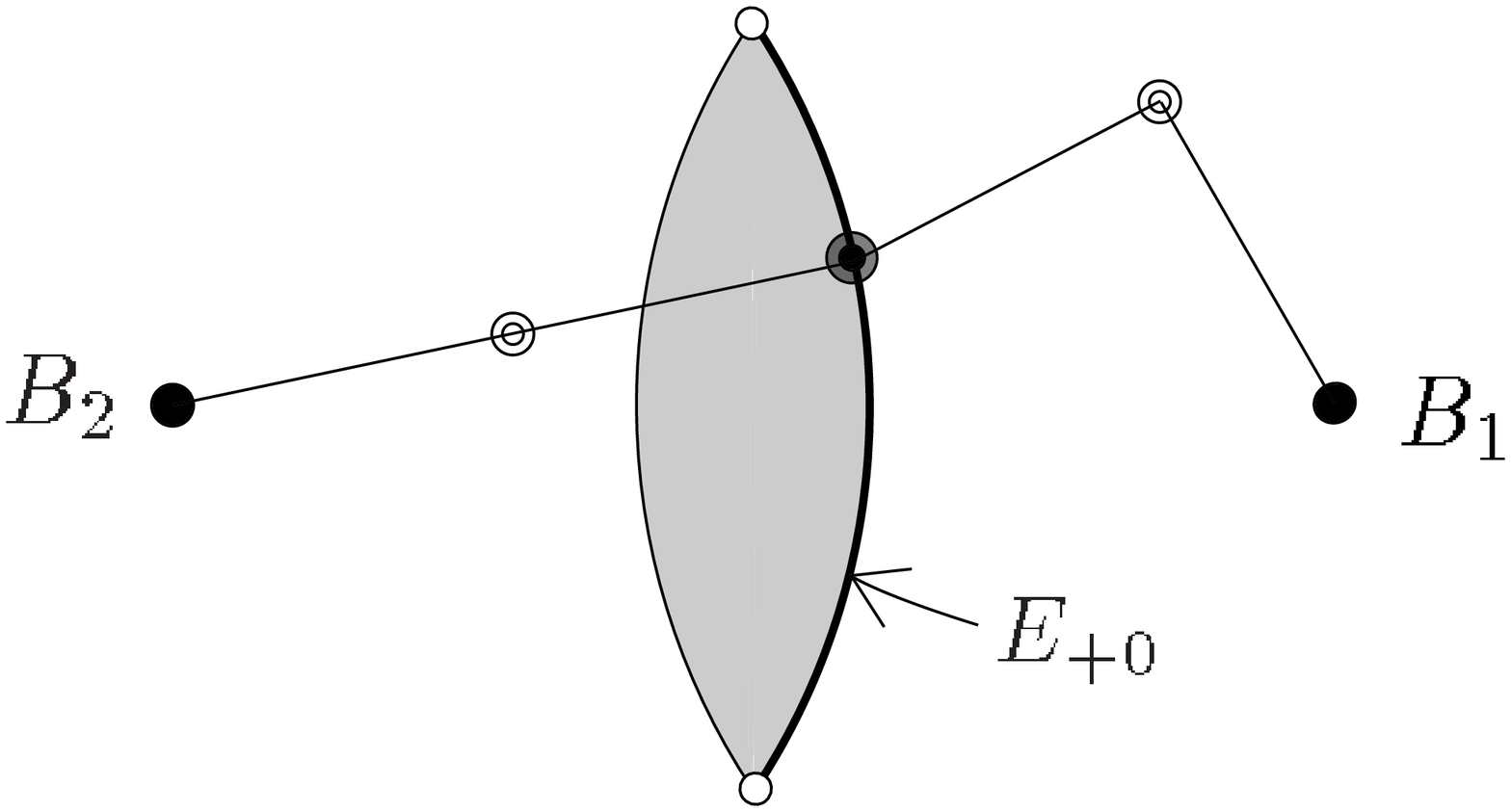}
\caption{The edge $E_{+0}$}
\label{fig_E+0}
\end{center}
\end{minipage}
\hskip 0.1cm
\begin{minipage}{.315\linewidth}
\begin{center}
\includegraphics[width=\linewidth]{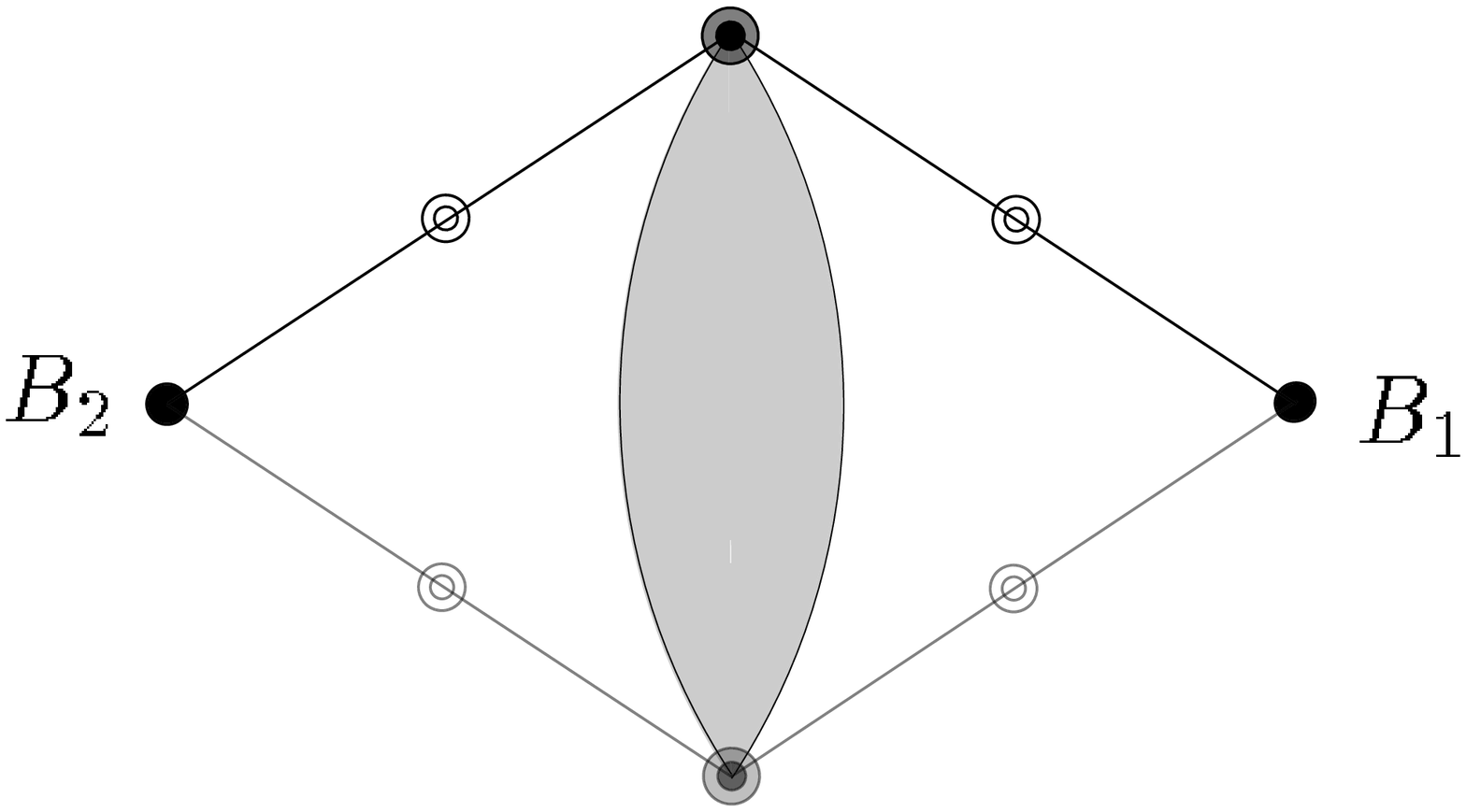}
\caption{A pair of points of $V_{00}$}
\label{fig_V00}
\end{center}
\end{minipage}
\end{center}
\end{figure}

Let $\vect{\varepsilon}=(\e_1, \cdots, \e_n)\in {\mathcal{I}}$. 
The closures of $D_{\vect{\varepsilon}}$ in ${\mathcal{M}}_n(R)$ is given by 
\begin{equation}\label{f_closure_D_e}
\overline{D_{\vect{\varepsilon}}}=D_{\vect{\varepsilon}}
\cup\bigcup_{k=1}^nE_{\vect{\varepsilon}_k^{\p}}
\cup\bigcup_{k=1}^nV_{\vect{\varepsilon}_{k,k+1}^{\p\p}},
\end{equation}
where $\vect{\varepsilon}_k^{\p}$ and $\vect{\varepsilon}_{k,k+1}^{\p\p}$ are given by 
\begin{equation}\label{e^p_k-e^pp_k}
\begin{array}{rl}
\vect{\varepsilon}_k^{\p}=(\e_1^{\p}, \cdots, \e_n^{\p}) &\displaystyle \mbox{with}\>\>
\left\{\begin{array}{l}
\e_j^{\p}=\e_j \>\>\>\mbox{if}\>\>\> j\ne k,\\[1mm]
\e_k^{\p}=0,
\end{array}
\right.\\[5mm]
\vect{\varepsilon}_{k,k+1}^{\p\p}=(\e_1^{\p\p}, \cdots, \e_n^{\p\p}) &\displaystyle \mbox{with}\>\>
\left\{\begin{array}{l}
\e_j^{\p\p}=\e_j \>\>\>\mbox{if}\>\>\> j\ne k,k+1,\\[1mm]
\e_k^{\p\p}=\e_{k+1}^{\p\p}=0,
\end{array}
\right.
\end{array}
\end{equation}
where the suffix is considered modulo $n$. 
%
With this notation, the closures of $E_{\vect{\varepsilon}_k^{\p}}$ in ${\mathcal{M}}_n(R)$ is given by 
$$
\overline{E_{\vect{\varepsilon}_k^{\p}}}=E_{\vect{\varepsilon}_k^{\p}}
\cup V_{\vect{\varepsilon}_{k-1,k}^{\p\p}}\cup V_{\vect{\varepsilon}_{k,k+1}^{\p\p}}. 
$$
%
\begin{figure}[htbp]
\begin{center}
\begin{minipage}{.5\linewidth}
\begin{center}
\includegraphics[width=\linewidth]{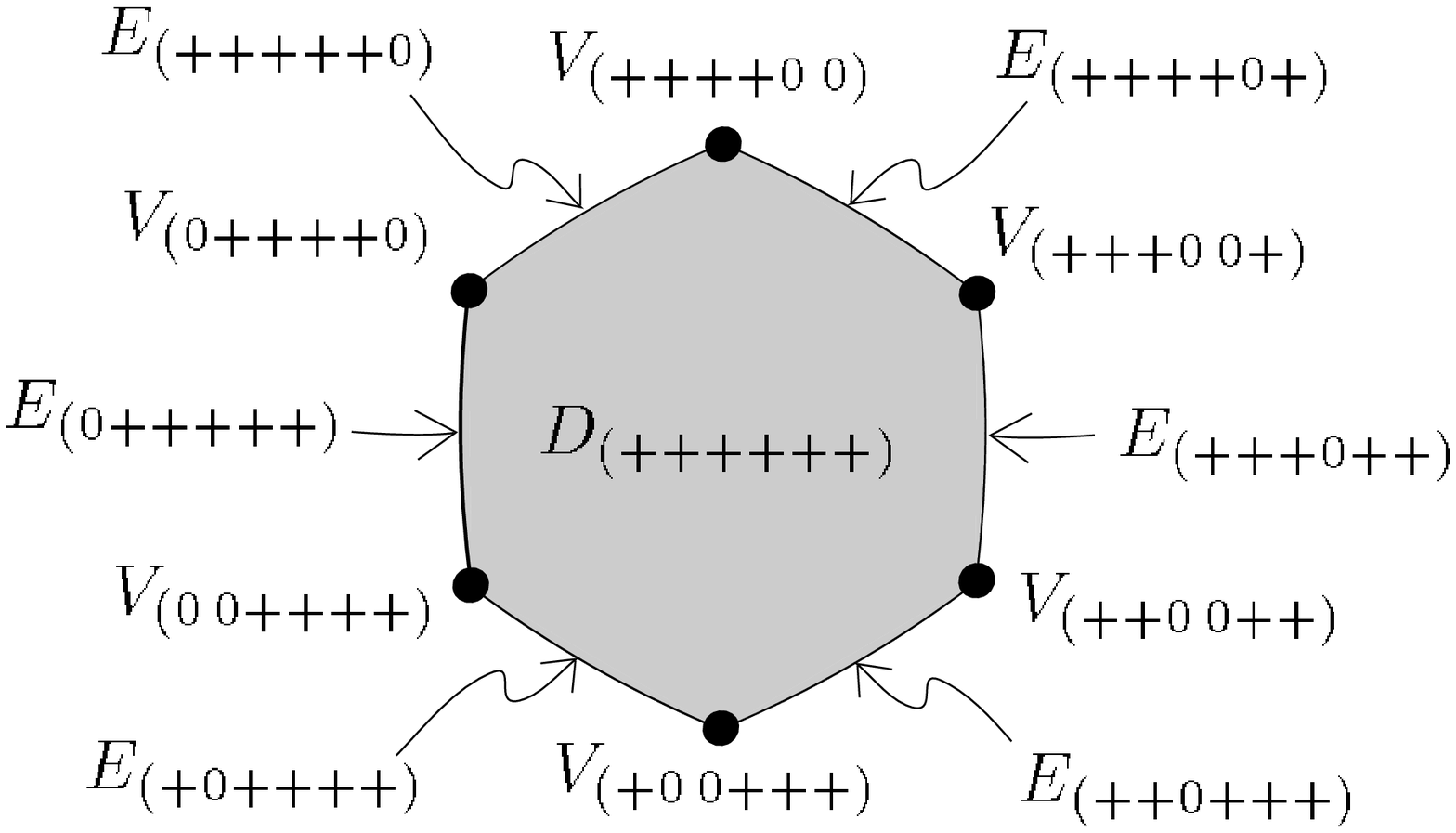}
\caption{A $2$-cell $D_{(++++++)}$}
\label{6arms-DEV}
\end{center}
\end{minipage}
\hskip 0.4cm
\begin{minipage}{.4\linewidth}
\begin{center}
\includegraphics[width=\linewidth]{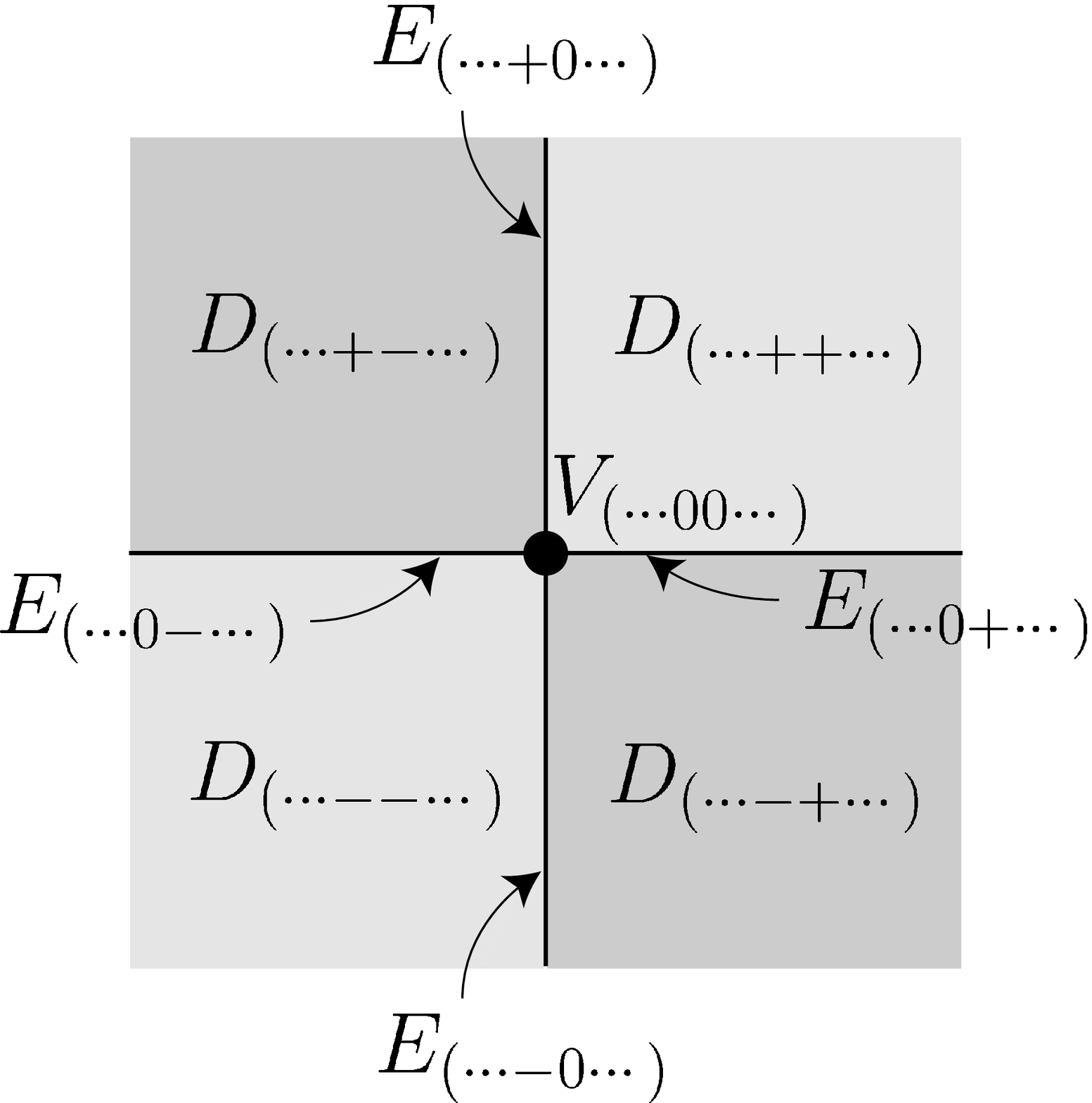}
\caption{Around a vertex}
\label{around_vertex}
\end{center}
\end{minipage}
\end{center}
\end{figure}
Therefore, the decomposition to the disjoint union: 
$$
{\mathcal{M}}_n(R)=\bigcup_{\vect{\varepsilon}\in{\mathcal{I}}_2}D_{\vect{\varepsilon}}
\cup\bigcup_{\vect{\varepsilon}^{\p}\in{\mathcal{I}}_1}E_{\vect{\varepsilon}^{\p}}
\cup\bigcup_{\vect{\varepsilon}^{\p\p}\in{\mathcal{I}}_0}V_{\vect{\varepsilon}^{\p\p}},
$$
gives a cell decomposition of ${\mathcal{M}}_n(R)$. 
The formula (\ref{f_closure_D_e}) implies that any vertex $V_{(\cdots 00 \cdots)}$ is contained in the closures of exactly four $2$-cells, 
$$\mbox{$D_{(\cdots ++ \cdots)}$, $D_{(\cdots +- \cdots)}$, $D_{(\cdots -+ \cdots)}$, and $D_{(\cdots -- \cdots)}$, }$$
where we agre that the other indices are the same (Figure \ref{around_vertex}). 

\medskip
(ii) {\bf The {\boldmath $0<R<R_n$} case. }

Let ${\mathcal{I}}^{\p}$ denote the set of the multi-indices of the points in ${\mathcal{M}}_n(R)$: 
$${\mathcal{I}}^{\p}={\mathcal{I}}\cup{\mathcal{I}}_{S^1},$$
where ${\mathcal{I}}$ is same as (\ref{f_multi-indices}) and ${\mathcal{I}}_{S^1}$ is the set of the multi-indices of a spider with a folded arm: 
\begin{equation}\label{f_I_S}
{\mathcal{I}}_{S^1}=\left\{\vect{\varepsilon}^{\circ}=(\e^{\circ}_1, \cdots, \e^{\circ}_n)\left|
\begin{array}{l}
\e^{\circ}_k\in\{+,-,\infty\}, \\
\sharp\{j:\e^{\circ}_j=\infty\}=1
\end{array}\right.
\right\}.
\end{equation}
Let ${\mathcal{I}}_m$ $(0\le m\le 2)$, $D_{\vect{\varepsilon}}$ $(\vect{\varepsilon}\in{\mathcal{I}}_2)$, $E_{\vect{\varepsilon}^{\p}}$ $(\vect{\varepsilon}^{\p}\in{\mathcal{I}}_1)$, and $V_{\vect{\varepsilon}^{\p\p}}$ $(\vect{\varepsilon}^{\p\p}\in{\mathcal{I}}_0)$ be given by (\ref{f_I_m}),  (\ref{f_D_e}), (\ref{f_E_e}), and (\ref{f_V_e}) as before. 
Each $D_{\vect{\varepsilon}}$ is homeomorphic to $\stackrel{\circ}D=\textrm{Int} D\setminus\{B_1, \cdots, B_n\}$ (Figure \ref{dilateral_R>1}). 
\begin{figure}[htbp]
\begin{center}
\includegraphics[width=.2\linewidth]{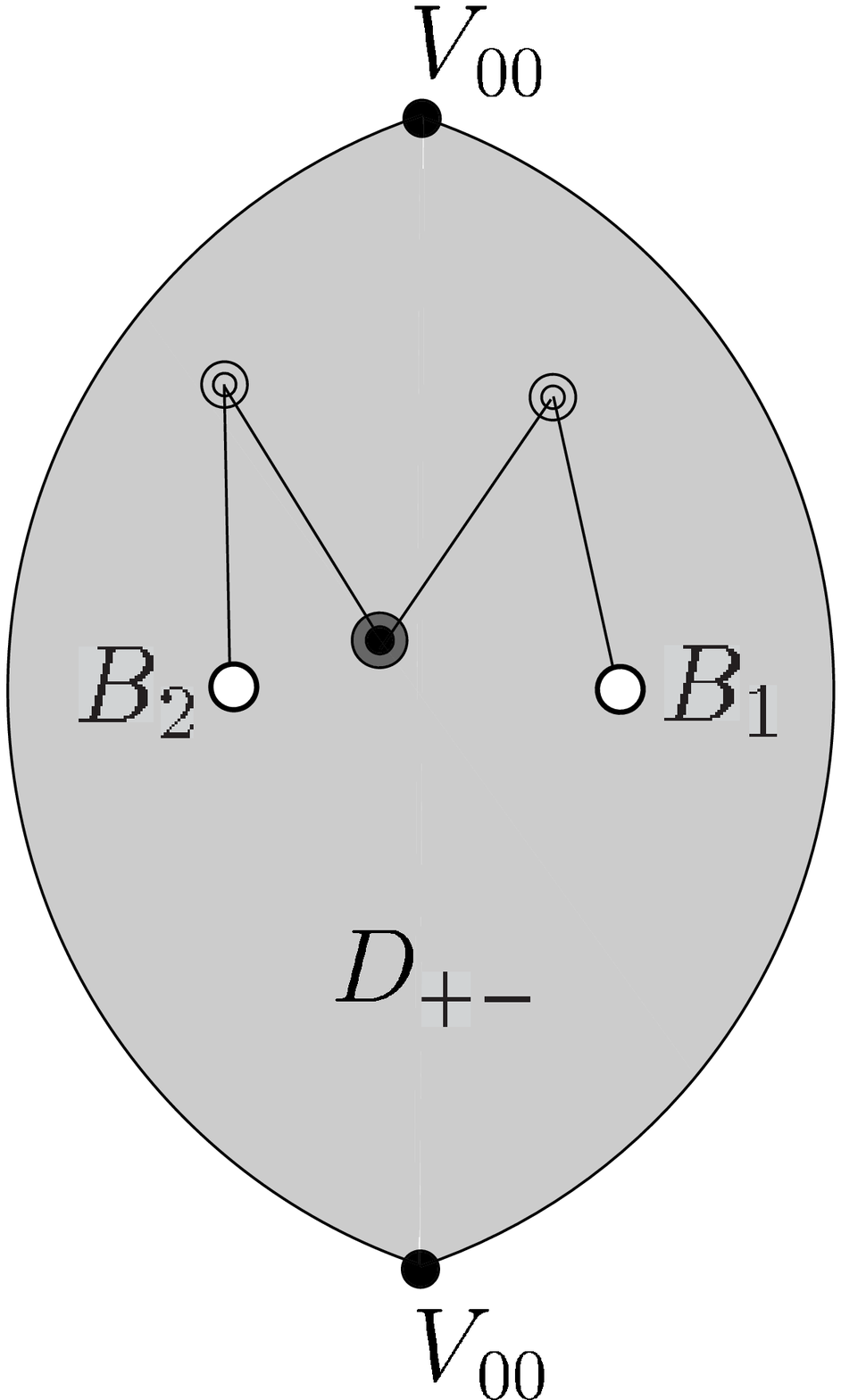}
\caption{The domain ${D}_{+-}$}
\label{dilateral_R>1}
\end{center}
\end{figure}
Put, for $\vect{\varepsilon}^{\circ}\in{\mathcal{I}}_{S^1}$, 
$$
S^1_{\vect{\varepsilon}^{\circ}}=\left\{\vect x\in {\mathcal{M}}_n(R):\,\vect{\varepsilon}(\vect x)=\vect{\varepsilon}^{\circ}\right\}. 
$$
The configuration space ${\mathcal{M}}_n(R)$ can be decomposed as the disjoint union: 
\begin{equation}\label{f_decomp_M_small_R}
{\mathcal{M}}_n(R)=\left(\bigcup_{\vect{\varepsilon}\in{\mathcal{I}}_2}D_{\vect{\varepsilon}}
\cup\bigcup_{\vect{\varepsilon}^{\p}\in{\mathcal{I}}_1}E_{\vect{\varepsilon}^{\p}}
\cup\bigcup_{\vect{\varepsilon}^{\p\p}\in{\mathcal{I}}_0}V_{\vect{\varepsilon}^{\p\p}}\right)
\cup\bigcup_{\vect{\varepsilon}^{\circ}\in{\mathcal{I}}_{S^1}}S^1_{\vect{\varepsilon}^{\circ}}. 
\end{equation}
The first term of the right hand side is homeomorphic to $n2^n$-times punctured orietable surface of genus $1-2^{n-1}+n2^{n-3}$. 
We see how $\cup_{\vect{\varepsilon}^{\circ}\in{\mathcal{I}}_{S^1}}S^1_{\vect{\varepsilon}^{\circ}}$ is glued to it in what follows. 

Let $e_{\vect{\varepsilon}^{\circ}}^{i\theta}$ $(\vect{\varepsilon}^{\circ}\in{\mathcal{I}}_{S^1}, 0\le\theta<2\pi)$ denote a point in $S^1_{\vect{\varepsilon}^{\circ}}\subset{\mathcal{M}}_n(R)$ where the folded arm has angle $\theta$ from the positive direction of the $x$-axis (Figure \ref{eit_in_S^1}). 
\begin{figure}[htbp]
\begin{center}
\includegraphics[width=.3\linewidth]{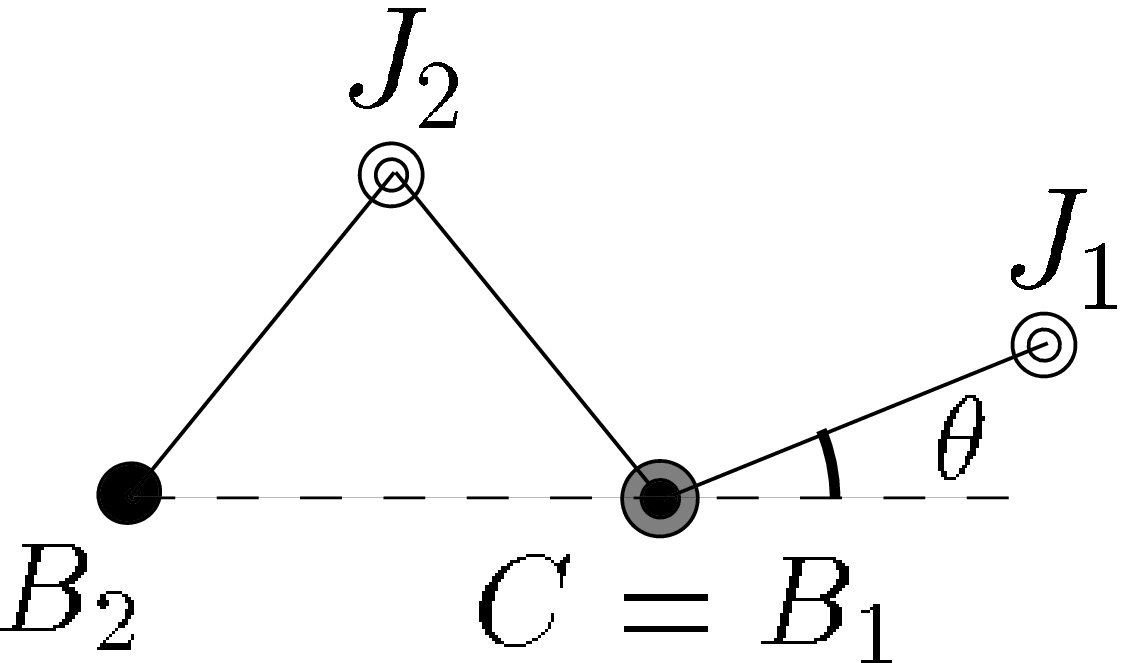}
\caption{$e^{i\theta}_{\infty \,-}\in S^1_{\infty \,-}$. A picture when $n=2$. }
\label{eit_in_S^1}
\end{center}
\end{figure}

Let $\vect{\varepsilon}=(\e_1, \cdots, \e_n)\in{\mathcal{I}}_{2}$. 
Put 
$$
\vect{\varepsilon}_k^{\circ}=(\e_1^{\circ}, \cdots, \e_n^{\circ})\in{\mathcal{I}}_{S^1} \>\>\> \mbox{with}\>\>
\left\{\begin{array}{l}
\e_j^{\circ}=\e_j \>\>\>\mbox{if}\>\>\> j\ne k,\\[1mm]
\e_k^{\circ}=\infty.
\end{array}
\right.
$$
%
Then $e_{\vect{\varepsilon}_k^{\circ}}^{i\theta}\in S^1_{\vect{\varepsilon}^{\circ}}$ is the limit of a sequence of points in $D_{\vect{\varepsilon}}$ whose bodies are located at 
\begin{equation}\label{f_S1_as_limit_from_D}
\left\{
\begin{array}{ll}
B_k+\delta\left(\cos(\theta+\frac{\pi}2), \sin(\theta+\frac{\pi}2)\right) &(\delta>0)\>\>\> \textrm{if} \hspace{0.2cm} \e_k=+ \\[1mm]
B_k+\delta\left(\cos(\theta-\frac{\pi}2), \sin(\theta-\frac{\pi}2)\right) &(\delta>0)\>\>\> \textrm{if} \hspace{0.2cm} \e_k=- 
\end{array}
\right.
\end{equation}
as $\delta$ goes down to $+0$ (Figure \ref{S1_dD}). 
\begin{figure}[htbp]
\begin{center}
\includegraphics[width=1\linewidth]{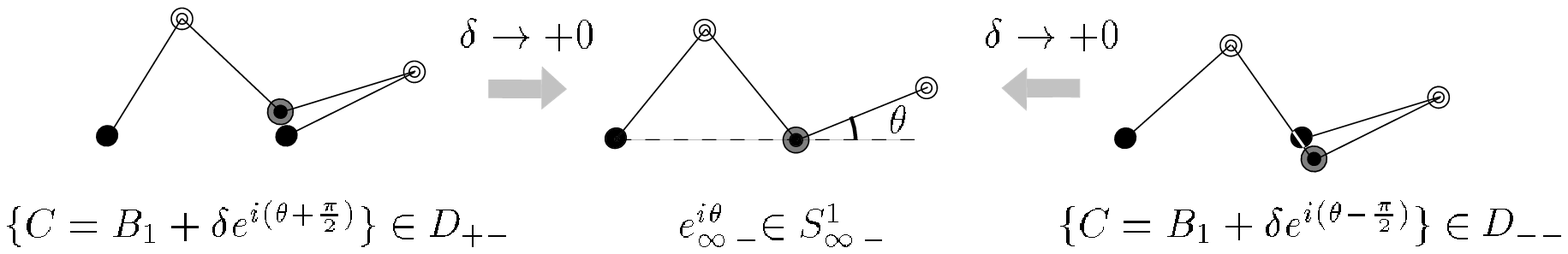}
\caption{}
\label{S1_dD}
\end{center}
\end{figure}

Suppose $\vect{\varepsilon}_k^{\p}$ and $\vect{\varepsilon}_{k,k+1}^{\p\p}$ are given by (\ref{e^p_k-e^pp_k}) as in the previous case. 
Then the closure of $D_{\vect{\varepsilon}}$ in ${\mathcal{M}}_n(R)$ is given by 
$$
\overline{D_{\vect{\varepsilon}}}=D_{\vect{\varepsilon}}
\cup\bigcup_{k=1}^nE_{\vect{\varepsilon}_k^{\p}}
\cup\bigcup_{k=1}^nV_{\vect{\varepsilon}_{k,k+1}^{\p\p}}
\cup\bigcup_{k=1}^nS^1_{\vect{\varepsilon}_k^{\circ}}.
$$
It implies that each $S^1_{\vect{\varepsilon}^{\circ}}$ is contained in exactly two $\overline{D_{\vect{\varepsilon}}}$'s. 
Since ${\mathcal{M}}_n(R)$ is orientable by Corollary \ref{conf_sp=surface} (or by the argument in the Remark below), it means that the decomposition (\ref{f_decomp_M_small_R}) can be considered topologically as attatching $n2^{n-1}$ $1$-handles to $\varSigma_{1-2^{n-1}+n2^{n-3}}$ minus $n2^n$ open discs at the boundary circles (Figure \ref{M_2R_fine}). 
\begin{figure}[htbp]
\begin{center}
\begin{minipage}{.45\linewidth}
\begin{center}
\includegraphics[width=.6\linewidth]{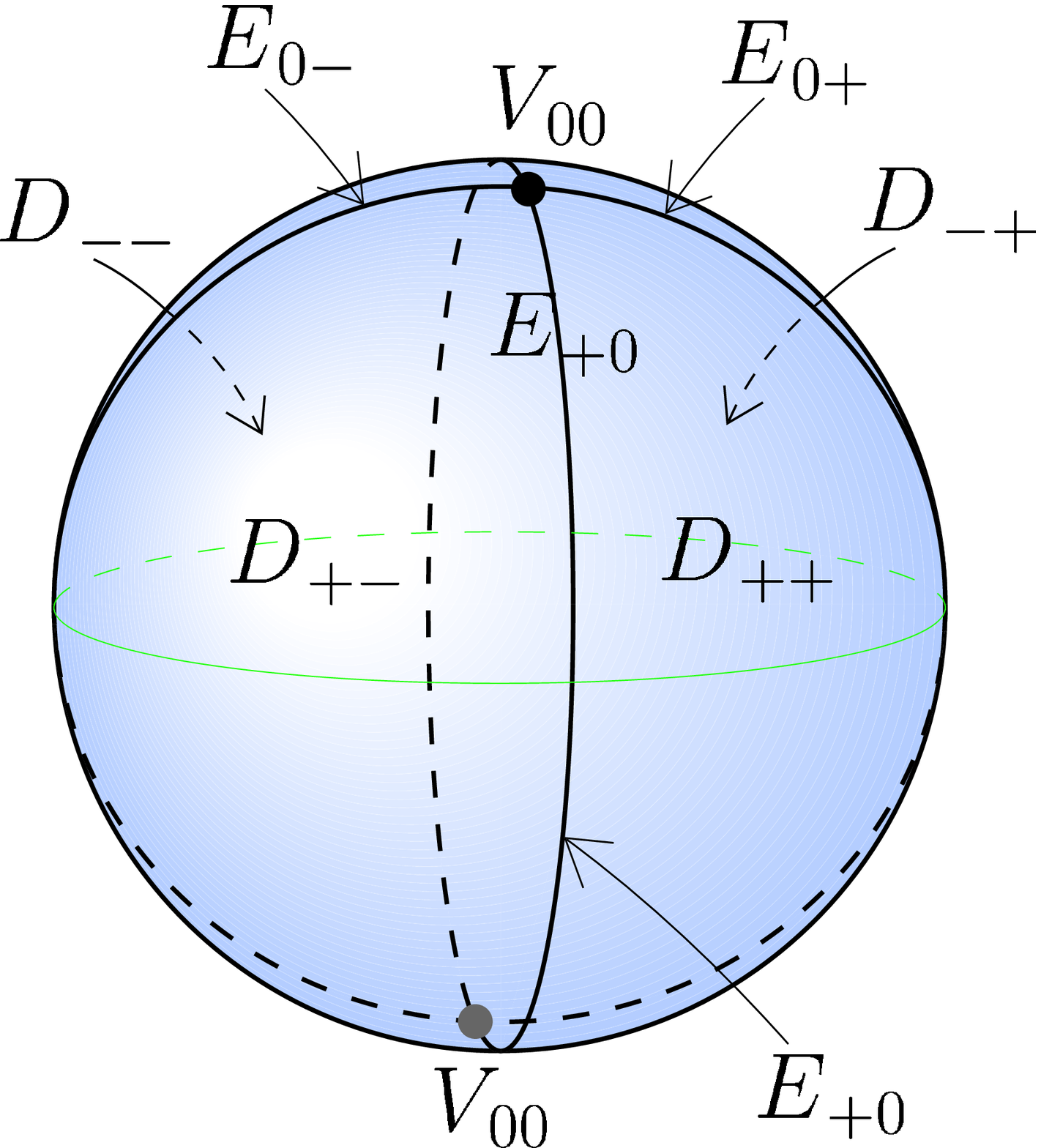}
\caption{A cell decomposition of ${\mathcal{M}}_2(R)$ when $1<R<2$}
\label{M_2Rb_fine}
\end{center}
\end{minipage}
\hskip 0.4cm
\begin{minipage}{.45\linewidth}
\begin{center}
\includegraphics[width=.6\linewidth]{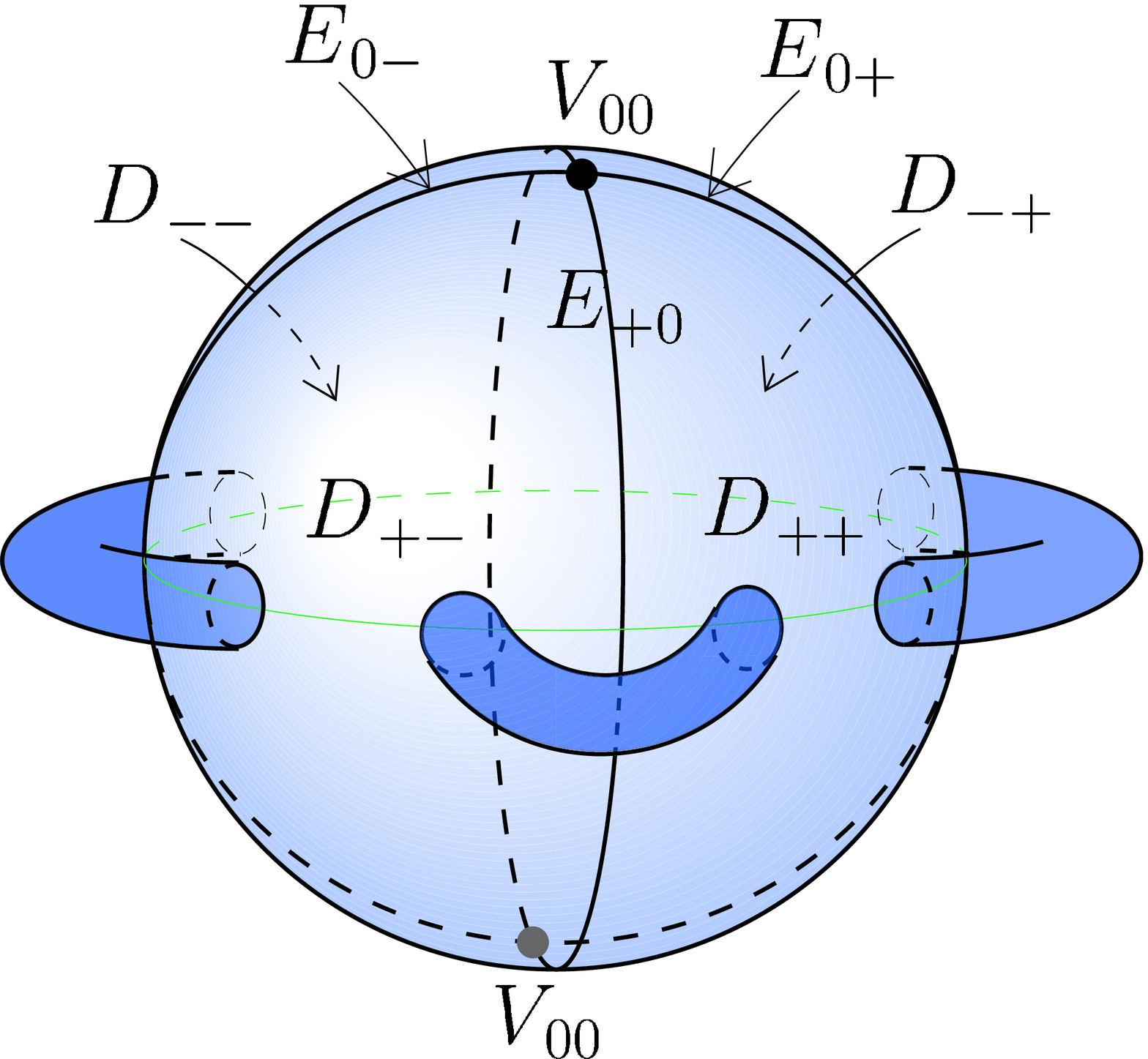}
\caption{${\mathcal{M}}_2(R^{\p})$ $(0<R^{\p}<1)$ can be obtained from ${\mathcal{M}}_2(R)$ $(1<R<2)$ by ataching four $1$-handles. }
\label{M_2R_fine}
\end{center}
\end{minipage}
\end{center}
\end{figure}
\end{proofofproposition}

\medskip
\begin{remark}
The above cell decomposition and cut-and-paste type argument give an alternative proof of Theorem \ref{thm_M_n} in topological category without Propositions \ref{rank_of_Jacobian} and \ref{connectedness}. 

\smallskip
(i) {\bf The {\boldmath $R_n<R<2$} case. }

The formula (\ref{f_closure_D_e}) implies that each edge $E_{\vect{\varepsilon}^{\p}}$ is contained in exactly two faces. 
It means that ${\mathcal{M}}_n(R)$ consisits of a union of closed surfaces. 

A point in $D_{\vect{\varepsilon}}$ and another point in $D_{\vect{\varepsilon}^{\p}}$ can be joined by a path which passes over $\sharp\{j:\e_j\ne\e^{\p}_j\}$ edges. 
Since ${\mathcal{M}}_n(R)$ is a union of the closures of $D_{\vect{\varepsilon}}$'s, it implies that ${\mathcal{M}}_n(R)$ is connected. 

As was mentioned before, $D_{\vect{\varepsilon}}$ can be identified with a copy of $\textrm{Int}D$ by the position of the body. 
Suppose the orientation of $D_{\vect{\varepsilon}}$ is given through this identification by that of $\textrm{Int}D$ multiplied by 
$$
(-1)^{m(\vect{\varepsilon})}, \hspace{0.2cm}\textrm{where}\hspace{0.2cm} m(\vect{\varepsilon})=\sharp\{j:\e_j=-\}. 
$$
Two faces $\overline{D_{\vect{\varepsilon}}}$ and $\overline{D_{\vect{\varepsilon}^{\p}}}$ are adjacent if and only if 
$$
\sharp\{j:\e_j\ne\e^{\p}_j\}=1.
$$
In this case they meet at an edge which corresponds to the same edge of $D$. 
The two orientations of $D_{\vect{\varepsilon}}$ and $D_{\vect{\varepsilon}^{\p}}$ fit at this edge. 
Therefore, ${\mathcal{M}}_n(R)$ is orientable. 

\medskip
(ii) {\bf The {\boldmath $0<R<R_n$} case. }

Let $\varphi_{\vect{\varepsilon}}:\stackrel{\circ}{D}\to D_{\vect{\varepsilon}}$ be the homeomorphism given by the position of the body. 

Define the compactification $\overline{D}$ of $\stackrel{\circ}{D}$ by 
$$
\overline{D}=\stackrel{\circ}{D}\cup\,\partial_1D\cup\partial_2D, 
$$
where $\partial_1D$ is the union of the $n$ edges of $D$, i.e. the boundary of $D\subset\mathbb{R}^2$ in the usual sense, and $\partial_2D$ is the union of $n$ $S^1$'s, 
where a point $e^{i\theta}$ in the $k$-th $S^1$ is the limit of a point $B_k+\delta(\cos\theta, \sin\theta)$ as $\delta$ goes down to $+0$. 
Then $\varphi_{\vect{\varepsilon}}$ can be extended to 
$$\overline{\varphi}_{\vect{\varepsilon}}:\overline{D}\to \overline{D_{\vect{\varepsilon}}}.$$

Assume the orientation of $D_{\vect{\varepsilon}}$ is given in the same way as in the previous case. 
Suppose two faces $\overline{D_{\vect{\varepsilon}}}$ and $\overline{D_{\vect{\varepsilon}^{\p}}}$ meet at some $S^1_{\vect{\varepsilon}^{\circ}}$ which is the image of the $k$-th $S^1$ in $\partial_2D$ by $\overline{\varphi}_{\vect{\varepsilon}}$ and $\overline{\varphi}_{\vect{\varepsilon}^{\p}}$. 
Then (\ref{f_S1_as_limit_from_D}) implies that the restriction of $\overline{\varphi}_{\vect{\varepsilon}^{\p}}{}^{-1}\circ\overline{\varphi}_{\vect{\varepsilon}}$ to the $k$-th $S^1$ is the antipodal map, which is isotopic to the identity. 
Since $\sharp\{j:\e_j\ne\e^{\p}_j\}=1$, $\overline{D_{\vect{\varepsilon}}}$ and $\overline{D_{\vect{\varepsilon}^{\p}}}$ inherit opposite orientations from $D$ through $\overline{\varphi}_{\vect{\varepsilon}}$ and $\overline{\varphi}_{\vect{\varepsilon}^{\p}}$. 
It means that, through $\overline{\varphi}_{\vect{\varepsilon}}$ and $\overline{\varphi}_{\vect{\varepsilon}^{\p}}$, $\overline{D}$ is glued to its copy with the opposite orientation at the $k$-th $S^1$ in $\partial_2D$ by the identity map. 
The two orientations fit at the $S^1$. 
Therefore, ${\mathcal{M}}_n(R)$ is orientable. 
\end{remark}

\subsection{Morse Theoretical method to determine the genus}

\begin{figure}[htbp]
\begin{center}
\includegraphics[width=.3\linewidth]{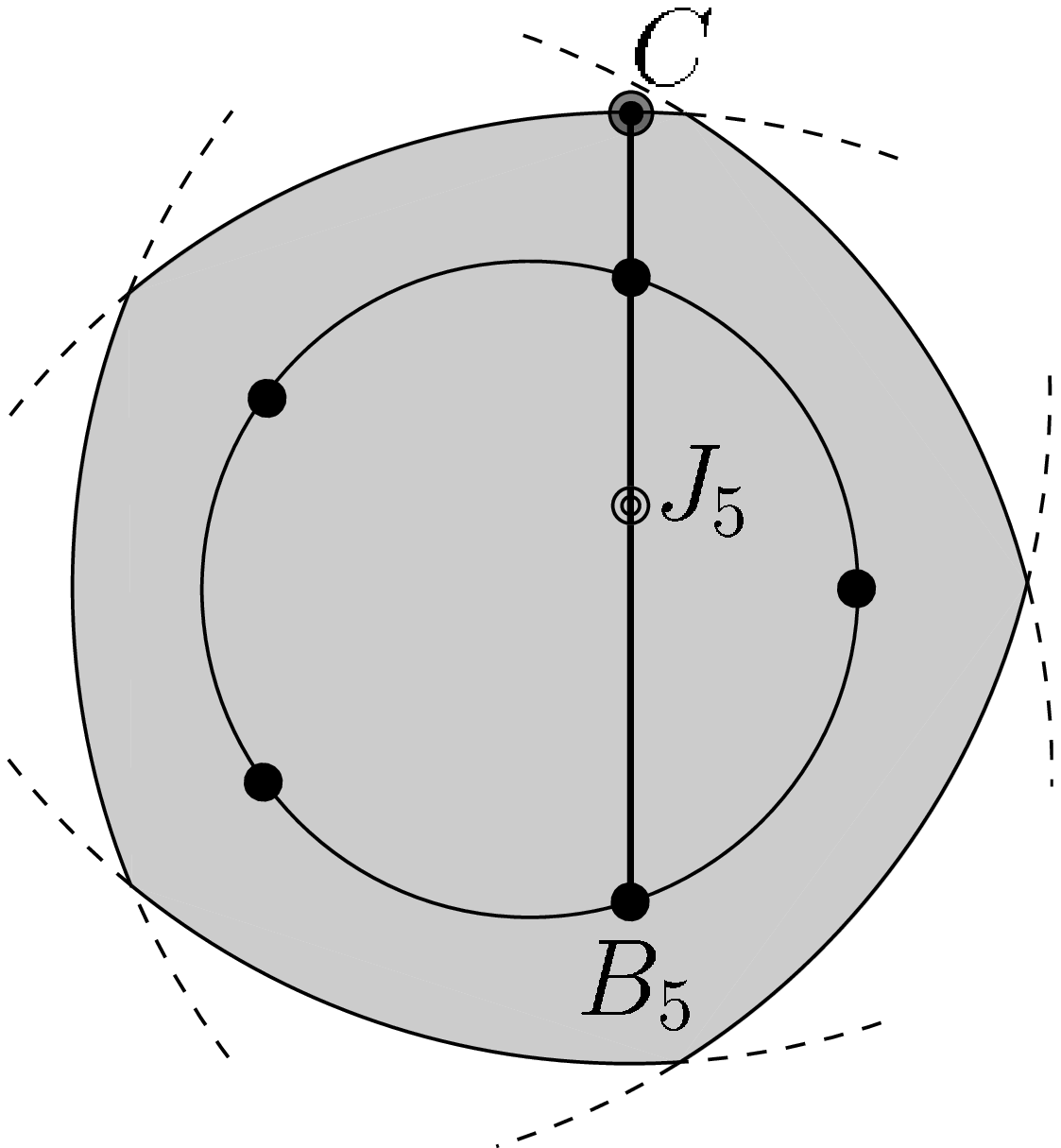}
\caption{The fourth arm is stretched out. The other arms are not drawn.}
\label{pentagon}
\end{center}
\end{figure}

We prove Proposition \ref{Morse_theoretical_proof_of_the_genus}. 
We assume that $R$ satisfies $0<R<R_n$ or $R_n<R<2$ in this Subsecton. 
\begin{lemma}\label{critical_points}
Suppose $R$ satisfies $0<R<R_n$ or $R_n<R<2$. 
A point $\vect x\in{\mathcal{M}}_n(R)$ is a critical point of $\psi(\vect x)=y$ if and only if one of the following conditions is satisfied: 
\begin{enumerate}
\item[{\rm (i)}] Two adjacent arms are stretced out inward {\rm (Figures \ref{6armsV_index1-case2}, \ref{6armsV_index1-case1}, and \ref{6armsV_index1-case3-part1})}. 
\item[{\rm (ii)}] Exactly one arm, say, the $k$-th, is folded, which is parallel to the $y$-axis {\rm (Figures \ref{6arms-ll} and \ref{6arms^ll_at_2})}. 
The body is located at $B_k$. 
It can occur only when $0<R<R_n$. 
\item[{\rm (iii)}] Exactly one arm, say, the $k$-th, is stretced out, which is parallel to the $y$-axis. 
It can occur if and only if $n$ is odd, $0<R<R_n$, and $B_k$ is either the highest or the lowest {\rm (Figure \ref{pentagon})}. 
\end{enumerate}
\end{lemma}
\begin{proof}
Recall that 
$${\mathcal{M}}_n(R)=\left\{\vect x:f_1(\vect x)=\cdots=f_{2n}(\vect x)=0\right\}.$$ 
Let $\mbox{\sl Span}\,\langle \partial f_k(\vect x)\rangle$ denote the linear subspace of $\mathbb{R}^{2n+2}$ spanned by $\partial f_1(\vect x), \cdots,$ $\partial f_{2n}(\vect x)$. 
It is codimension $2$. 
The tangent space $T_{\vect x}{\mathcal{M}}_n(R)$ of ${\mathcal{M}}_n(R)$ at $\vect x$ is equal to the orthogonal complement $\{\mbox{\sl Span}\,\langle \partial f_k(\vect x)\rangle\}^{\perp}$ of $\mbox{\sl Span}\,\langle \partial f_k(\vect x)\rangle$. 
A point $\vect x=(x,y,p_1,q_1,\cdots, p_n,q_n)\in{\mathcal{M}}_n(R)$ is a critical point of $\psi(\vect x)=y$ if and only if its differential vector $\partial \psi(\vect x)=(\vect e_2, \vect 0, \cdots, \vect 0)$ is orthogonal to $T_{\vect x}{\mathcal{M}}_n(R)=\{\mbox{\sl Span}\,\langle \partial f_k(\vect x)\rangle\}^{\perp}$, which occurs if and only if $(\vect e_2, \vect 0, \cdots, \vect 0)$ is contained in $\mbox{\sl Span}\,\langle \partial f_k(\vect x)\rangle$. 

Suppose 
$$\sum_{k=1}^{2n}c_k\partial f_k=(\vect e_2, \vect 0, \cdots, \vect 0) \hspace{0.3cm}.$$ 
The differential vectors $\partial f_k(\vect x)$ of $f_k$ are give by 
$$\begin{array}{rclcccr}
\partial f_1(\vect x)\!\!&\!\!=\!\!&\!\!2(\vect a_1, &-\vect a_1, &\vect 0, &\cdots, &\vect 0),\\
\partial f_2(\vect x)\!\!&\!\!=\!\!&\!\!2(\vect 0, &\vect b_1, &\vect 0, &\cdots, &\vect 0),\\
&\vdots & \\
\partial f_{2n-1}(\vect x)\!\!&\!\!=\!\!&\!\!2(\vect a_n, &\vect 0, & \cdots, &-\vect a_n, &\vect 0),\\
\partial f_{2n}(\vect x)\!\!&\!\!=\!\!&\!\!2(\vect 0, & \vect 0, & \cdots, &\vect 0&\vect b_n), 
\end{array}$$
where $\vect a_k=\overrightarrow{J_kC}$ and $\vect b_k=\overrightarrow{B_kJ_k}$.
At least one of $c_{2k-1}$'s is not equal to $0$. 
If $c_{2k-1}\ne0$ then $c_{2k}=\pm c_{2k-1}$ and $\vect a_k=\pm \vect b_k$, i.e. the $k$-th arm is either stretced out or folded. 

If it is folded then Lemma \ref{d_n} implies that there are no other non-zero $c_{2j-1}$'s as $R\ne R_n$. 
It is the case (ii). 

If there are no folded arms then Lemma \ref{stretched-out_arms} implies the number of stretced-out arms is either one or two. 
The latter case corresponds to the case (i). 

Suppose there is exactly one stretced-out arm, say the $k$-th arm. 
Then $\vect a_k=\vect b_k=\pm \vect e_2$, i.e. the $k$-th arm is parallel to the $y$-axis. 
If $B_k$ is not the unique highest (or the lowest) point, then at least one of $|CB_{k-1}|$ and $|CB_{k+1}|$ is bigger than $2$, which is a contradiction. 
Therefore, $n$ cannot be even. 
Suppose $n$ is odd. 
A line segment of length $2$ and parallel to the $y$-axis which starts from the highest (or the lowest) $B_k$ is contained in the the curved $n$-gon $D$ if and only if $0<R<R_n$. 
This is the case (iii). 
\end{proof}

\bigskip
\begin{proofofproposition}\ref{Morse_theoretical_proof_of_the_genus}. 
Lemma \ref{critical_points} implies the Proposition under the assumption that $\psi$ is a Morse function on ${\mathcal{M}}_n(R)$, which will be proved in Proposition \ref{non-degenerate} below. 

Since ${\mathcal{M}}_n(R)$ is $2$-dimensional, the index of a critical point $\vect x$ of $\psi$ is $2$ if $\psi(\vect x)$ is local maximum, $0$ if $\psi(\vect x)$ is local minimum, and $1$ otherwise. 

\smallskip
(1) Suppose $R_n<R<2$. 
All the critical points are of type (i) of Lemma \ref{critical_points}, i.e. with two adjacent arms stretced-out inward. 
\begin{figure}[htbp]
\begin{center}
\begin{minipage}{.3\linewidth}
\begin{center}
\includegraphics[width=\linewidth]{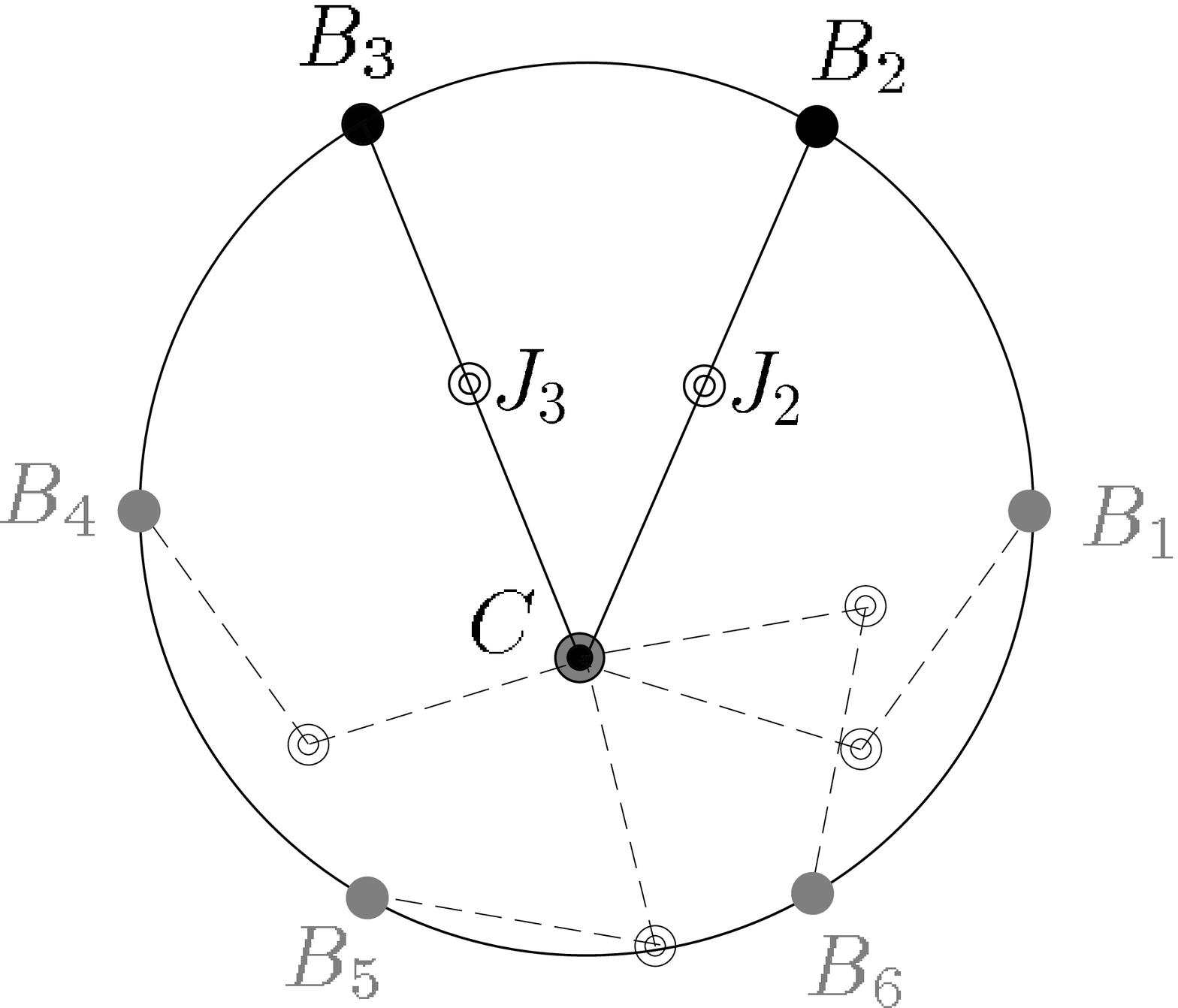}
\caption{Type (i) critical point of index $0$}
\label{6armsV_index1-case2}
\end{center}
\end{minipage}
\hskip 0.1cm
\begin{minipage}{.3\linewidth}
\begin{center}
\includegraphics[width=\linewidth]{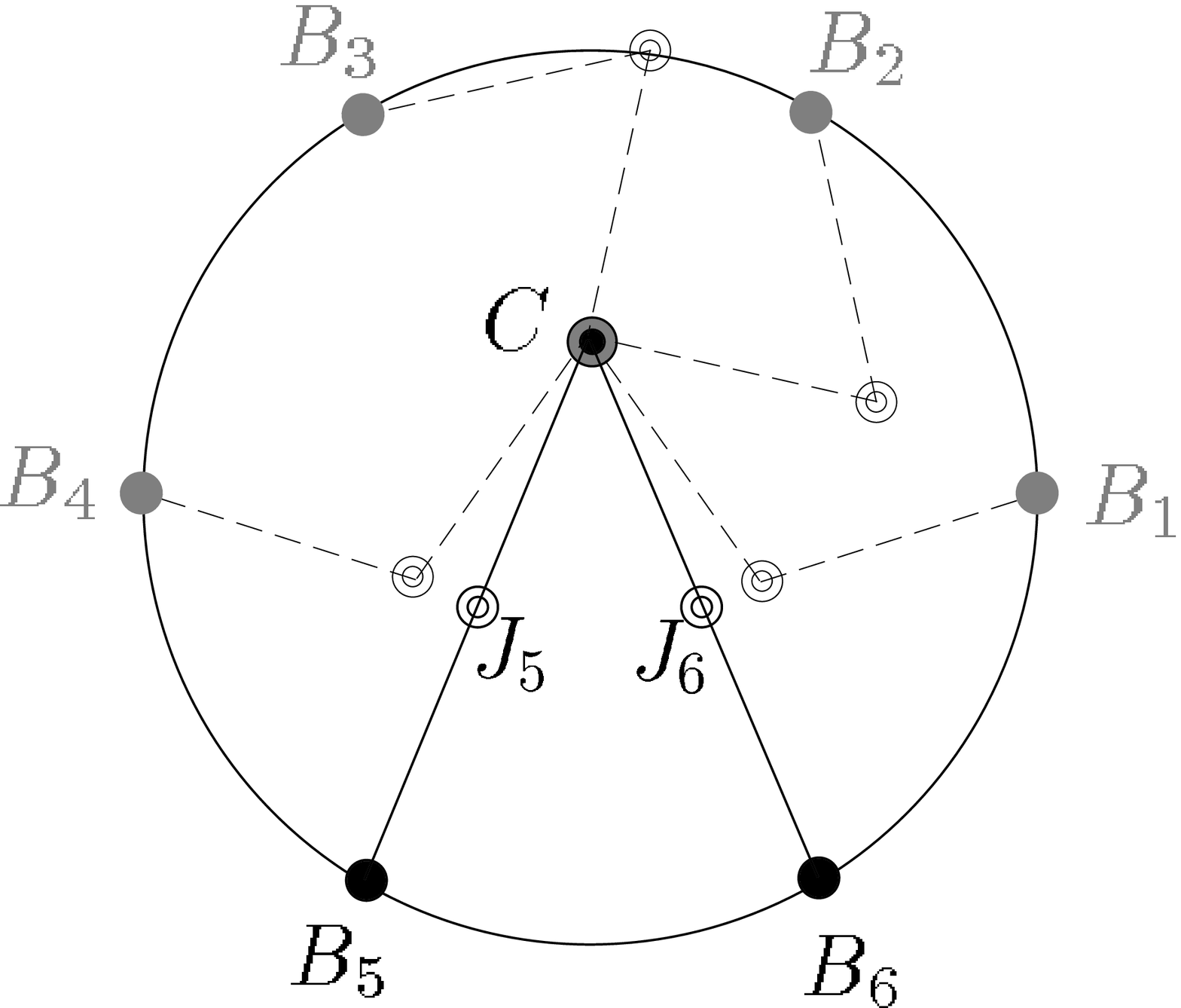}
\caption{Type (i) critical point of index $2$}
\label{6armsV_index1-case1}
\end{center}
\end{minipage}
\hskip 0.1cm
\begin{minipage}{.3\linewidth}
\begin{center}
\includegraphics[width=\linewidth]{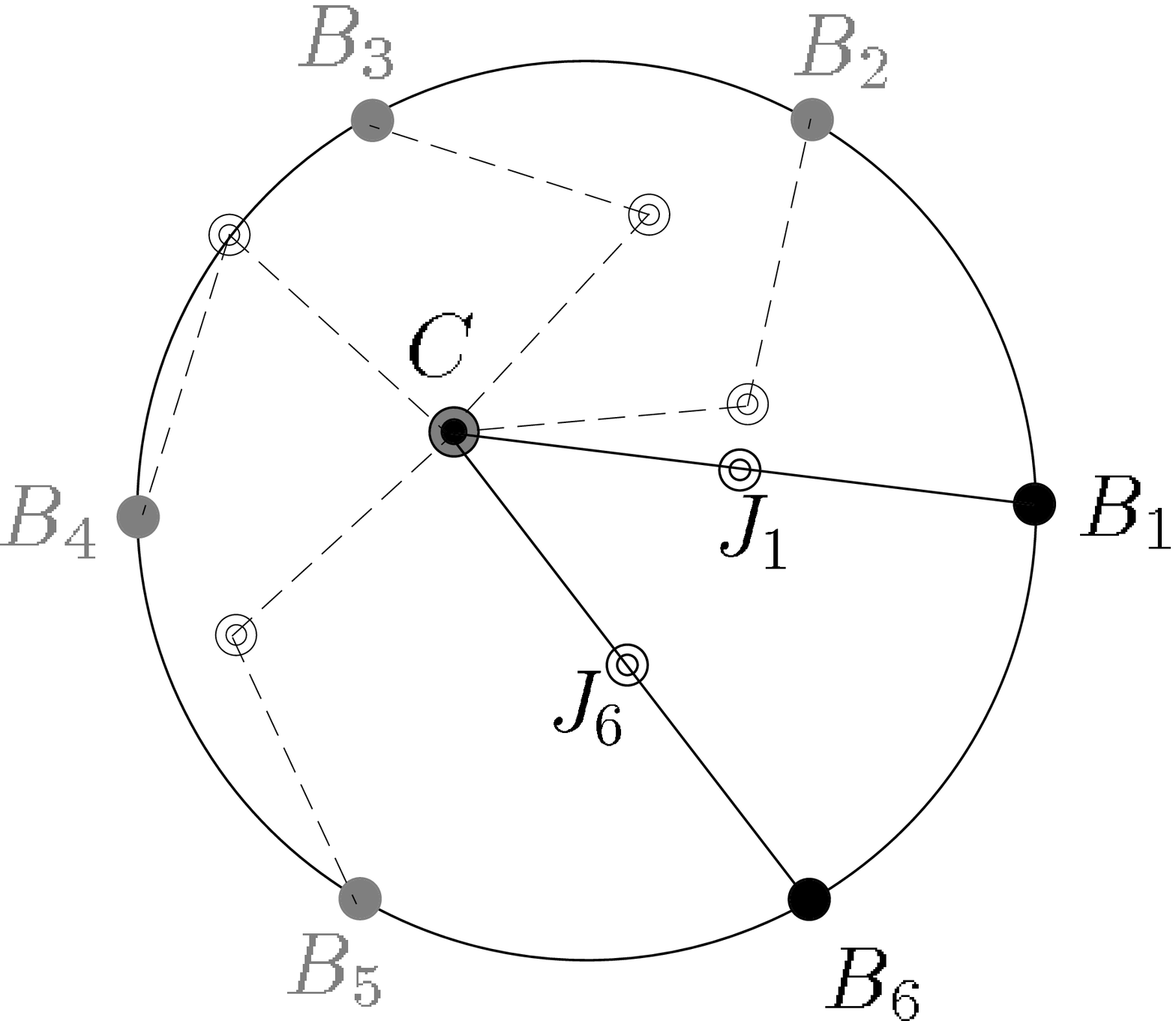}
\caption{Type (i) critical point of index $1$}
\label{6armsV_index1-case3-part1}
\end{center}
\end{minipage}
\end{center}
\end{figure}
There are $2^{n-2}$ critical points of index $0$ (Figure \ref{6armsV_index1-case2}), same number of critical points of index $2$ (Figure \ref{6armsV_index1-case1}), and $(n-2)2^{n-2}$ critical points of index $1$ (Figure \ref{6armsV_index1-case3-part1}). 

\smallskip
(2) Suppose $0<R<R_n$ and $n$ is even. 
All the critical points are either of type (i) or type (ii) of Lemma \ref{critical_points}. 
The number and the indices of type (i) critical points are same as in the previous case (1). 
There are $n\cdot 2\cdot 2^{n-1}=n2^n$ critical points of type (ii). 
They all have index $1$ (Figures \ref{6arms-ll} and \ref{6arms^ll_at_2}). 

\begin{figure}[htbp]
\begin{center}
\begin{minipage}{.45\linewidth}
\begin{center}
\includegraphics[width=.8\linewidth]{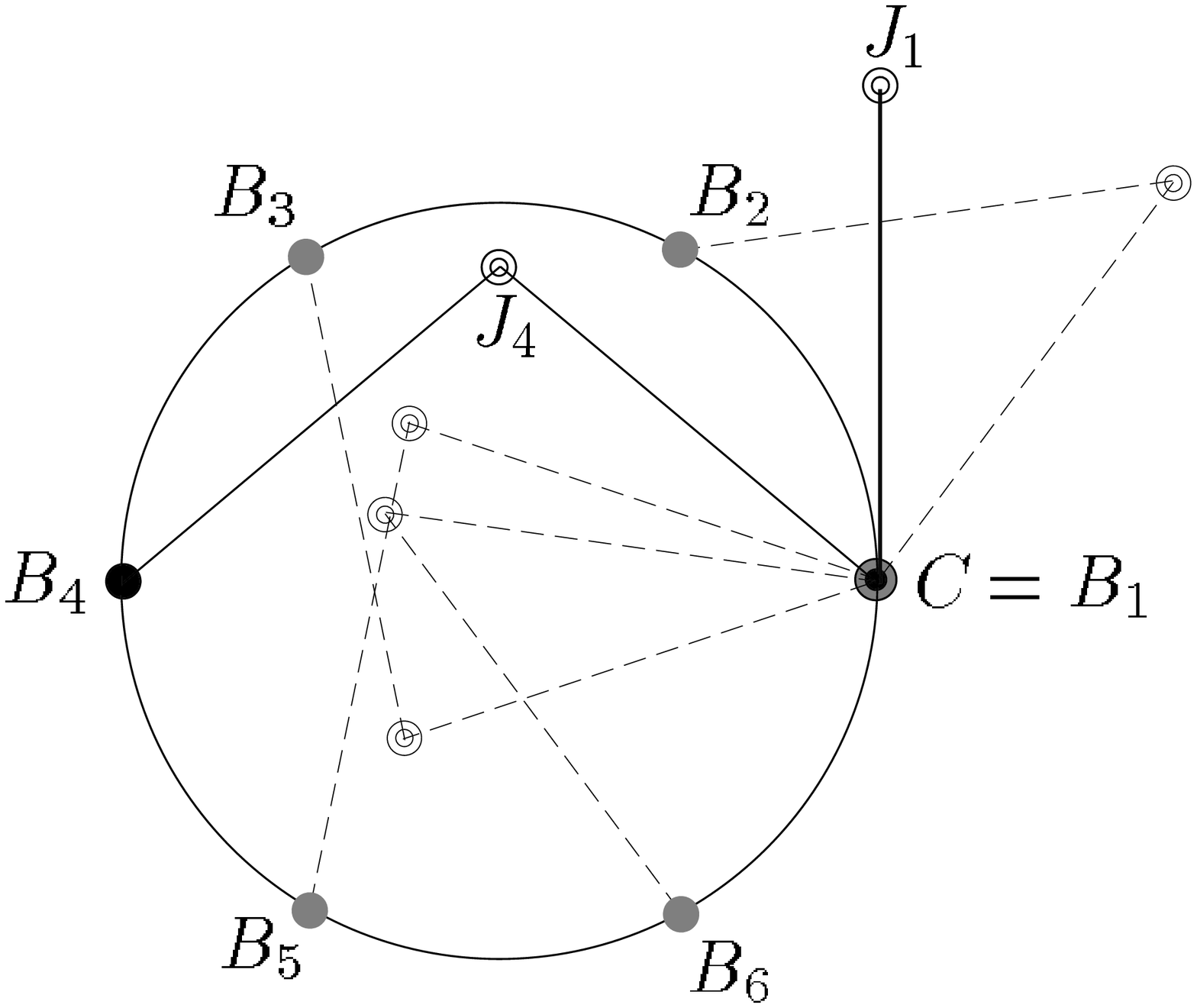}
\caption{Type (ii) critical point of index $1$}
\label{6arms-ll}
\end{center}
\end{minipage}
\hskip 0.4cm
\begin{minipage}{.45\linewidth}
\begin{center}
\includegraphics[width=.72\linewidth]{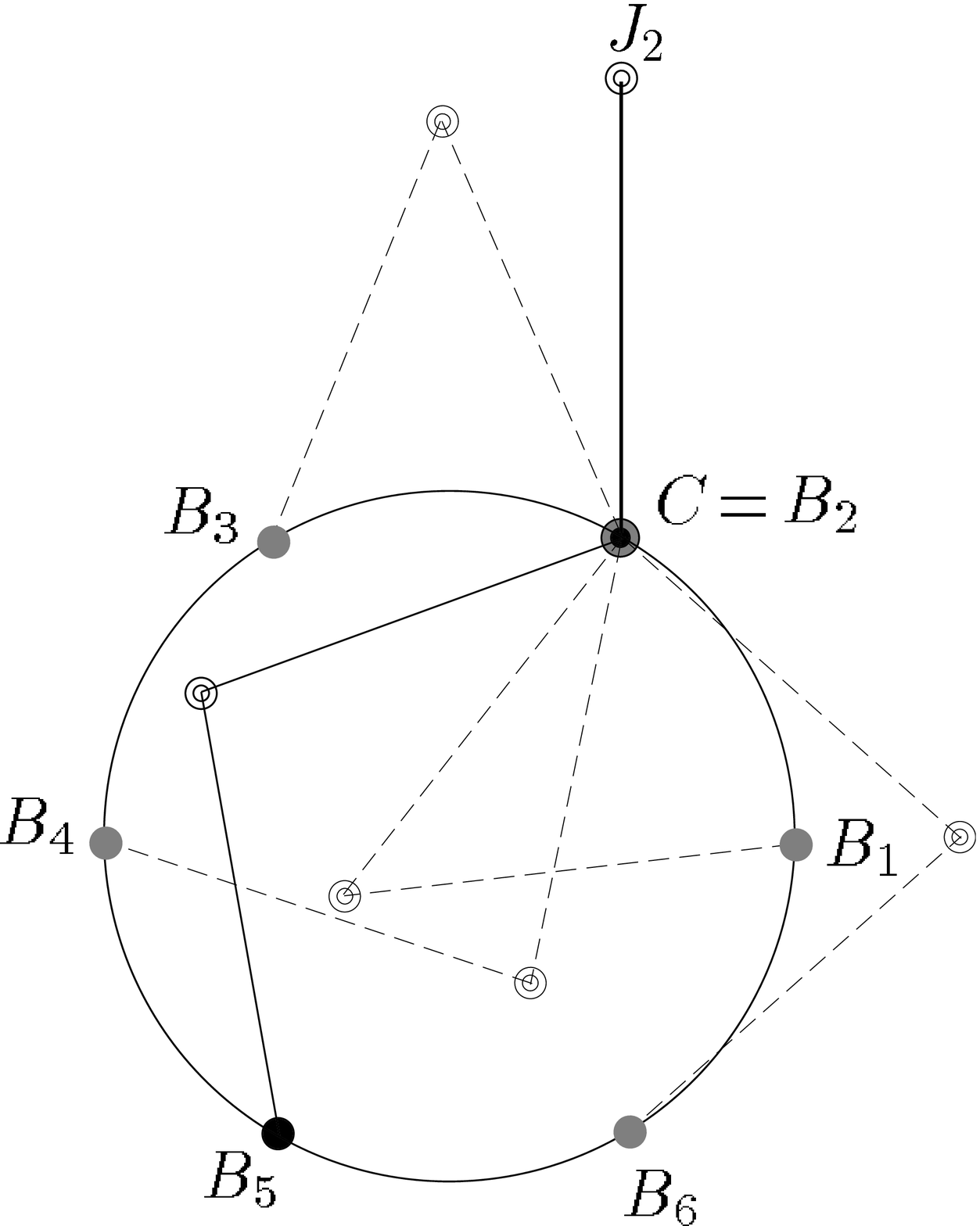}
\caption{Type (ii) critical point of index $1$}
\label{6arms^ll_at_2}
\end{center}
\end{minipage}
\end{center}
\end{figure}

\smallskip
(3) Suppose $0<R<R_n$ and $n$ is odd. 
The three types, (i), (ii), and (iii) of Lemma \ref{critical_points} appear as critical points. 
Unlike in the previous two cases, all the critical points of type (i) have index $1$ since any vertex of the curved $n$-gon $D$ cannot be a highest or a lowest point in $D$ (Figure \ref{3-gon_R_small-i}). 
\begin{figure}[htbp]
\begin{center}
\begin{minipage}{.45\linewidth}
\begin{center}
\includegraphics[width=.6\linewidth]{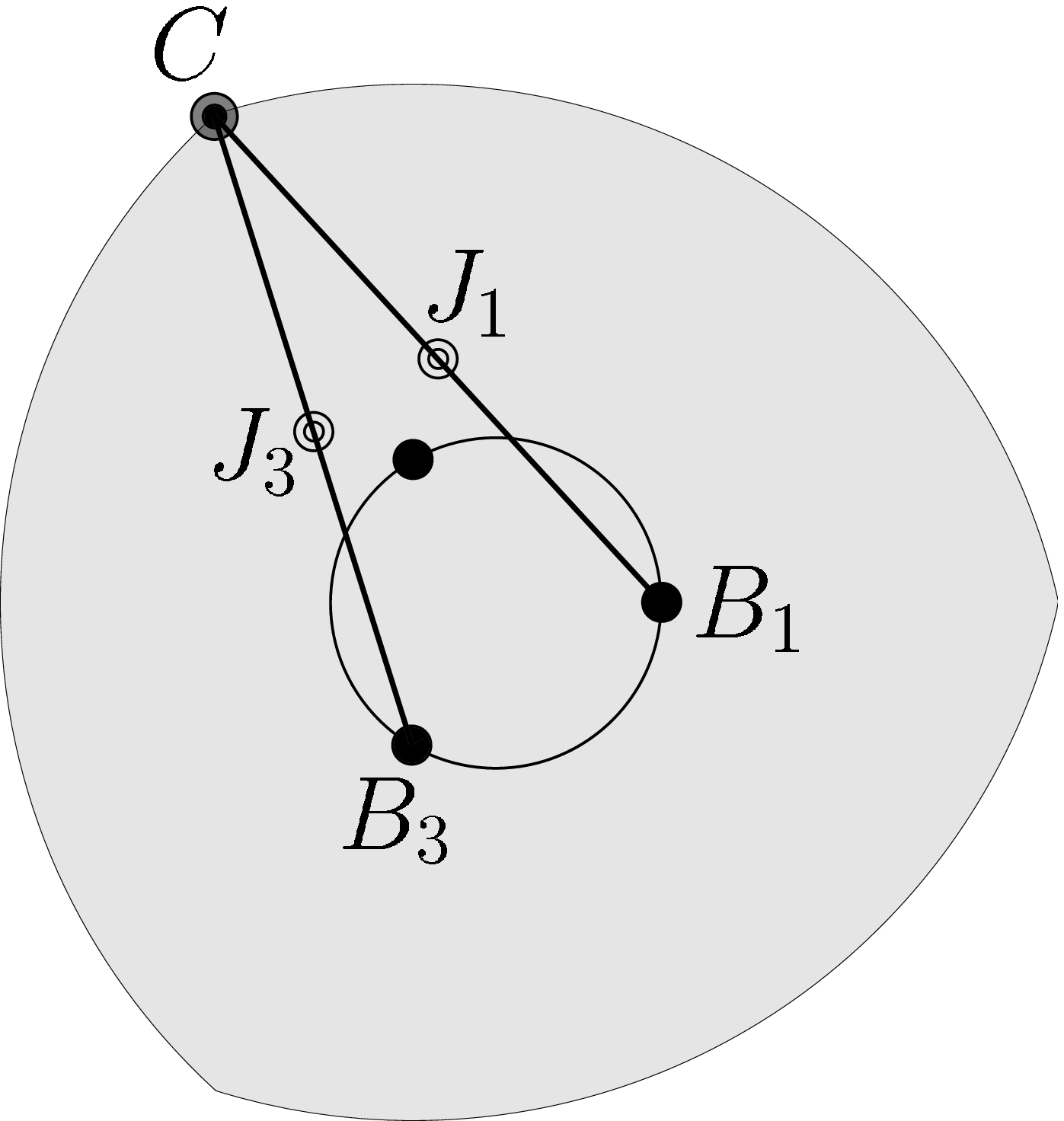}
\caption{Type (i) critical point with index $1$. The second arm is not drawn.}
\label{3-gon_R_small-i}
\end{center}
\end{minipage}
\hskip 0.4cm
\begin{minipage}{.45\linewidth}
\begin{center}
\includegraphics[width=.6\linewidth]{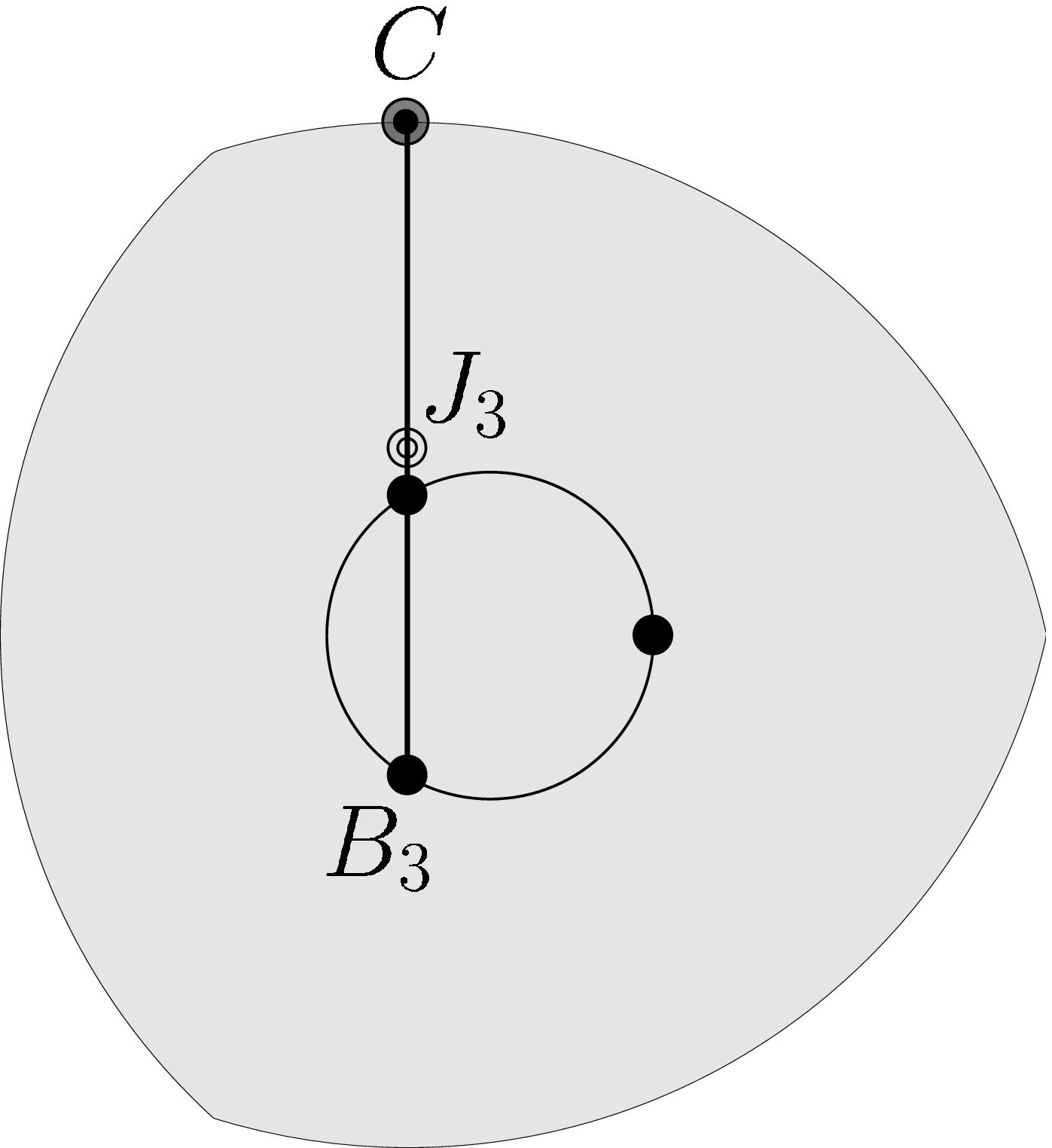}
\caption{Type (iii) critical point with index $2$. Only the third arm is drawn.}
\label{3-gon_R_small-ii}
\end{center}
\end{minipage}
\end{center}
\end{figure}
Therefore, summing up critical points of types (i) and (ii), we can find $n2^{n-2}+n2^n$ critical points of index $1$. 
A critical point of type (iii) has index $0$ or $2$, each case has $2^{n-1}$ critical points. 
\end{proofofproposition}
\begin{proposition}\label{non-degenerate}
The critical points of $\psi$ which are given in {\rm Lemma \ref{critical_points}} are non-degenerate. 
\end{proposition}
\begin{proof}
We can give local coordinates around a critical point using the stretched-out arms or the folded arm, since they determine the position of the body, which determines the position of all the other bended arms in turn. 

The proof is devided into three cases according to the types of critical points. 

\medskip
{\bf Type (i) critical points of Lemma \ref{critical_points}}. 

\begin{lemma}\label{parallel_to_x_y_axis}
Each of any pair of inward stretched-out adjacent arms is not parallel either to the $x$-axis or to the $y$-axis if $0<R<R_n$ or $R_n<R<2$. 
\end{lemma}
\begin{proof}
Suppose the $k$-th and $(k+1)$-th arms are stretched out. 
Let $\rho$ be the angle $(0\le\rho<\pi)$ of one of the two stretched-out arms from the $x$-axis. 
Then $\rho$ or $\rho+\pi$ belongs to $\left(\frac{2(k-1)}{n}\pi, \frac{2k}{n}\pi\right)$, where the two boundary values correspond to the case of $R=2$. 
It tends to $\frac{2k-1}{n}\pi$ as $R$ approaches $+0$. 
If $n$ is odd, it tends to $\frac{4k-2\pm 1}{2n}\pi$ as $R$ approaches $R_n$. 
\begin{figure}[htbp]
\begin{center}
\includegraphics[width=.7\linewidth]{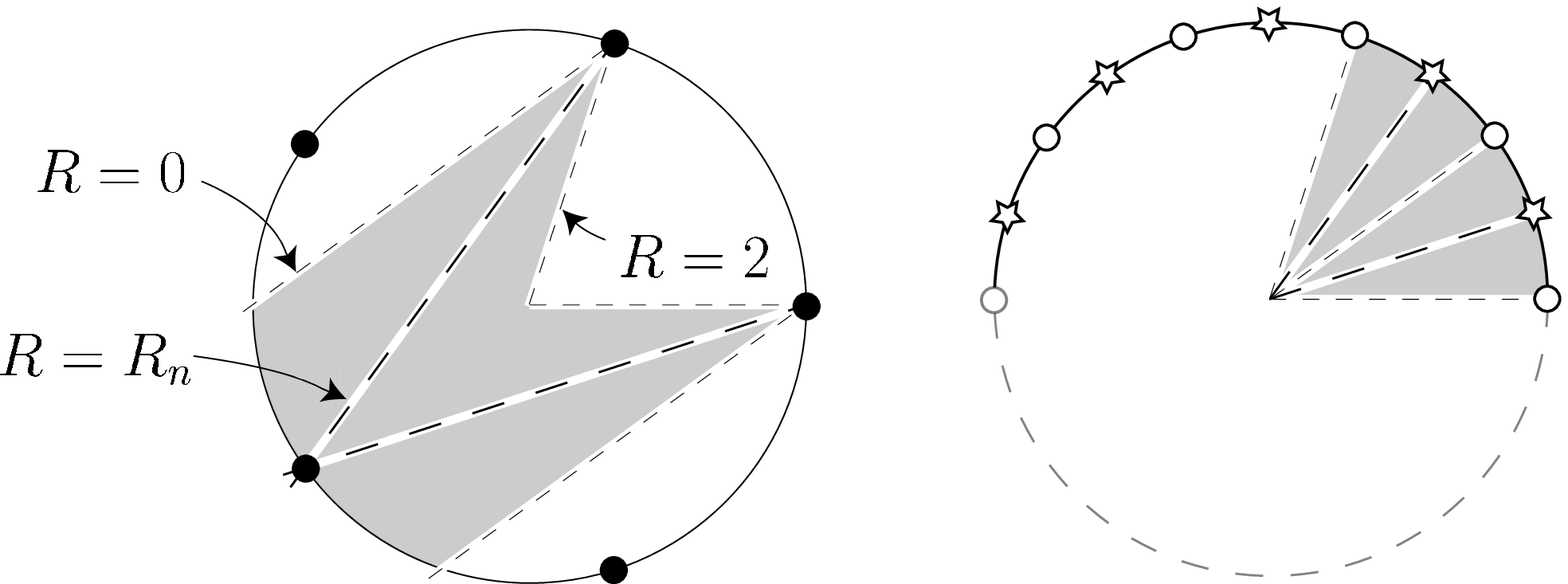}
\caption{The excluded angles}
\label{range_angle_5}
\end{center}
\end{figure}

Let $\mathcal{A}$ be the set of the such angles $\rho$ $(0\le\rho<\pi)$ as $R$ varies in $0<R<R_n$ and $R_n<R<2$. 
Then $\mathcal{A}$ misses $2n$ points in $[0,\pi)$; $n$ points corresponding to the case of $R=0$ or $R=2$, and another $n$ points to the case of $R=R_n$. 
The former are $0, \frac{\pi}n, \cdots, \frac{n-1}n\pi$, which include $0$, and furthermore, $\frac{\pi}2$ if $n$ is even. 
If $n$ is odd, the latter are $\frac{\pi}{2n}, \frac{3\pi}{2n}, \cdots, \frac{2n-1}{2n}\pi$, which include $\frac{\pi}2$ (Figure \ref{range_angle_5}). 
\end{proof}

\medskip
Suppose the $k$-th and $(k+1)$-th arms are stretched out. 
Since neither is parallel to the $x$-axis we have: 
\begin{lemma}\label{local_coordinates}
The $x$-coordinates $p_k$ and $p_{k+1}$ of the two joints $J_k$ and $J_{k+1}$ can serve as local coordinates. 
\end{lemma}
\begin{proof}
The inverse function theorem implies that a pair of functions $\xi(\vect x)=p_1$ and $\eta(\vect x)=p_2$ serves as a system of local coordinates of ${\mathcal{M}}_n(R)=F^{-1}(\vect 0)$ in a neighbourhood of a point $\vect x\in{\mathcal{M}}_n(R)$ if and only if the matrix 
\begin{equation}\label{matrix_p_1_p_2}
\left(
\begin{array}{c}
\partial f_1(\vect x)\\
\partial f_2(\vect x)\\
\vdots\\
\partial f_{2k-1}(\vect x)\\
\partial f_{2k}(\vect x)\\
\partial f_{2(k+1)-1}(\vect x)\\
\partial f_{2(k+1)}(\vect x)\\
\vdots\\
\partial f_{2n-1}(\vect x)\\
\partial f_{2n}(\vect x)\\
\partial \xi(\vect x)\\
\partial \eta(\vect x)
\end{array}
\right)=
\left(
\begin{array}{ccccc}
\> 2\vect a_1 \>&\> -2\vect a_1 \> & &&\\
& 2\vect b_1 & &&\\
\vdots&&\ddots&&\\
\> 2\vect a_k \>&&\> -2\vect a_k \> && \\
&& 2\vect b_k && \\
\> 2\vect a_{k+1} \>&&&\> -2\vect a_{k+1} \> & \\
&&& 2\vect b_{k+1} & \\
\vdots&&&\ddots&\\
\> 2\vect a_{n} \>&&&&\> -2\vect a_{n} \>  \\
&&&& 2\vect b_{n}  \\
&& \vect e_1 && \\
&&& \vect e_1 & 
\end{array}
\right)
\end{equation}
is non-singular. 

Suppose 
$$\sum_{j=1}^{2n}c_j\partial f_j(\vect x)+d_k\partial \xi(\vect x)+d_{k+1}\partial \eta(\vect x)=\vect 0.$$ 
We have $c_{2j-1}=c_{2j}=0$ if $j\ne k,k+1$ since $\vect a_j\ne\pm\vect b_j$. 
We have $c_{2k-1}=c_{2(k+1)-1}=0$ since $\vect a_k\ne\pm\vect a_{k+1}$. 
Since Lemma \ref{parallel_to_x_y_axis} implies that $\vect b_k\ne\pm\vect e_1\ne\vect b_{k+1}$ we have $c_{2k}=d_k=c_{2(k+1)}=d_{k+1}=0$, which completes the proof. 
\end{proof}

\medskip
Recall $(u_j, v_j)$, $(p_j, q_j)$, and $(x, y)$ denote the coordinates of $B_j$ given by (\ref{coord_B_k}), $J_j$, and $C$ respectively. 
We show that the Hessian of $\psi$ at $\vect x$ does not vanish: 
$$
\det H(\psi)(\vect x)=\left|
\begin{array}{cc}
\displaystyle \frac{\partial^2 y}{\partial p_k{}^2} & \displaystyle \frac{\partial^2 y}{\partial p_k\partial p_{k+1}}\\[4mm]
\displaystyle \frac{\partial^2 y}{\partial p_{k+1}\partial p_k} & \displaystyle \frac{\partial^2 y}{\partial p_{k+1}{}^2}
\end{array}
\right|\ne0. 
$$
We have 
\begin{eqnarray}\label{}
&&(x-p_k)^2+(y-q_k)^2-1\equiv 0, \label{f_d_1}\\[0.5mm]
&&(p_k-u_k)^2+(q_k-v_k)^2-1\equiv 0, \label{f_d_2}\\[0.5mm]
&&(x-p_{k+1})^2+(y-q_{k+1})^2-1\equiv 0, \label{f_d_3}\\[0.5mm]
&&(p_{k+1}-u_{k+1})^2+(q_{k+1}-v_{k+1})^2-1\equiv 0. \label{f_d_4}
\end{eqnarray}
By differentiating (\ref{f_d_2}) and (\ref{f_d_4}) by $p_k$ and $p_{k+1}$ we have 
\begin{equation}\label{q_k-p_k}
\begin{array}{c}
\displaystyle 
\frac{\partial q_k}{\partial p_k}=-\frac{p_k-u_k}{q_k-v_k},\>\>\>
\frac{\partial q_{k+1}}{\partial p_{k+1}}=-\frac{p_{k+1}-u_{k+1}}{q_{k+1}-v_{k+1}},\>\>\>
\frac{\partial q_k}{\partial p_{k+1}}=\frac{\partial q_{k+1}}{\partial p_k}=0, \\[4mm]
\displaystyle \frac{\partial^2 q_k}{\partial p_k{}^2}=-\frac1{(q_k-v_k)^3}, \>\>\>
\frac{\partial^2 q_{k+1}}{\partial p_{k+1}{}^2}=-\frac1{(q_{k+1}-v_{k+1})^3}. 
\end{array}
\end{equation}
By differentiating (\ref{f_d_1}) and (\ref{f_d_3}) by $p_k$ and by applying (\ref{q_k-p_k}) we get 
$$
\left\{
\begin{array}{rcl}
\displaystyle (x-p_k)\frac{\partial x}{\partial p_k}+(y-q_k)\frac{\partial y}{\partial p_k}&=&\displaystyle (x-p_k)+(y-q_k)\frac{\partial q_k}{\partial p_k}\\[4mm]
&=&\displaystyle  \frac{(x-p_k)(q_k-v_k)-(y-q_k)(p_k-u_k)}{q_k-v_k},\\
\displaystyle (x-p_{k+1})\frac{\partial x}{\partial p_k}+(y-q_{k+1})\frac{\partial y}{\partial p_k}&=&0, 
\end{array}
\right.
$$
which implies 
$$
\left(\begin{array}{c}
\displaystyle \frac{\partial x}{\partial p_k}\\[4mm]
\displaystyle \frac{\partial y}{\partial p_k}
\end{array}\right)
=\frac1{q_k-v_k}
\left|
\begin{array}{cc}
x-p_k\,&\,y-q_k\\
p_k-u_k&q_j-v_k
\end{array}
\right|
\left(\begin{array}{cc}
x-p_k\,&\,y-q_k\\
x-p_{k+1}\,&\,y-q_{k+1}\\
\end{array}\right)^{-1}
\left(\begin{array}{c}
1\\[1mm]
0
\end{array}\right). 
$$
Similarly we have 
$$
\begin{array}{l}
\left(\begin{array}{c}
\displaystyle \frac{\partial x}{\partial p_{k+1}}\\[4mm]
\displaystyle \frac{\partial y}{\partial p_{k+1}}
\end{array}\right)
=\displaystyle \frac1{q_{k+1}-v_{k+1}}
\left|
\begin{array}{cc}
x-p_{k+1}\,&\,y-q_{k+1}\\
p_{k+1}-u_{k+1}&q_j-v_{k+1}
\end{array}
\right|\\[8mm]
\hspace{4.4cm}
\left(\begin{array}{cc}
x-p_k\,&\,y-q_k\\
x-p_{k+1}\,&\,y-q_{k+1}\\
\end{array}\right)^{-1}
\left(\begin{array}{c}
0\\[1mm]
1
\end{array}\right). 
\end{array}
$$
Since $\vect a_k=\vect b_k$ and $\vect a_{k+1}=\vect b_{k+1}$ we have 
$$
\begin{array}{l}
\left|
\begin{array}{cc}
x-p_k\,&\,y-q_k\\
p_k-u_k&q_j-v_k
\end{array}
\right|
=\left|
\begin{array}{c}
\,\vect a_k\,\\[1mm]
\vect b_k
\end{array}
\right|
=0, \\[5mm]
\left|
\begin{array}{cc}
x-p_{k+1}\,&\,y-q_{k+1}\\
p_{k+1}-u_{k+1}&q_j-v_{k+1}
\end{array}
\right|
=\left|
\begin{array}{c}
\,\vect a_{k+1}\,\\[1mm]
\vect b_{k+1}
\end{array}
\right|
=0, 
\end{array}
$$
which imply
\begin{equation}\label{X_k=0}
\frac{\partial x}{\partial p_k}=\frac{\partial y}{\partial p_k}=\frac{\partial x}{\partial p_{k+1}}=\frac{\partial y}{\partial p_{k+1}}=0.
\end{equation}

By differentiating (\ref{f_d_1}) and (\ref{f_d_3}) by $p_k$ twice we get 
$$
\left\{
\begin{array}{rcl}
\displaystyle (x-p_k)\frac{\partial^2 x}{\partial p_k{}^2}+(y-q_k)\frac{\partial^2 y}{\partial p_k{}^2}
&=&\displaystyle -\left(\frac{\partial x}{\partial p_k}-1\right)^2-\left(\frac{\partial y}{\partial p_k}-\frac{\partial q_k}{\partial p_k}\right)^2\\[4mm]
&&\displaystyle \hspace{2.5cm} +(y-q_k)\frac{\partial^2 q_k}{\partial p_k{}^2},
\\[4mm]
\displaystyle (x-p_{k+1})\frac{\partial^2 x}{\partial p_k{}^2}+(y-q_{k+1})\frac{\partial^2 y}{\partial p_k{}^2}&=&\displaystyle -\left(\frac{\partial x}{\partial p_k}\right)^2-\left(\frac{\partial y}{\partial p_k}\right)^2,
\end{array}
\right.
$$
which implies 
%
\begin{equation}\label{f_differentials1}
\begin{array}{rcl}
\left(\begin{array}{c}
\displaystyle \frac{\partial^2 x}{\partial p_k{}^2}\\[3mm]
\displaystyle \frac{\partial^2 y}{\partial p_k{}^2}
\end{array}\right)
\!\!&\!\!=\!\!&\!\!\displaystyle 
\left(\begin{array}{cc}
x-p_k\,&\,y-q_k\\
&\\
x-p_{k+1}\,&\,y-q_{k+1}\\
\end{array}\right)^{-1}{}\cdot\\[4mm]
&&\displaystyle 
\left(\begin{array}{c}
\displaystyle -\left(\frac{\partial x}{\partial p_k}-1\right)^2-\left(\frac{\partial y}{\partial p_k}-\frac{\partial q_k}{\partial p_k}\right)^2+(y-q_k)\frac{\partial^2 q_k}{\partial p_k{}^2}\\[4mm]
\displaystyle -\left(\frac{\partial x}{\partial p_k}\right)^2-\left(\frac{\partial y}{\partial p_k}\right)^2
\end{array}\right). 
\end{array} 
\end{equation}
Similarly we have 
%
\begin{equation}\label{f_differentials2}
\begin{array}{l}
\left(\begin{array}{c}
\displaystyle \frac{\partial^2 x}{\partial p_{k+1}{}^2}\\[3mm]
\displaystyle \frac{\partial^2 y}{\partial p_{k+1}{}^2}
\end{array}\right)
=\left(\begin{array}{cc}
x-p_k\,&\,y-q_k\\
&\\
x-p_{k+1}\,&\,y-q_{k+1}\\
\end{array}\right)^{-1}{}\cdot\\[8mm]
\hspace{1.4cm}
\displaystyle 
\left(\begin{array}{c}
\displaystyle -\left(\frac{\partial x}{\partial p_{k+1}}\right)^2-\left(\frac{\partial y}{\partial p_{k+1}}\right)^2\\[4mm]
\displaystyle -\left(\frac{\partial x}{\partial p_{k+1}}-1\right)^2-\left(\frac{\partial y}{\partial p_{k+1}}-\frac{\partial q_{k+1}}{\partial p_{k+1}}\right)^2+(y-q_{k+1})\frac{\partial^2 q_{k+1}}{\partial p_{k+1}{}^2}\\[4mm]
\end{array}\right), \\[8mm]
\end{array}
\end{equation}
%
and
$$
\begin{array}{rcl}
\left(\begin{array}{c}
\displaystyle \frac{\partial^2 x}{\partial p_k\partial p_{k+1}}\\[3mm]
\displaystyle \frac{\partial^2 y}{\partial p_k\partial p_{k+1}}
\end{array}\right)
\!\!&\!\!=\!\!&\!\!\displaystyle 
\left(\begin{array}{cc}
x-p_k\,&\,y-q_k\\
&\\
x-p_{k+1}\,&\,y-q_{k+1}\\
\end{array}\right)^{-1}{}\cdot\\[8mm]
&&\left(\begin{array}{c}
\displaystyle -\left(\frac{\partial x}{\partial p_k}-1\right)\frac{\partial x}{\partial p_{k+1}}-\left(\frac{\partial y}{\partial p_k}-\frac{\partial^2 q_k}{\partial p_k{}^2}\right)\frac{\partial y}{\partial p_{k+1}}\\[4mm]
\displaystyle -\left(\frac{\partial x}{\partial p_{k+1}}-1\right)\frac{\partial x}{\partial p_k}-\left(\frac{\partial y}{\partial p_{k+1}}-\frac{\partial^2 q_{k+1}}{\partial p_{k+1}{}^2}\right)\frac{\partial y}{\partial p_k}
\end{array}\right).
\end{array}
$$
%
%
%
%
\smallskip
Since $\displaystyle \frac{\partial x}{\partial p_k}=\frac{\partial y}{\partial p_k}=\frac{\partial x}{\partial p_{k+1}}=\frac{\partial y}{\partial p_{k+1}}=0$ the above formula implies 
\begin{equation}\label{d^2y_dp_kdp_k+1}
\displaystyle \frac{\partial^2 y}{\partial p_k\partial p_{k+1}}=0. 
\end{equation}
Let $\theta$ and $\theta^{\p}$ be the angles of $\vect a_k=\vect b_k$ and $\vect a_{k+1}=\vect b_{k+1}$ from the $x$-axis respectively. 
Then 
$$
\displaystyle 
\left(\begin{array}{cc}
x-p_k\,&\,y-q_k\\
x-p_{k+1}\,&\,y-q_{k+1}\\
\end{array}\right)
=\left(\begin{array}{cc}
\cos\theta\,&\,\sin\theta\\
\cos\theta^{\p}\,&\,\sin\theta^{\p}\\
\end{array}\right),
$$
and (\ref{q_k-p_k}) implies 
$$
\frac{\partial q_k}{\partial p_k}=-\frac{\cos\theta}{\sin\theta}, \>\>
\frac{\partial^2 q_k}{\partial p_k{}^2}=-\frac{1}{\sin^3\theta}, \>\>
\frac{\partial q_{k+1}}{\partial p_{k+1}}=-\frac{\cos\theta^{\p}}{\sin\theta^{\p}}, \>\>
\frac{\partial^2 q_{k+1}}{\partial p_{k+1}{}^2}=-\frac{1}{\sin^3\theta^{\p}}.
$$
Therefore, (\ref{f_differentials1}) and (\ref{f_differentials2}) imply 
$$
\begin{array}{rcl}
\left(\begin{array}{c}
\displaystyle \frac{\partial^2 x}{\partial p_k{}^2}\\[3mm]
\displaystyle \frac{\partial^2 y}{\partial p_k{}^2}
\end{array}\right)
\!\!&\!\!=\!\!&\!\!\displaystyle \frac1{\sin(\theta^{\p}-\theta)}
\left(\begin{array}{cc}
\sin\theta^{\p}&-\sin\theta\\[1mm]
-\cos\theta^{\p}&\cos\theta
\end{array}\right)
\left(
\begin{array}{c}
\displaystyle -\,\frac2{\,\sin^2\theta\,}\\[4mm]
0
\end{array}
\right),
\\[9mm]
\left(\begin{array}{c}
\displaystyle \frac{\partial^2 x}{\partial p_{k+1}{}^2}\\[3mm]
\displaystyle \frac{\partial^2 y}{\partial p_{k+1}{}^2}
\end{array}\right)
\!\!&\!\!=\!\!&\!\!\displaystyle \frac1{\sin(\theta^{\p}-\theta)}
\left(\begin{array}{cc}
\sin\theta^{\p}&-\sin\theta\\[1mm]
-\cos\theta^{\p}&\cos\theta
\end{array}\right)
\left(
\begin{array}{c}
0\\
\displaystyle -\,\frac2{\,\sin^2\theta^{\p}\,}
\end{array}
\right),
\end{array}
$$
which imply 
\begin{equation}\label{d^2y_dp^2}
\begin{array}{rcl}
\displaystyle \frac{\partial^2 y}{\partial p_k{}^2}
\!\!&\!\!=\!\!&\!\! \displaystyle 
\frac{2\cos\theta^{\p}}{\sin(\theta^{\p}-\theta)\sin^2\theta},\\[4mm]
\displaystyle \frac{\partial^2 y}{\partial p_{k+1}{}^2}
\!\!&\!\!=\!\!&\!\! \displaystyle 
-\frac{2\cos\theta}{\sin(\theta^{\p}-\theta)\sin^2\theta^{\p}}.
\end{array}
\end{equation}
Since $\theta, \theta^{\p}\ne 0, \pi, \frac{\pi}2, \frac{3\pi}2$ by Lemma \ref{parallel_to_x_y_axis}, and $\theta^{\p}-\theta\ne 0,\pm\pi$, (\ref{d^2y_dp_kdp_k+1}) and (\ref{d^2y_dp^2}) imply that the Hessian is not equal to $0$. 

\medskip
{\bf Type (ii) critical points of Lemma \ref{critical_points}}. 

\smallskip
Suppose the $k$-th arm is folded. 
Then $\vect a_k=-\vect b_k=\pm\vect e_2$. 
A slight modification of Lemma \ref{local_coordinates} implies that $x$ and $p_k$ can serve as local coordinates. 

By differentiating 
$$
(p_k-u_k)^2+(q_k-v_k)^2-1\equiv 0
$$
by $x$ and $p_k$ we get 
\begin{equation}\label{f_pd_ii}
\frac{\partial q_k}{\partial p_k}=-\frac{p_k-u_k}{q_k-v_k},\>\>\>
\frac{\partial^2 q_k}{\partial p_k{}^2}=-\frac1{(q_k-v_k)^3}, \>\>\>
\frac{\partial q_k}{\partial x}\equiv 0. 
\end{equation}
By differentiating 
$$
(x-p_k)^2+(y-q_k)^2-1\equiv 0
$$
by $x$ and $p_k$, and by applying 
\begin{equation}\label{f_ii_01}
x-p_k=p_k-u_k=0, \>y-q_k=-(q_k-v_k)=\pm 1, 
\end{equation}
and (\ref{f_pd_ii}), we obtain $\displaystyle \frac{\partial y}{\partial x}=\frac{\partial y}{\partial p_k}=0$ and 
\begin{equation}\label{d2y_dx^2-ii}
\begin{array}{rcl}
\displaystyle \frac{\partial^2y}{\partial x^2}\!\!&\!\!=\!\!&\!\! \displaystyle -\frac1{\,y-q_k\,}=\mp1,\\[4mm]
\displaystyle \frac{\partial^2y}{\partial p_k{}^2}\!\!&\!\!=\!\!&\!\! \displaystyle -\frac1{\,y-q_k\,}-\frac1{\,(q_k-v_k)^3\,}=0,\\[4mm]
\displaystyle \frac{\partial^2y}{\partial x\partial p_k}\!\!&\!\!=\!\!&\!\! \displaystyle \frac1{\,y-q_k\,}=\pm1,
\end{array}
\end{equation}
which implies that the Hessian is equal to $-1$. 

\medskip
{\bf Type (iii) critical points of Lemma \ref{critical_points}}. 

Suppose the $k$-th arm is stretched out. 
Then $\vect a_k=\vect b_k=\pm\vect e_2$. 
The argument goes parallel to the previous case. 
We can take $x$ and $p_k$ as local coordinates. 
What is different from the previous case is that (\ref{f_ii_01}) is replaced by 
\begin{equation}\label{f_iii_01}
x-p_k=p_k-u_k=0, \>y-q_k=q_k-v_k=\pm 1, 
\end{equation}
and hence (\ref{d2y_dx^2-ii}) is replaced by 
\begin{equation}\label{d2y_dx^2-iii}
\begin{array}{rcl}
\displaystyle \frac{\partial^2y}{\partial x^2}\!\!&\!\!=\!\!&\!\! \displaystyle -\frac1{\,y-q_k\,}=\mp1,\\[4mm]
\displaystyle \frac{\partial^2y}{\partial p_k{}^2}\!\!&\!\!=\!\!&\!\! \displaystyle -\frac1{\,y-q_k\,}-\frac1{\,(q_k-v_k)^3\,}=\mp2,\\[4mm]
\displaystyle \frac{\partial^2y}{\partial x\partial p_k}\!\!&\!\!=\!\!&\!\! \displaystyle \frac1{\,y-q_k\,}=\pm1,
\end{array}
\end{equation}
which implies that the determinant of the Hessian matrix is equal to $1$. 

This completes the proof of Proposition \ref{non-degenerate}. 
\end{proof} 

We remark that the Proposition can also be proved by expressing $\psi$ explicitly in terms of $p_k$ and $p_{k+1}$ (or other coordinates). 
The calculation becomes much more complicated. 

\section{Proof for the singular case}\label{sec_sing}

In this section we study the configuration space ${\mathcal{M}}_n(R)$ of the spiders with $n$ arms of radius $R$ when it is not a smooth surface. 

\smallskip
\begin{proofoftheorem} \ref{thm_M_n_sing}. 

(4) When $R=2$ there is a unique configuraion of a spider where all the arms are stretched out and the body of the spider is located at the origin. 

\smallskip
{\bf (1) The {\boldmath $R=0$} case.}

When $R=0$ all the fixed endpoints $B_i$'s coincide with the origin. 
As was noticed in Remark \ref{symmetry}, $S^1$ acts on the configuration space ${\mathcal{M}}_n(0)$ as rotation. 
We can choose as ${\mathcal{M}}^{\p}_n(0)$ the configuration space of the spiders when $\vect b_1=\overrightarrow{B_1J_1}$ is fixed to be $\vect e_1$. 
When $n=2$, ${\mathcal{M}}^{\p}_2(0)$ is the configuration space of rhombics, which was proved to be homeomorphic to the union of three circles any two of which are tangent at a pair of distict points (\cite{To}). 

\smallskip
Suppose $\vect b_1=\vect e_1$. 
Then the body domain (i.e. the domain where the body can be located) is a circle 
$$\{C_{\theta}=(1+\cos\theta, \sin\theta):-\pi<\theta\le\pi\},$$
where $\theta$ is the angle of $\vect a_1=\overrightarrow{J_1C}$ form the $x$-axis. 

(i) When $\theta=0$ all the arms are stretched out. 
The configuration corresponds to a unique point $S$ in ${\mathcal{M}}^{\p}_n(0)$. 

(ii) When $\theta\ne0,\pi$ all the arms are bended. 
The space of the configurations can be given by 
$$\displaystyle %
A=\left\{\sigma_{(\theta;\e_2,\cdots,\e_n)}:\theta\in(-\pi, 0)\cup(0, \pi),\,\e_j\in\{+,-\}	\right\}\subset{\mathcal{M}}^{\p}_n(0),$$
where $\e_j$ denotes the index of the $j$-th arm. 
The space $A$ is homeomorphic to the disjoint union of $2^{n-1}$ copies of $(-\pi, 0)\cup(0, \pi)$. 
%
%

(iii) When $\theta=\pi$ all the arms are folded, which can rotate arond the origin except for the first arm. 
This configuration corresponds to a point in an $(n-1)$-torus 
$$
T^{n-1}=\{\tau_{(\theta_2, \cdots, \theta_n)}:0\le\theta_j<2\pi\}\subset{\mathcal{M}}^{\p}_n(0), 
$$
where $\theta_j$ denote the angle of $\vect b_j$ from the $x$-axis. 

\smallskip
Now let us how they are glued together. 

If the body approaches $(2,0)$ then all the arms tend to be stretched out. 
Therefore, 
$$
\lim_{\theta\to0}\sigma_{(\theta;\e_2,\cdots,\e_n)}=S,  
$$
which implies that $A\cup\{S\}$ is homeomorphic to a join of $2^{n-1}$ open intervals $\displaystyle \stackrel{2^{n-1}}{\vee} (-\pi, \pi)$. 

On the other hand, if the body approaches the origin then all the arms tend to be folded. 
The angle of $\vect b_j$ from the $x$-axis tends to be equal to either $0$ or $\pi$; $0$ if $\theta$ approaches $\pi$ from below and $\e_j=+$ or $\theta$ approaches $-\pi$ from above and $\e_j=-$, and $\pi$ otherwise. 
Therefore, 
$$
\lim_{\e\to +0}\sigma_{(\e_1(\pi-\e)\,;\,\e_2,\cdots,\e_n)}=\tau_{(\theta_2,\cdots,\theta_n)},
$$
where $\e_1\in\{+,-\}$ is the index of the first arm and $\theta_j$ is given by 
$$\theta_j=\left(1-(-1)^{\e_1\e_j}\right)\frac{\pi}2\in\{0,\pi\}. $$
It means that $2^{n-1}$ pairs of ``boundary points'' of $A\cup\{S\}\cong\displaystyle \stackrel{2^{n-1}}{\vee} (-\pi, \pi)$ are glued to mutually distinct $2^{n-1}$ points in $T^{n-1}$ respectively to produce ${\mathcal{M}}^{\p}_n(0)$. 

\smallskip
{\bf (2) The {\boldmath $n$} being even and {\boldmath $R=R_n=1$} case. }

Suppose $n$ is even $n=2m$ and $R=R_n=1$. 
We agree that the suffixes are considered modulo $n$ in what follows, i.e. $k+m$ means $k-m$ if $k+m>n$. 
The body domain is a curved $n$-gon $D$, where the $k$-th edge of $\partial D$ contains the endpoint $B_{k+m}$ in its interior (Figure \ref{6arms_R=1_D}). 
\begin{figure}[htbp]
\begin{center}
\includegraphics[width=.4\linewidth]{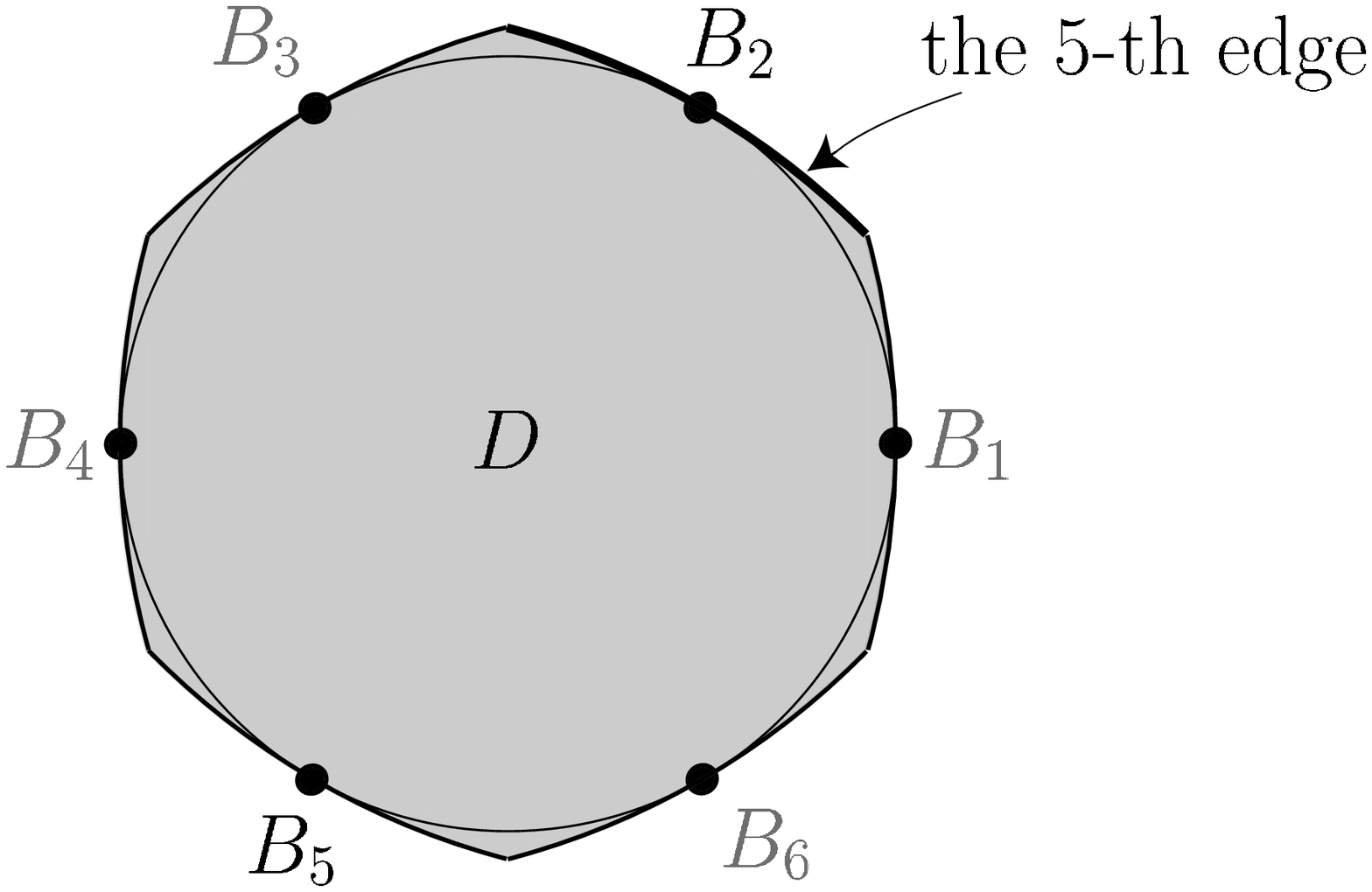}
\caption{The curved hexagon when $R=1$}
\label{6arms_R=1_D}
\end{center}
\end{figure}

Let $\overline{\mathcal{I}}$ denote the set of the multi-indices of the points in ${\mathcal{M}}_n(R)$: 
$$\overline{\mathcal{I}}={\mathcal{I}}\cup{\mathcal{I}}_{SE},$$
where ${\mathcal{I}}$ is same as (\ref{f_multi-indices}) and ${\mathcal{I}}_{SE}$ is given by 
\begin{equation}\label{f_I_SE}
{\mathcal{I}}_{SE}=\left\{\bar{\vect{\varepsilon}}^{\circ}=(\bar{\e}^{\circ}_1, \cdots, \bar{\e}^{\circ}_n)\left|
\begin{array}{l}
\bar{\e}^{\circ}_k\in\{+,-,0,\infty\}, \\[1mm]
\sharp\{i:\bar{\e}^{\circ}_i=0\}=1, \>\sharp\{j:\bar{\e}^{\circ}_j=\infty\}=1,\\[1mm]
\mbox{if}\>\>\bar{\e}^{\circ}_k=\infty\>\>\mbox{then}\>\>\bar{\e}^{\circ}_{k+m}=0
\end{array}\right.
\right\}.
\end{equation}
Let ${\mathcal{I}}_m$ $(0\le m\le 2)$, $D_{\vect{\varepsilon}}$ $(\vect{\varepsilon}\in{\mathcal{I}}_2)$, $E_{\vect{\varepsilon}^{\p}}$ $(\vect{\varepsilon}^{\p}\in{\mathcal{I}}_1)$, and $V_{\vect{\varepsilon}^{\p\p}}$ $(\vect{\varepsilon}^{\p\p}\in{\mathcal{I}}_0)$ are given by (\ref{f_I_m}),  (\ref{f_D_e}), (\ref{f_E_e}), and (\ref{f_V_e}) as in the non-singular case. 
Each $D_{\vect{\varepsilon}}$ is homeomorphic to $\textrm{Int} D$, where $D$ is the curved $n$-gon given by (\ref{curved_n-gon}). 
We remark that, unlike in the non-singular case, $E_{\vect{\varepsilon}^{\p}}$ is homeomorphic to an open interval minus one point; 
if $\e_k^{\p}=0$ then $E_{\vect{\varepsilon}^{\p}}$ is homeomorphic to the interior of the $k$-th edge of $\partial D$ minus $B_{k+m}$. 

Put, for $\bar{\vect{\varepsilon}}^{\circ}\in{\mathcal{I}}_{SE}$, 
$$
S^1_{\bar{\vect{\varepsilon}}^{\circ}}=\left\{\vect x\in {\mathcal{M}}_n(R):\,\vect{\varepsilon}(\vect x)=\bar{\vect{\varepsilon}}^{\circ}\right\}. 
$$
The configuration space ${\mathcal{M}}_n(R)$ can be decomposed as the disjoint union: 
\begin{equation}\label{f_decomp_M_n_even_R=1}
{\mathcal{M}}_n(R)=\left(\bigcup_{\vect{\varepsilon}\in{\mathcal{I}}_2}D_{\vect{\varepsilon}}
\cup\bigcup_{\vect{\varepsilon}^{\p}\in{\mathcal{I}}_1}E_{\vect{\varepsilon}^{\p}}
\cup\bigcup_{\vect{\varepsilon}^{\p\p}\in{\mathcal{I}}_0}V_{\vect{\varepsilon}^{\p\p}}\right)
\cup\bigcup_{\bar{\vect{\varepsilon}}^{\circ}\in{\mathcal{I}}_{SE}}S^1_{\bar{\vect{\varepsilon}}^{\circ}}. 
\end{equation}
The first term of the right hand side is homeomorphic to $n2^{n-1}$-times punctured orietable surface of genus $1-2^{n-1}+n2^{n-3}$, and $\bigcup_{\bar{\vect{\varepsilon}}^{\circ}\in{\mathcal{I}}_{SE}}S^1_{\bar{\vect{\varepsilon}}^{\circ}}$ is the disjoint union of $n2^{n-2}$ circles. 
We see how $\cup_{\bar{\vect{\varepsilon}}^{\circ}\in{\mathcal{I}}_{SE}}S^1_{\bar{\vect{\varepsilon}}^{\circ}}$ is glued to it in what follows. 
The argument in the non-singular and $0<R<R_n$ case runs parallel after modification according to the following differences: 
\begin{enumerate}
\item[(i)] Since $B_k$ is not located in $\mbox{\rm Int}D$ but in the interior of an edge of $\partial D$, the body cannot approach $B_k$ from all the directions, but from the ``half'' of them. 
\item[(ii)] Since not only the folded $k$-th arm but also the stretched-out $(k+m)$-th arm can be relaxed to bended arms, $S^1_{\bar{\vect{\varepsilon}}^{\circ}}$ intersects the closure of four $D_{\vect{\varepsilon}}$'s. 
\end{enumerate}

Let $e_{\bar{\vect{\varepsilon}}^{\circ}}^{i\theta}$ $(\bar{\vect{\varepsilon}}^{\circ}\in{\mathcal{I}}_{SE}, \theta\in\mathbb{R}/(2\pi\mathbb{Z}))$ denote a point in $S^1_{\bar{\vect{\varepsilon}}^{\circ}}\subset{\mathcal{M}}_n(R)$ where the folded arm has angle $\theta$ from the positive direction of the $x$-axis. 

Suppose $\bar{\vect{\varepsilon}}^{\circ}=(\bar{\e}_1^{\circ}, \cdots, \bar{\e}_n^{\circ})\in{\mathcal{I}}_{SE}$ satisfies $\bar{\e}_k^{\circ}=\infty$ and $\bar{\e}_{k+m}^{\circ}=0$, i.e. the $k$-th arm is folded. 

Define $\vect{\varepsilon}_{\sigma \tau}\in{\mathcal{I}}_{2}\cup{\mathcal{I}}_1$ $(\sigma\in\{+,-\}, \tau\in\{+,-,0\})$ by 
$$
\vect{\varepsilon}_{\sigma \tau}=(\e_1, \cdots, \e_n) \>\>\> \mbox{with}\>\>
\left\{\begin{array}{l}
\e_j=\bar{\e}_j^{\circ} \>\>\>\mbox{if}\>\>\> j\ne k,k+m\\[1mm]
\e_k=\sigma, \,{\e}_{k+m}=\tau.
\end{array}
\right.
$$
Let $S$ (or $F$) be the configuration in $S^1_{\bar{\vect{\varepsilon}}^{\circ}}$ where the folded $k$-th arm and the stretched-out $(k+m)$-th arm are collinear and the folded arm is outside (or respectively, inside) the curved $n$-gon $D$ (Figures \ref{6arms_R=1_S} and \ref{6arms_R=1_F}):
$$
S=e_{\bar{\vect{\varepsilon}}^{\circ}}^{i\,\frac{2(k-1)}n\pi},\>
F=e_{\bar{\vect{\varepsilon}}^{\circ}}^{i\left(\frac{2(k-1)}n\pi+\pi\right)}.
$$
\begin{figure}[htbp]
\begin{center}
\begin{minipage}{.45\linewidth}
\begin{center}
\includegraphics[width=0.6\linewidth]{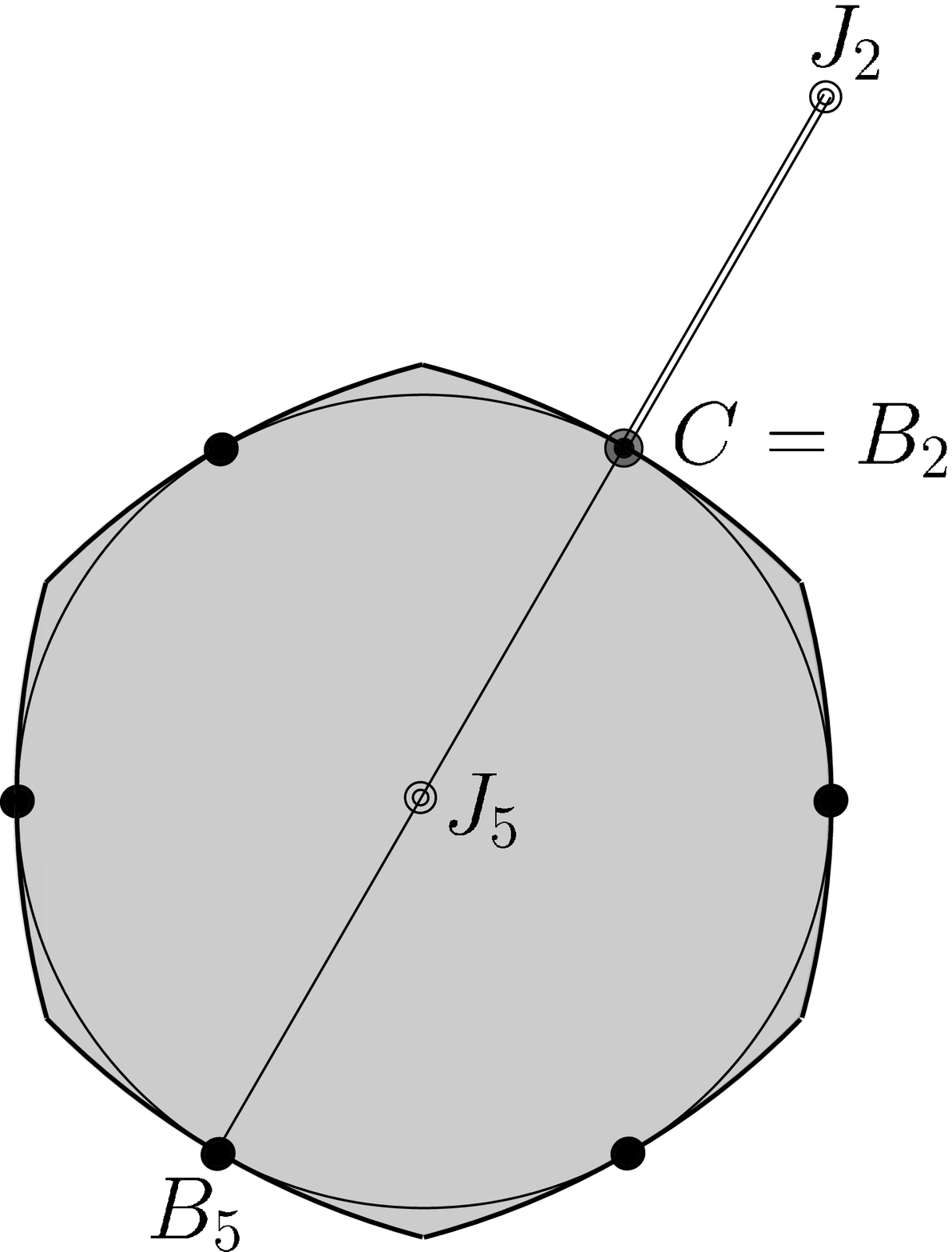}
\caption{$S=e_{\bar{\vect{\varepsilon}}^{\circ}}^{i\,\frac{2(k-1)}n\pi}$. Only the second and the $5$-th arms are drawn. }
\label{6arms_R=1_S}
\end{center}
\end{minipage}
\hskip 0.4cm
\begin{minipage}{.45\linewidth}
\begin{center}
\includegraphics[width=0.6\linewidth]{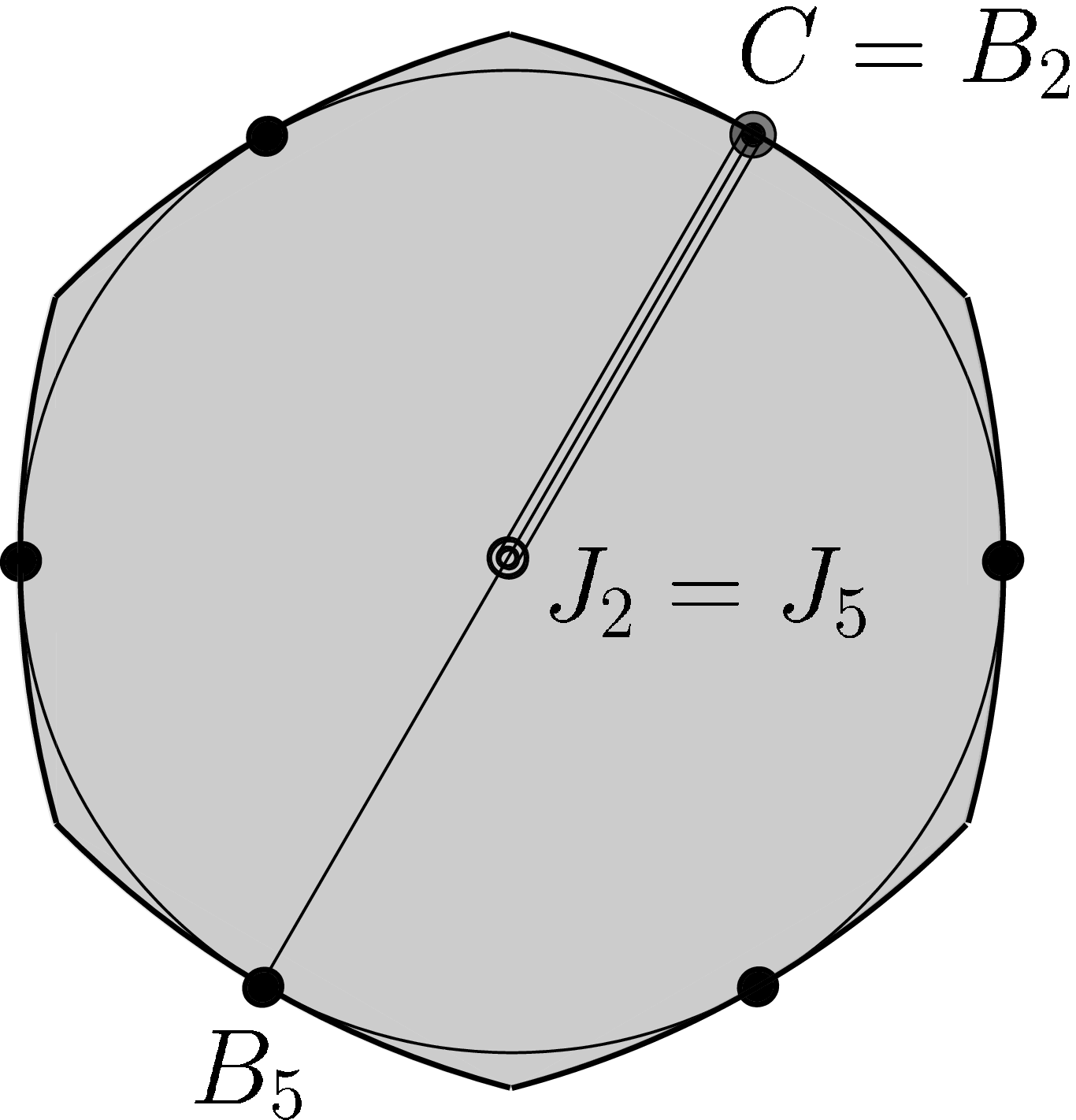}
\caption{$F=e_{\bar{\vect{\varepsilon}}^{\circ}}^{i\left(\frac{2(k-1)}n\pi+\pi\right)}$. Only the second and the $5$-th arms are drawn. }
\label{6arms_R=1_F}
\end{center}
\end{minipage}
\end{center}
\end{figure}
Let $\Gamma_+$ (or $\Gamma_-$) be an open subarc of $S^1_{\bar{\vect{\varepsilon}}^{\circ}}$ from $S$ to $F$ (or respectively, from $F$ to $S$):
$$
\begin{array}{l}
\displaystyle \Gamma_+=\left\{e_{\bar{\vect{\varepsilon}}^{\circ}}^{i\theta}:\frac{2(k-1)}n\pi<\theta<\frac{2(k-1)}n\pi+\pi\right\},\\[4mm]
\displaystyle \Gamma_-=\left\{e_{\bar{\vect{\varepsilon}}^{\circ}}^{i\theta}:\frac{2(k-1)}n\pi-\pi<\theta<\frac{2(k-1)}n\pi\right\}.
\end{array}
$$

Let us show that a point in $\Gamma_+$ (or $\Gamma_-$) is the limit of a sequence of points in $D_{\vect{\varepsilon}_{++}}$ or $D_{\vect{\varepsilon}_{+-}}$ (or respectively, $D_{\vect{\varepsilon}_{-+}}$ or $D_{\vect{\varepsilon}_{--}}$ ). 
Remark first that if $C\in \mbox{\rm Int}D$ then the angle $\theta^{\p}$ of $\overrightarrow{B_kC}$ from the $x$-axis satisfies 
$$\frac{2(k-1)}n\pi+\frac{1}2\pi<\theta^{\p}<\frac{2(k-1)}n\pi+\frac{3}2\pi.$$
Consider a sequence of points in $D_{\vect{\varepsilon}_{++}}$ or $D_{\vect{\varepsilon}_{+-}}$ whose bodies are located at 
$$B_k+\delta\left(\cos\theta^{\p}, \sin\theta^{\p}\right) \>\>\>\left(\delta>0, \>\frac{2(k-1)}n\pi+\frac{1}2\pi<\theta^{\p}<\frac{2(k-1)}n\pi+\frac{3}2\pi\right).$$
Then Figure \ref{6arms_R=1_theta}, which can be obtained by a slight modification from Figure \ref{S1_dD}, implies that the formula (\ref{f_S1_as_limit_from_D}) also holds in this case. 
\begin{figure}[htbp]
\begin{center}
\begin{minipage}{.45\linewidth}
\begin{center}
\includegraphics[width=.6\linewidth]{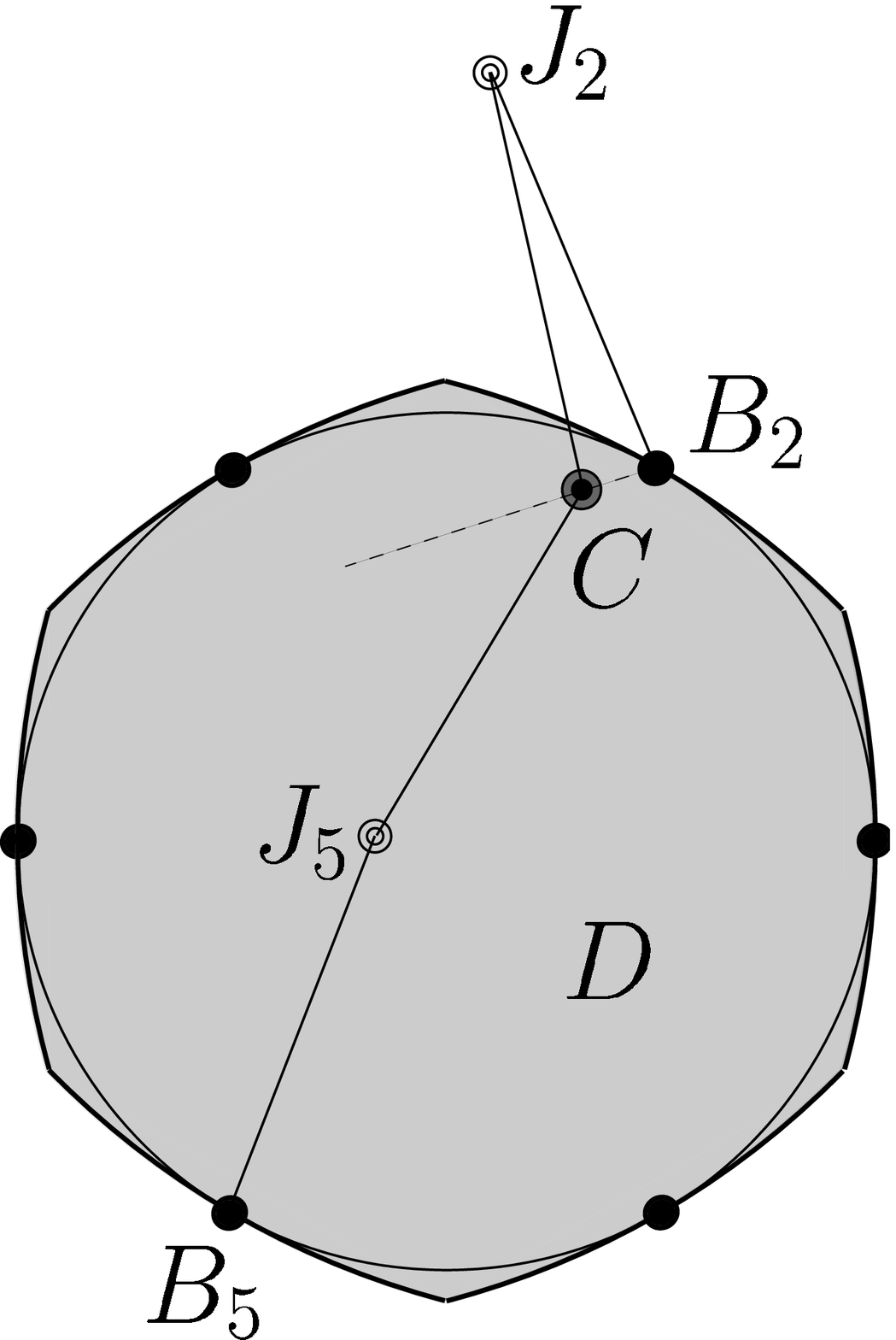}
\caption{The body $C$ approaches $B_2$ from inside the curved hexagon $D$}
\label{6arms_R=1_theta}
\end{center}
\end{minipage}
\hskip 0.4cm
\begin{minipage}{.45\linewidth}
\begin{center}
\includegraphics[width=.8\linewidth]{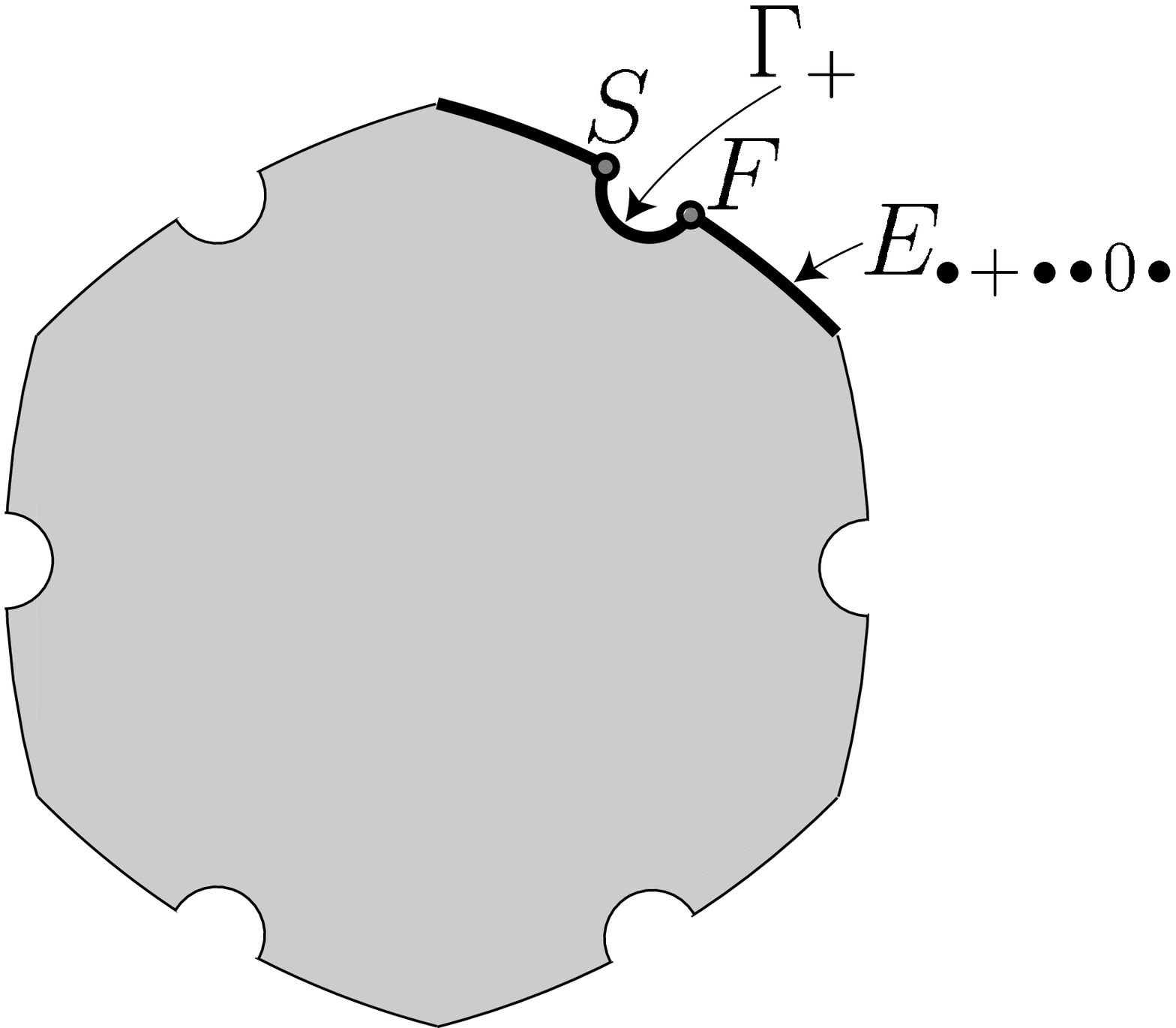}
\caption{The $5$-th edge of $D_{\bullet + \bullet \bullet + \bullet}$ or $D_{\bullet + \bullet \bullet - \bullet}$}
\label{6arms_R=1_dD}
\end{center}
\end{minipage}
\end{center}
\end{figure}
It follows that the limet of this sequence as $\delta$ goes down to $+0$ is the point $e_{\bar{\vect{\varepsilon}}^{\circ}}^{i\left(\theta^{\p}-\frac{\pi}2\right)}$, where $\theta^{\p}-\frac{\pi}2$ satisfies 
$$\frac{2(k-1)}n\pi<\theta^{\p}-\frac{\pi}2<\frac{2(k-1)}n\pi+\pi.$$
It implies that the limit $e_{\bar{\vect{\varepsilon}}^{\circ}}^{i\left(\theta^{\p}-\frac{\pi}2\right)}$ belongs to $\Gamma_+$, and conversely that any point in $\Gamma_+$ can be expressed as a limit of this kind. 

On the other hand, the point $S$ can be expressed as the limit in two ways; as the limit of a sequence of points in $E_{\vect{\varepsilon}_{+0}}\subset \overline{D_{\vect{\varepsilon}_{++}}}\cap\overline{D_{\vect{\varepsilon}_{+-}}}$ whose bodies are located at 
$$
B_{k+m}+2\left(\cos\left(\frac{2(k-1)}n\pi+\delta\right),\sin\left(\frac{2(k-1)}n\pi+\delta\right)\right) \>\>\>(\delta>0)
$$
as $\delta$ goes down to $+0$, i.e. the body approaches $B_k$ from the ``front'' side, and as the limit of a sequence of points the points in $E_{\vect{\varepsilon}_{-0}}\subset \overline{D_{\vect{\varepsilon}_{-+}}}\cap\overline{D_{\vect{\varepsilon}_{--}}}$ whose bodies approach $B_k$ from the ``back'' side. 

It follows that the $(k+m)$-th edges of $\overline{D_{\vect{\varepsilon}_{++}}}$ and of $\overline{D_{\vect{\varepsilon}_{+-}}}$ are both given by 
$$\overline{D_{\vect{\varepsilon}_{++}}}\cap\overline{D_{\vect{\varepsilon}_{+-}}}=E_{\vect{\varepsilon}_{+0}}\cup\Gamma_+\cup\{S,F\}$$ 
(Figure \ref{6arms_R=1_dD}). 
Similarly 
$$\overline{D_{\vect{\varepsilon}_{-+}}}\cap\overline{D_{\vect{\varepsilon}_{--}}}=E_{\vect{\varepsilon}_{-0}}\cup\Gamma_-\cup\{S,F\}.$$ 
%
%
Therefore, 
$$
\begin{array}{rcl}
S^1_{\bar{\vect{\varepsilon}}^{\circ}}=\Gamma_+\cup\Gamma_-\cup\{S,F\}\!\!&\!\!\subset\!\!&\!\!\left(\overline{D_{\vect{\varepsilon}_{++}}}\cap\overline{D_{\vect{\varepsilon}_{+-}}}\right)\cup\left(\overline{D_{\vect{\varepsilon}_{-+}}}\cap\overline{D_{\vect{\varepsilon}_{--}}}\right),\\[2mm]
\{S,F\}\!\!&\!\!=\!\!&\!\!\left(\overline{D_{\vect{\varepsilon}_{++}}}\cap\overline{D_{\vect{\varepsilon}_{+-}}}\right)\cap\left(\overline{D_{\vect{\varepsilon}_{-+}}}\cap\overline{D_{\vect{\varepsilon}_{--}}}\right),
\end{array}
$$
which implies that $S^1_{\bar{\vect{\varepsilon}}^{\circ}}$ passes through two $1$-handles which are pinched at the middle, $S$ and $F$ (Figure \ref{pinched_1-handle}). 
\begin{figure}[htbp]
\begin{center}
\includegraphics[width=.5\linewidth]{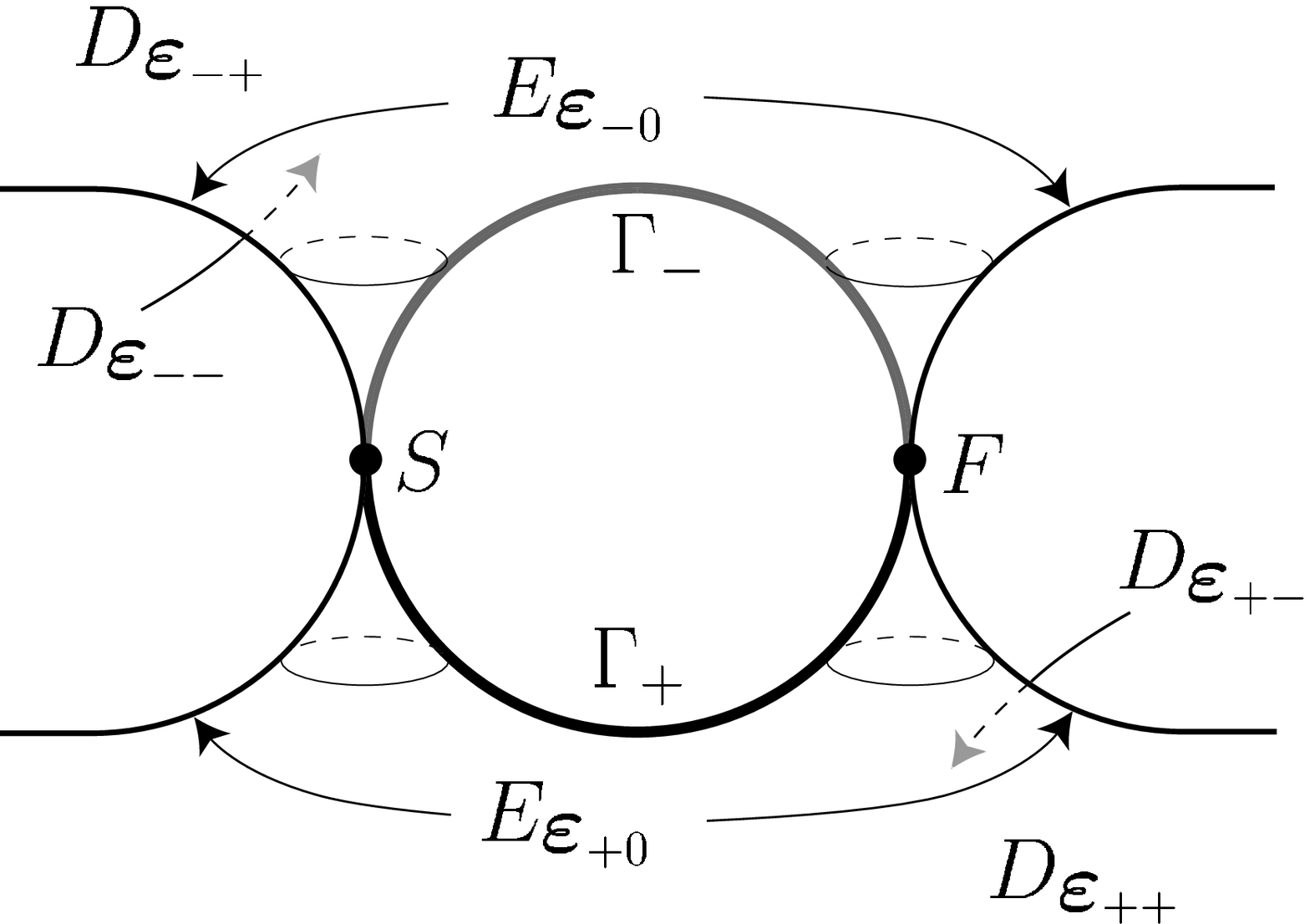}
\caption{$S^1_{\bar{\vect{\varepsilon}}^{\circ}}$ and $D_{\vect{\varepsilon}_{++}}, D_{\vect{\varepsilon}_{+-}}, D_{\vect{\varepsilon}_{-+}}$, and $D_{\vect{\varepsilon}_{--}}$. }
\label{pinched_1-handle}
\end{center}
\end{figure}

\smallskip
{\bf (3) The {\boldmath $n$} being odd and {\boldmath $R=R_n$} case.} 

Suppose $n$ is odd $n=2m+1$ and $R=R_n$. 
We agree that the suffixes are considered modulo $n$ in what follows. 
The body domain is a curved $n$-gon $D$ whose vertices are $B_1, \cdots, B_n$. 
Let $\widehat{\mathcal{I}}$ denote the set of the multi-indices of the points in ${\mathcal{M}}_n(R)$: 
$$\widehat{\mathcal{I}}={\mathcal{I}}_2\cup{\mathcal{I}}_1\cup{\mathcal{I}}_{SO},$$
where ${\mathcal{I}}_2$ and ${\mathcal{I}}_1$ are given by (\ref{f_I_m}) and ${\mathcal{I}}_{SO}$ is given by 
\begin{equation}\label{f_I_SO}
{\mathcal{I}}_{SO}=\left\{\widehat{\vect{\varepsilon}}^{\circ}=(\hat{\e}^{\circ}_1, \cdots, \hat{\e}^{\circ}_n)\left|
\begin{array}{l}
\hat{\e}^{\circ}_k\in\{+,-,0,\infty\}, \\[1mm]
\sharp\{i:\hat{\e}^{\circ}_i=0\}=2, \>\sharp\{j:\hat{\e}^{\circ}_j=\infty\}=1,\\[1mm]
\mbox{if}\>\>\hat{\e}^{\circ}_k=\infty\>\>\mbox{then}\>\>\hat{\e}^{\circ}_{k+m}=\hat{\e}^{\circ}_{k+m+1}=0
\end{array}\right.
\right\}.
\end{equation}
Let $D_{\vect{\varepsilon}}$ $(\vect{\varepsilon}\in{\mathcal{I}}_2)$, $E_{\vect{\varepsilon}^{\p}}$ $(\vect{\varepsilon}^{\p}\in{\mathcal{I}}_1)$ are given by (\ref{f_D_e}) and  (\ref{f_E_e}) as in the non-singular case. 
Each $D_{\vect{\varepsilon}}$ is homeomorphic to $\textrm{Int} D$, where $D$ is the curved $n$-gon given by (\ref{curved_n-gon}). 
We remark that, unlike in the non-singular case, if the body is located at a vertex of $D$ then there is a folded arm which can rotate. 
Therefore, the set of $0$-cells $V_{\vect{\varepsilon}^{\p\p}}$ $(\vect{\varepsilon}^{\p\p}\in{\mathcal{I}}_0)$ in the non-singular case (\ref{f_V_e}) should be replaced by circles 
$$
S^1_{\widehat{\vect{\varepsilon}}^{\circ}}=\left\{\vect x\in {\mathcal{M}}_n(R):\,\vect{\varepsilon}(\vect x)=\widehat{\vect{\varepsilon}}^{\circ}\right\} \>\>\>(\widehat{\vect{\varepsilon}}^{\circ}\in{\mathcal{I}}_{SO}). 
$$
The configuration space ${\mathcal{M}}_n(R)$ can be decomposed as the disjoint union: 
\begin{equation}\label{f_decomp_M_n_odd_R=1}
{\mathcal{M}}_n(R)=\left(\bigcup_{\vect{\varepsilon}\in{\mathcal{I}}_2}D_{\vect{\varepsilon}}
\cup\bigcup_{\vect{\varepsilon}^{\p}\in{\mathcal{I}}_1}E_{\vect{\varepsilon}^{\p}}\right)
\cup\bigcup_{\widehat{\vect{\varepsilon}}^{\circ}\in{\mathcal{I}}_{SO}}S^1_{\widehat{\vect{\varepsilon}}^{\circ}}. 
\end{equation}
The first term of the right hand side is homeomorphic to $n2^{n-2}$-times punctured orietable surface of genus $1-2^{n-1}+n2^{n-3}$, and $\bigcup_{\widehat{\vect{\varepsilon}}^{\circ}\in{\mathcal{I}}_{SO}}S^1_{\widehat{\vect{\varepsilon}}^{\circ}}$ is the disjoint union of $n2^{n-3}$ circles. 
We see how $\cup_{\widehat{\vect{\varepsilon}}^{\circ}\in{\mathcal{I}}_{SO}}S^1_{\widehat{\vect{\varepsilon}}^{\circ}}$ is glued to it in what follows. 
The argument in the previous case runs parallel. 
What is different is that since $B_k$ is located at a vertex of $D$ the range of the possible directions of approaches of the body to $B_k$ is restricted to $(1-\frac1n)\frac{\pi}2$. 

Let $e_{\widehat{\vect{\varepsilon}}^{\circ}}^{i\theta}$ $(\widehat{\vect{\varepsilon}}^{\circ}\in{\mathcal{I}}_{SO}, \theta\in\mathbb{R}/(2\pi\mathbb{Z}))$ denote a point in $S^1_{\widehat{\vect{\varepsilon}}^{\circ}}\subset{\mathcal{M}}_n(R)$ where the folded arm has angle $\theta$ from the positive direction of the $x$-axis. 

Suppose $\widehat{\vect{\varepsilon}}^{\circ}=(\hat{\e}_1^{\circ}, \cdots, \hat{\e}_n^{\circ})\in{\mathcal{I}}_{SO}$ satisfies $\bar{\e}_k^{\circ}=\infty$ and $\hat{\e}_{k+m}^{\circ}=\hat{\e}_{k+m+1}^{\circ}=0$, i.e. the $k$-th arm is folded. 

Define $\vect{\varepsilon}_{\sigma \tau \tau^{\p}}\in{\mathcal{I}}_{2}\cup{\mathcal{I}}_1$ $(\sigma\in\{+,-\}, \tau,\tau^{\p}\in\{+,-,0\}, (\tau,\tau^{\p})\ne(0,0))$ by 
$$
\vect{\varepsilon}_{\sigma \tau \tau^{\p}}=(\e_1, \cdots, \e_n) \>\>\> \mbox{with}\>\>
\left\{\begin{array}{l}
\e_j=\bar{\e}_j^{\circ} \>\>\>\mbox{if}\>\>\> j\ne k,k+m,k+m+1\\[1mm]
\e_k=\sigma, \,{\e}_{k+m}=\tau, \,{\e}_{k+m+1}=\tau^{\p}.
\end{array}
\right.
$$
Put 
$$
\begin{array}{ll}
S_+=e_{\widehat{\vect{\varepsilon}}^{\circ}}^{i\left(\frac{2(k-1)}n\pi+\frac1{2n}\pi\right)}, 
&T_+=e_{\widehat{\vect{\varepsilon}}^{\circ}}^{i\left(\frac{2(k-1)}n\pi+\pi-\frac1{2n}\pi\right)}, \\[2mm]
S_-=e_{\widehat{\vect{\varepsilon}}^{\circ}}^{i\left(\frac{2(k-1)}n\pi-\pi+\frac1{2n}\pi\right)}, 
&T_-=e_{\widehat{\vect{\varepsilon}}^{\circ}}^{i\left(\frac{2(k-1)}n\pi-\frac1{2n}\pi\right)}. 
\end{array}
$$
Let $\Gamma_+$, $\Gamma_-$, $\Gamma_S$, and $\Gamma_F$ be open subarcs of $S^1_{\widehat{\vect{\varepsilon}}^{\circ}}$ from $S_+$ to $T_+$, from $S_-$ to $T_-$, from $T_-$ to $S_+$, and from $T_+$ to $S_-$ respectively (Figure \ref{3arms_R=R3}):
$$
\begin{array}{l}
\displaystyle \Gamma_+=\left\{e_{\widehat{\vect{\varepsilon}}^{\circ}}^{i\theta}:\frac{2(k-1)}n\pi+\frac1{2n}\pi<\theta<\frac{2(k-1)}n\pi+\pi-\frac1{2n}\pi\right\},\\[4mm]
\displaystyle \Gamma_-=\left\{e_{\widehat{\vect{\varepsilon}}^{\circ}}^{i\theta}:\frac{2(k-1)}n\pi-\pi+\frac1{2n}\pi<\theta<\frac{2(k-1)}n\pi-\frac1{2n}\pi\right\},\\[4mm]
\displaystyle \Gamma_S=\left\{e_{\widehat{\vect{\varepsilon}}^{\circ}}^{i\theta}:\frac{2(k-1)}n\pi-\frac1{2n}\pi<\theta<\frac{2(k-1)}n\pi+\frac1{2n}\pi\right\},\\[4mm]
\displaystyle \Gamma_F=\left\{e_{\widehat{\vect{\varepsilon}}^{\circ}}^{i\theta}:\frac{2(k-1)}n\pi+\pi-\frac1{2n}\pi<\theta<\frac{2(k-1)}n\pi+\pi+\frac1{2n}\pi\right\}.
\end{array}
$$

\begin{figure}[htbp]
\begin{center}
\includegraphics[width=.4\linewidth]{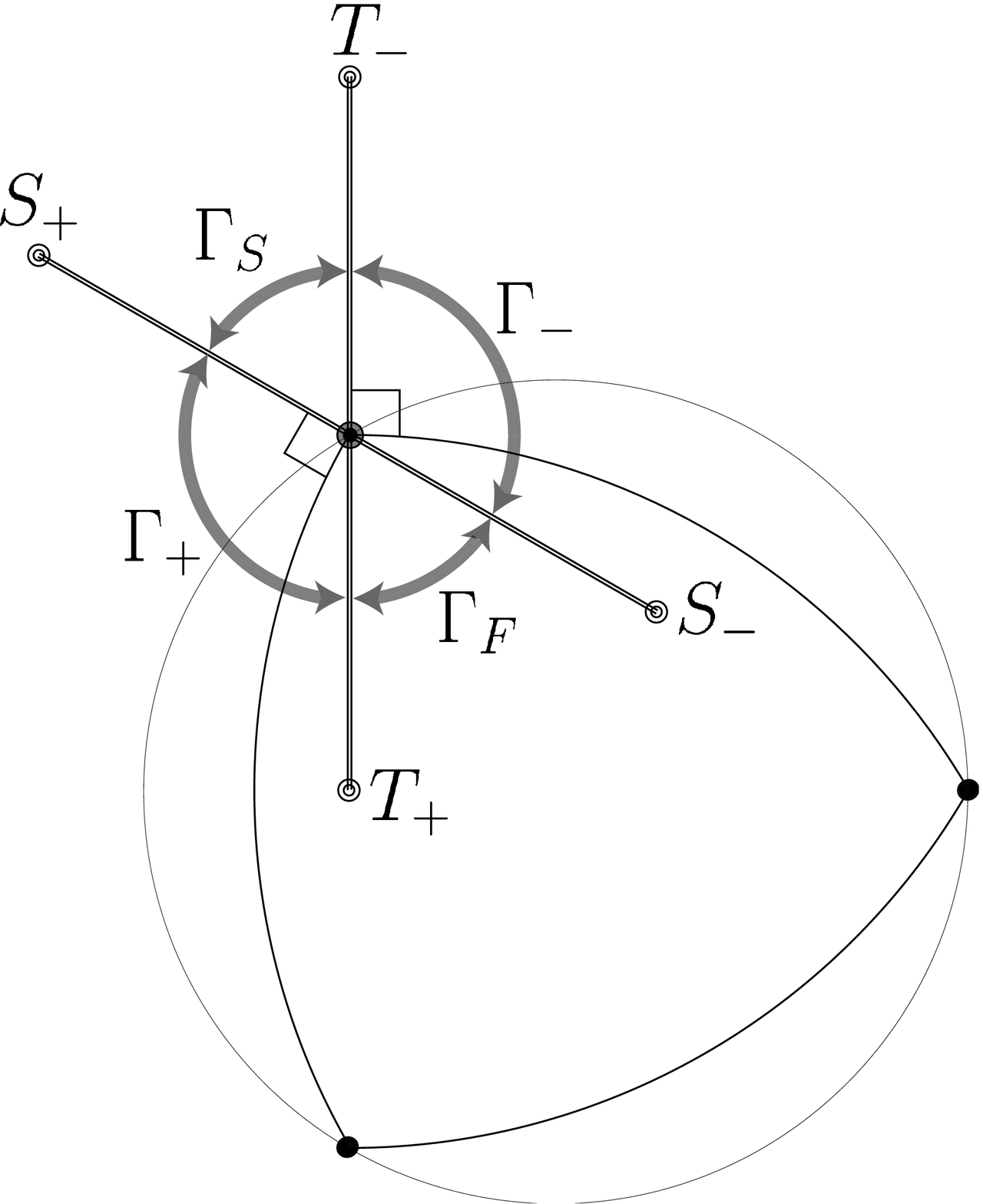}
\caption{$S_{\pm}, T_{\pm}$, $\Gamma_{\pm}$, and $\Gamma_S$, $\Gamma_F$ when $n=3$}
\label{3arms_R=R3}
\end{center}
\end{figure}

Just like in the previous case, a point in $\Gamma_+$ (or $\Gamma_-$) is the limit of a sequence of points in $D_{\vect{\varepsilon}_{+ \tau \tau^{\p}}}$ (or respectively, $D_{\vect{\varepsilon}_{- \tau \tau^{\p}}}$) $(\tau,\tau^{\p}\in\{+,-\})$. 
The point $S_+$ (or $T_+$) can be expressed as the limit of a sequence of points in $E_{\vect{\varepsilon}_{+ \tau 0}}\subset \overline{D_{\vect{\varepsilon}_{+ \tau +}}}\cap\overline{D_{\vect{\varepsilon}_{+ \tau -}}}$ (or respectively, $E_{\vect{\varepsilon}_{+ 0 \tau^{\p}}}\subset \overline{D_{\vect{\varepsilon}_{+ + \tau^{\p}}}}\cap\overline{D_{\vect{\varepsilon}_{+ - \tau^{\p}}}}$) $(\tau,\tau^{\p}\in\{+,-\})$. 
\begin{figure}[htbp]
\begin{center}
\begin{minipage}{.44\linewidth}
\begin{center}
\includegraphics[width=\linewidth]{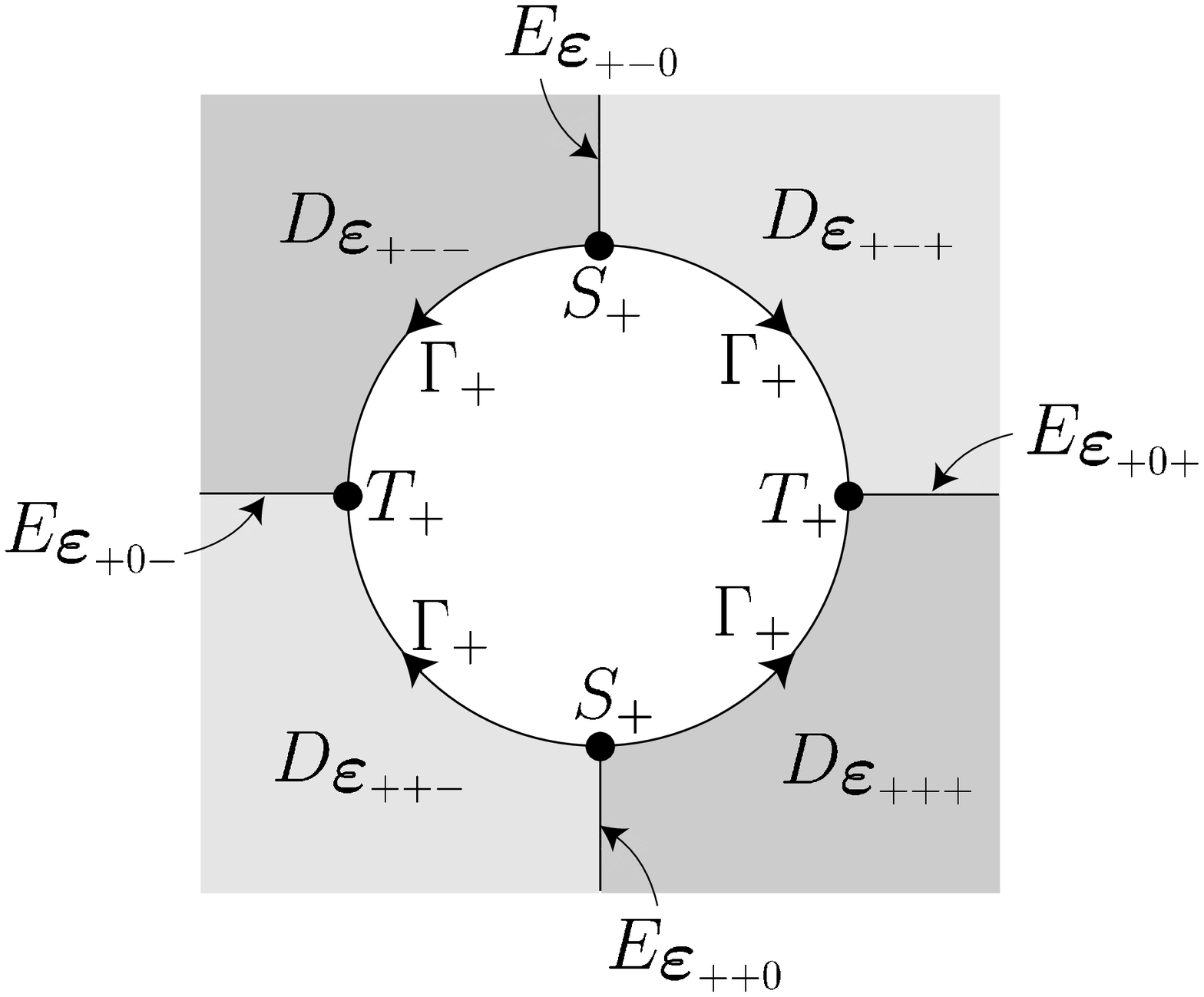}
\caption{A plane minus an open disc is being stitched up along $\Gamma_+$ to produce the space in the next Figure}
\label{idf_circle_2}
\end{center}
\end{minipage}
\hskip 0.4cm
\begin{minipage}{.48\linewidth}
\begin{center}
\includegraphics[width=\linewidth]{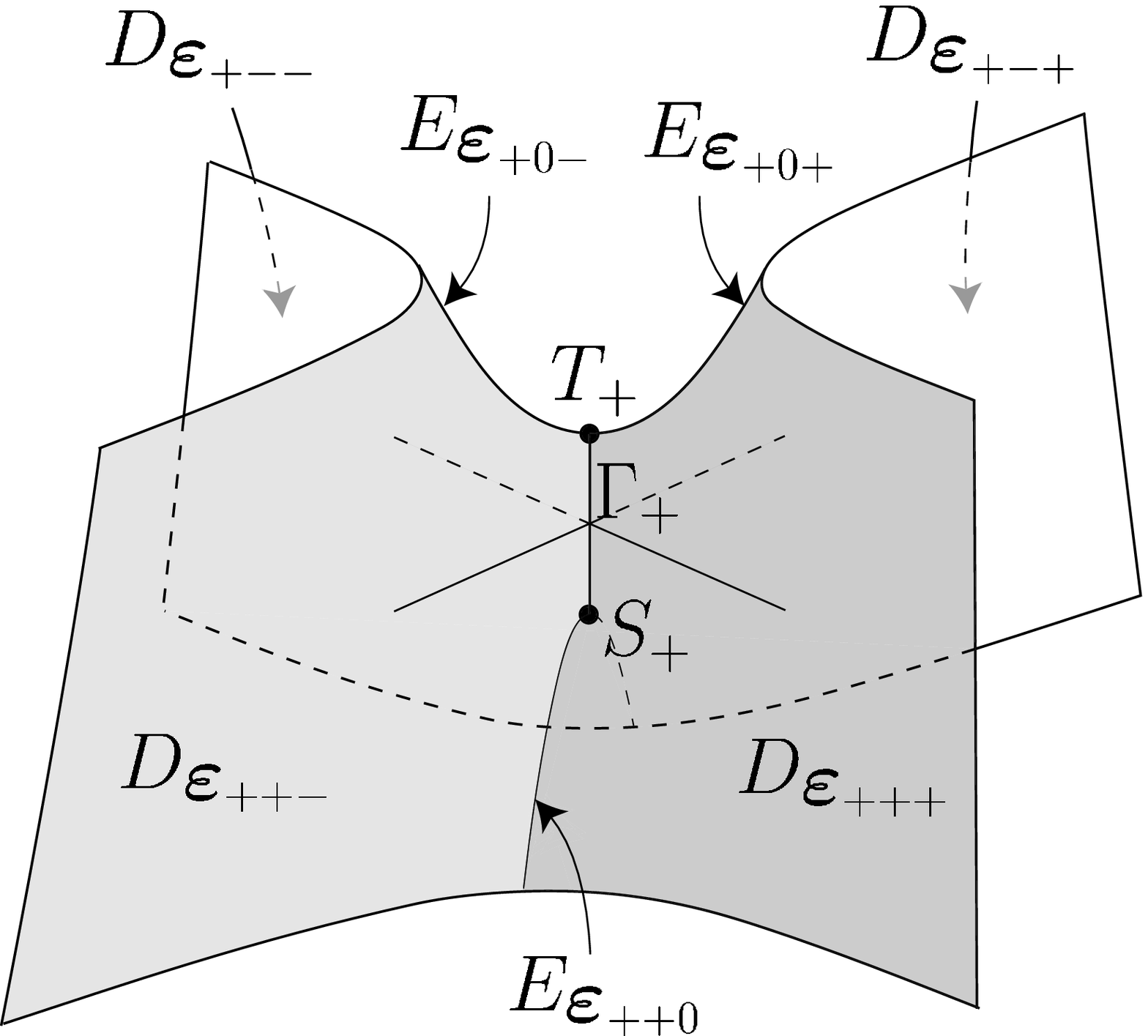}
\caption{A space homeomorphic to the stitched-up disc in Figure \ref{stitched_cap}}
\label{M_odd_sing_2}
\end{center}
\end{minipage}
\end{center}
\end{figure}

Figures \ref{idf_circle_2} and \ref{M_odd_sing_2} illustrate how the subarc $\overline{\Gamma_+}=\Gamma_+\cup\{S_+, T_+\}$ of $S^1_{\widehat{\vect{\varepsilon}}^{\circ}}$ is glued to $\displaystyle \cup_{\tau, \tau^{\p}\in\{+,-\}}\overline{D_{\vect{\varepsilon}_{+\tau\tau^{\p}}}}$. 
The subarc $\overline{\Gamma_-}=\Gamma_-\cup\{S_-, T_-\}$ of $S^1_{\widehat{\vect{\varepsilon}}^{\circ}}$ is glued to $\displaystyle \cup_{\tau, \tau^{\p}\in\{+,-\}}\overline{D_{\vect{\varepsilon}_{-\tau\tau^{\p}}}}$ similarly. 
Their endpoints $S_+, T_-$ and $S_-, T_+$ are joined by $\overline{\Gamma_S}$ and $\overline{\Gamma_F}$. 
This completes the proof. 
\end{proofoftheorem}

\begin{figure}[htbp]
\begin{center}
\includegraphics[width=.3\linewidth]{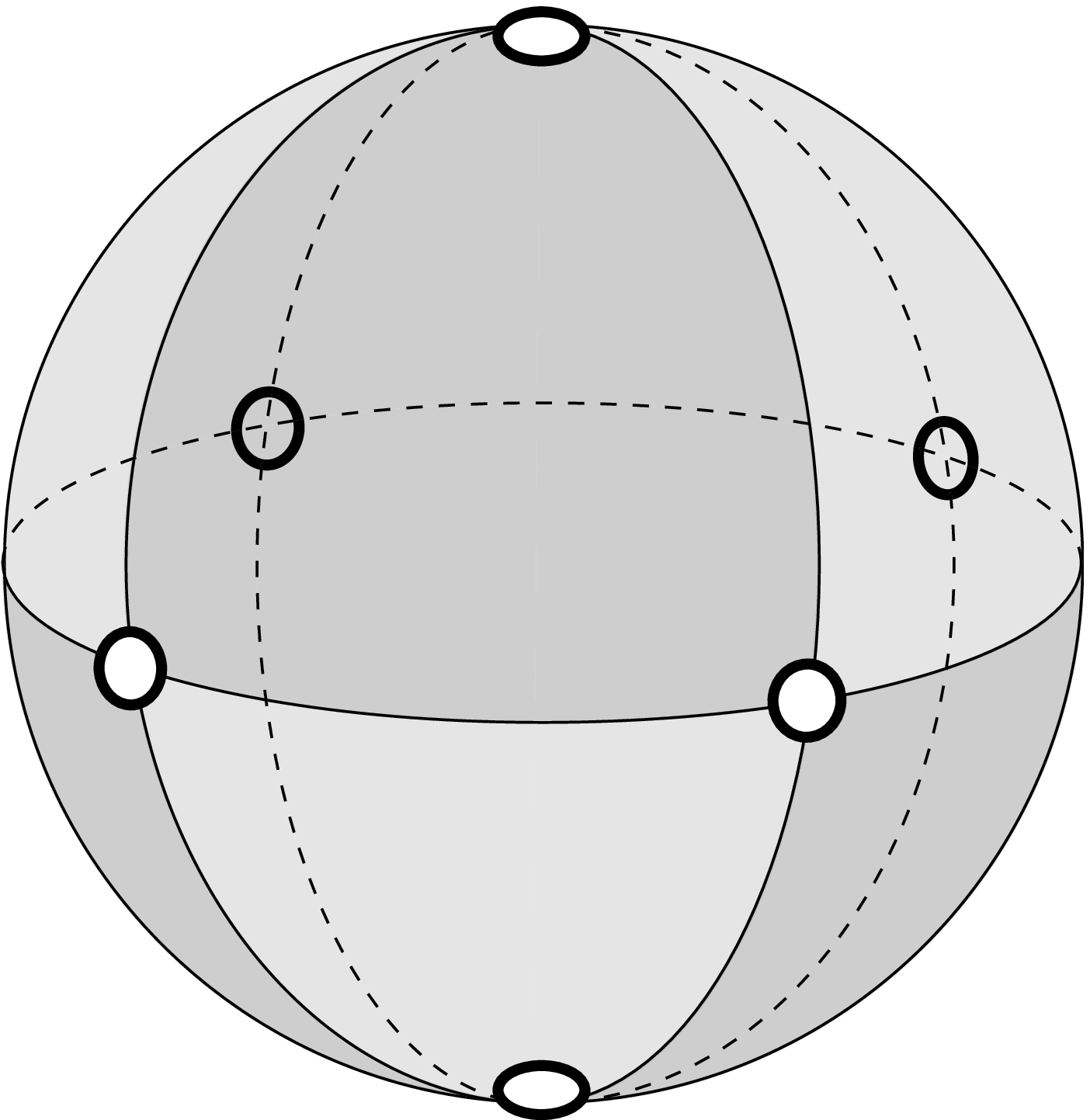}
\caption{${\mathcal{M}}_3(R)$ $(\frac{2}{\sqrt3}<R<2)$ minus $6$ open discs}
\label{octahedron}
\end{center}
\end{figure}

\medskip
For example, when $n=3$, the configuration space ${\mathcal{M}}_3(\frac{2}{\sqrt3})$ can be obtained by first replacing $6$ discs of an $S^2$ (Figure \ref{octahedron}) by $6$ copies of the space illustrated in Figure \ref{M_odd_sing_2}, and then joining $3$ pairs of pair of points (copies of $S_{\pm}$ and $T_{\pm}$) by $6$ arcs. 

\section{Appendix}\label{sec_app}
We like to end this article by proposing a problem. 
The linkages that we have studied in this paper have maximum symmetry. 
The configuration spaces of the spiders without the symmetry can produce other types of spaces. 
For example, when $n=2$, the configuration space is nothing but the moduli space of pentagons, which can produce connected orientable closed surfaces of genera from $0$ up to $4$ (reported in \cite{Ka-Mi1}), whereas only $S^2$ and $\Sigma_4$ can occur in our most symmetric cases. 

It seems to the author that the configuration spaces of the spiders do not cover all the genera even if the asymmetric cases are included. 
On the other hand, Kapovich and Millson showed that any smooth manifold can be obtained as a connected component of the configuration space of some planar linkage (\cite {Ka-Mi2}). 
Thus we are lead to: 
\begin{problem}\rm 
Find a family of planar linkages $\{\mathcal{L}_n\}_{n=0,1,2,\cdots}$ such that (a connected component of) the configuration space of $\mathcal{L}_n$ is homeomorphic to $\varSigma_n$. 
\end{problem}

\bigskip
\noindent
{\large{\bf Acknowledgement and comments}}

The author thanks John Crisp and Hiroaki Terao for helpful conversations. 

The author thanks the referee deeply for a lot of invaluable suggestions and for pointing out the difference between the cases $n$ even and $n$ odd, of which the author was not aware in the first draft. 

\medskip
Since there are a great number of references in this topic, it is far from being completed.

Department of Mathematics, Tokyo Metropolitan University, 

1-1 Minami-Ohsawa, Hachiouji-Shi, Tokyo 192-0397, JAPAN. 

E-mail: ohara@comp.metro-u.ac.jp
\end{document}